\documentclass[11pt]{article}
\usepackage[a4paper, margin=2cm, top=1cm, bottom=1.5cm]{geometry}
\usepackage{amsmath, amssymb, amsthm, bm, mathtools}
\usepackage{enumitem}
\usepackage{bbm}
\usepackage[hidelinks]{hyperref}
\usepackage{color}
\usepackage{comment}
\usepackage{booktabs,tabularx}
\usepackage{graphicx}
\usepackage{placeins}
\usepackage{capt-of}
\theoremstyle{plain}

\usepackage{authblk}

\newtheorem{theorem}{Theorem}

\newtheorem{proposition}{Proposition}

\newtheorem{lemma}{Lemma}
\newtheorem*{lemma*}{Lemma}

\newtheorem{definition}{Definition}

\newtheorem{corollary}{Corollary}

\newtheoremstyle{examplestyle}
  {\topsep}{\topsep}
  {\normalfont}{}
  {\itshape}{.}
  {0.5em}{}

\theoremstyle{examplestyle}

\theoremstyle{remark}
\newtheorem{remark}{Remark}

\makeatletter
\let\appendixsubsection\subsection
\renewcommand\subsection{\@startsection{subsection}{2}{\z@}%
  {3.25ex \@plus 1ex \@minus .2ex}%
  {-1em}%
  {\normalfont\large\bfseries}}
\let\originalappendix\appendix
\renewcommand\appendix{%
  \let\subsection\appendixsubsection
  \originalappendix}
\makeatother

\newcommand{\pa}[1]{\left(#1\right)}

\title{Functional Estimation under Local Differential Privacy: General Lower Bounds and Applications to Entropy and Power Sums}
\author{Yann Issartel\thanks{\texttt{yann.issartel@telecom-paris.fr}}\\
LTCI, Télécom Paris, Institut Polytechnique de Paris}
\date{}

\begin{document}

\maketitle

\begin{abstract}
We study the estimation of  functionals of discrete distributions under local differential privacy (LDP). Even for well-studied examples like Shannon entropy and power sums, the existing bounds left open how the minimax risk scales with the alphabet size,  the sample size, and the privacy parameter. To understand these dependences, we establish two general minimax lower bounds for the broader class of additive functionals. The first lower bound links the difficulty of estimating a  functional to how well it can be approximated by polynomials of a given degree. The second complements the first by capturing how sensitive the functional is to small changes in the underlying distribution. To illustrate the effectiveness of these two bounds, we apply them to the estimation of Shannon entropy and power sums. In both cases, they close the existing polynomial gaps and, perhaps surprisingly, show that simple plug-in type estimators already achieve the optimal polynomial rates. 
\end{abstract}


\setcounter{tocdepth}{2}
\makeatletter
\begingroup
\let\contentsaddvspace\addvspace
\renewcommand{\addvspace}[1]{\contentsaddvspace{.25em}}
\@starttoc{toc}
\endgroup
\makeatother

\section{Introduction}

Consider discrete probability distributions on a finite alphabet $[S]:=\{1,\ldots,S\}$, with $S\geq2$, represented by probability vectors in $\mathcal P_S:=\{\boldsymbol{p}\in[0,1]^S:\sum_{j\in[S]}p_j=1\}$. Let $X_1,\ldots,X_n$ be i.i.d. observations with common unknown distribution $\boldsymbol{p}\in\mathcal P_S$. Given a function $\varphi:[0,1]\to\mathbb R$, define the additive functional
\begin{equation}
\label{additive-fucntional}
F_\varphi(\boldsymbol{p}):=\sum_{j\in[S]}\varphi(p_j).
\end{equation}
The classical statistical problem is to estimate $F_\varphi(\boldsymbol{p})$ from the sample $X_{1:n}:=(X_1,\ldots,X_n)$, with estimation error measured by the quadratic risk $\mathbb E_{\boldsymbol{p}}[(\widehat F(X_{1:n})-F_\varphi(\boldsymbol{p}))^2]$.
Two widely studied examples are Shannon entropy and power sums,
\[
F_{\varphi_H}(\boldsymbol{p})=-\sum_{j\in[S]}p_j\log p_j,
\qquad \qquad
F_{\varphi_\gamma}(\boldsymbol{p})=\sum_{j\in[S]}p_j^\gamma,\quad \gamma>0,\]
which correspond to $\varphi_H(x):=-x\log x$, with the convention $0\log0=0$, and to $\varphi_\gamma(x):=x^\gamma$, respectively. Their estimation has been studied extensively, with sharp minimax rates for entropy~\cite{WY} and power sums~\cite{JiaoVenkatHanWeissman2015}, which describe how the estimation error depends on the sample size $n$ and the alphabet size $S$.
The lower bounds for these problems rely on a delicate comparison of mixture distributions, using a subtle link between best polynomial approximation and moment matching. However, this argument requires a specific form of independence that can be lost when the observations are privatized. We will return to this difficulty in Section~\ref{sec:intro-methodology}.

\textit{Functional estimation under local differential privacy (LDP).} The classical estimation problem above assumes that the statistician observes $X_{1:n}$ directly. But when the data are sensitive, sharing them can put individual privacy at risk. Differential privacy provides a mathematical framework for addressing this issue, by protecting individual information~\cite{Dwork2006}. In its local formulation, each individual $i$ releases a randomized version $Z_i$ of their data $X_i$, before any communication to the statistician~\cite{DuchiJordanWainwright2018}. Estimation of $F_{\varphi}(\boldsymbol{p})$ is then based only on the privatized observations $Z_{1:n}:=(Z_1,\ldots,Z_n)$.
This  problem of estimating an additive functional under local differential privacy (LDP) is the focus of the present paper.
It  has already been studied for several functionals, including Shannon entropy~\cite{BH,AZ} and power sums~\cite{ButuceaIssartel2021}. 
Despite these recent developments, the minimax rates remain poorly understood for both functionals. Even if we set aside log factors, their polynomial dependence on $n$ (sample size), $S$ (alphabet size) and $\alpha$ (privacy parameter) is still unresolved. As pointed out in~\cite{ButuceaIssartel2021}, methods that have proved effective in the classical setting do not readily extend to locally private observations. One example is the mixture approach mentioned in the previous paragraph.

\textit{Outline of contributions and paper organization.} Our main contribution is to develop general tools for proving minimax lower bounds for estimating the additive functional $F_\varphi$ in \eqref{additive-fucntional} under LDP. Specifically, we establish two complementary lower bounds, valid (in particular) for any continuous $\varphi$. To explain the challenges behind these results, we first discuss in Section~\ref{sec:intro-methodology} why classical lower-bound methods do not readily extend to LDP. Our two general lower bounds are then presented in Section~\ref{sec2:general-lb-thms}, where we introduce a special family $\mathcal P_{m,b}^{\mathrm{pair}}$ of discrete distributions designed to facilitate lower-bound proofs. This family may also be useful beyond the problems considered here.\\
\indent To illustrate the value of these tools, we apply them to Shannon entropy estimation in Section~\ref{sec:entropy} and to power-sum estimation in Section~\ref{sec-power-sum}. We bring together existing results and new lower and upper bounds to characterize their minimax rates. We also give a detailed interpretation of these rates, highlighting new features of functional estimation under LDP.
For a quick overview of the rates, see Table~\ref{tab:intro-entropy-bounds} for entropy, and Table~\ref{tab:power-sum-bounds-summary} for power sums.\\
\indent
We give the proofs of our two general lower bounds in Sections~\ref{sec:proof-general-lower-bound}--\ref{proof:orientation-lower-bound}, explaining the key ideas and how they address the difficulties discussed in Section~\ref{sec:intro-methodology}. We state key arguments separately and in a general form, so that they can be used beyond the functional estimation problems studied here. The remaining proofs are deferred to the appendices.
 Before giving the formal statistical model in Section~\ref{subsec:ldp-observation-model}, we briefly illustrate the open questions through the problem of entropy estimation.

\textit{A brief preview of entropy estimation $F_{\varphi_H}$.} Without privacy constraints, the minimax quadratic risk for entropy estimation is of order $\frac{S^2}{n^2}+\frac1n$ for $n\gtrsim S$, leaving aside log factors~\cite{WY}. A first question is whether imposing $\alpha$-LDP simply reduces the effective sample size from $n$ to $n_\alpha:=\alpha^2n$, as in many scalar estimation problems~\cite{DuchiJordanWainwright2018}. This would give the rate $\frac{S^2}{n_\alpha^2}+\frac1{n_\alpha}$. However, to our knowledge, the only existing $\alpha$-LDP estimator~\cite{BH} has a risk upper bound of order $\frac{S^2}{n_\alpha}\wedge1$, suggesting a much slower rate. It is therefore unclear whether this gap comes from the perhaps naive plug-in construction behind this estimator, or reflects an additional privacy cost for entropy estimation. 
Such a cost could point to a new phenomenon in scalar estimation problems under LDP.
The existing lower bound $\frac1{n_\alpha} \wedge 1$~\cite{AZ} does not capture the factor $S^2$, leaving this question open. Similar questions remain open for power sums~\cite{ButuceaIssartel2021}, with substantial polynomial gaps in $S$ and $n_\alpha$ between the known lower and upper bounds. The picture is more complex, as the rates depend on the value of the exponent $\gamma$ and can even change from one regime of $n_\alpha$ to another.



\subsection{Statistical model and minimax risk}
\label{subsec:ldp-observation-model}

As explained above, the statistician does not observe $X_{1:n}$ directly, but receives privatized observations $Z_{1:n}$. Conditionally on $X_{1:n}$, these observations are independent and satisfy $Z_i\sim K_i(\cdot\mid X_i)$, where each $K_i$ is a Markov kernel from $[S]$ to a measurable space\footnote{Throughout, the $\sigma$-algebra of each measurable space is omitted from the notation, and products of measurable spaces are equipped with their product $\sigma$-algebras.} $\mathcal Z_i$, subject to the following privacy constraint.
For any real number $\alpha>0$, a kernel $K_i$ is said to be \textit{$\alpha$-locally differentially private} ($\alpha$-LDP) if
\[
K_i(A\mid x)\leq e^\alpha K_i(A\mid x')
\]
for all $x,x'\in[S]$ and all measurable sets $A\subseteq\mathcal Z_i$. We denote by $\mathcal M_\alpha^n$ the class of mechanisms $K_{1:n}:=(K_1,\ldots,K_n)$ such that each kernel $K_i$, $i\in[n]$, satisfies this privacy constraint. Throughout the paper, we focus on the high-privacy regime $\alpha\in(0,1)$, as is common in the statistical LDP literature~\cite{DuchiJordanWainwright2018,ButuceaDuboisSaumard2023,ButuceaRohdeSteinberger2023}. 

\begin{remark}
In the literature, the privacy mechanisms $K_{1:n}\in\mathcal M_\alpha^n$ are often referred to as non-interactive mechanisms: each $K_i$ uses only $X_i$, without taking into account other users' data $X_j$ or privatized outputs $Z_j$, $j\neq i$. The class $\mathcal M_\alpha^n$ is arguably the simplest class of LDP mechanisms and has been widely studied~\cite{DuchiJordanWainwright2018,RohdeSteinberger2020}. We restrict our attention to this class, as estimating nonlinear functionals remains challenging even within $\mathcal M_\alpha^n$. We should still mention that the literature also considers more general mechanisms, allowing interaction between individuals~\cite{ButuceaIssartel2021,ButuceaRohdeSteinberger2023,DuchiRogers2019}.
\end{remark}

For $\boldsymbol{p}\in\mathcal P_S$ and $K_{1:n}\in\mathcal M_\alpha^n$, let $P_{\boldsymbol{p},K_{1:n}}$ denote the joint law of $(X_{1:n},Z_{1:n})$.
For the  functional $F_\varphi$ introduced in \eqref{additive-fucntional}, we define the minimax quadratic risk under $\alpha$-LDP by
\begin{equation}
\label{eq:minimax-risk}
R_{n,\alpha}(F_\varphi)
:=
\inf_{K_{1:n}\in\mathcal M_\alpha^n}
\inf_{\widehat F}
\sup_{\boldsymbol{p}\in\mathcal P_S}
\mathbb E_{\boldsymbol{p},K_{1:n}}
\left[
\left(\widehat F(Z_{1:n})-F_\varphi(\boldsymbol{p})\right)^2
\right],
\end{equation}
where $\mathbb E_{\boldsymbol{p},K_{1:n}}$ denotes expectation under $P_{\boldsymbol{p},K_{1:n}}$. Here, the first infimum is over all $\alpha$-LDP mechanisms $K_{1:n}$, and the second is over all estimators based on the  privatized outputs $Z_{1:n}$.

\section{Challenges in deriving lower bounds under LDP}
\label{sec:intro-methodology}

In this informal section, we motivate our two general lower bounds. We  first explain why classical lower-bound techniques are hard to adapt to LDP. Along the way, we also briefly recall some of these classical techniques. Formal statements come later, in Section~\ref{sec2:general-lb-thms} for our two lower bounds, and in Section~\ref{subsec:classical-tools} for the classical tools we use.

To prove our first general lower bound, we compare two mixture distributions of the privatized sample $Z_{1:n}$. Each mixture is obtained by drawing $\boldsymbol{p}$ from a prior, and then generating $Z_{1:n}$. Following the classical method of two fuzzy hypotheses~\cite{Tsybakov2009}, we choose two priors under which $F_\varphi(\boldsymbol{p})$ concentrates around different values, while the two mixture distributions of $Z_{1:n}$ remain close (in KL divergence or total variation, for example). Section~\ref{subsec:reduction-one-dimensional-mixtures} explains why comparing these mixtures is more difficult under LDP than in classical observation models. Section~\ref{subsec:mixture-approach-limitations} highlights a limitation of this mixture-based approach and motivates a second, complementary approach to proving lower bounds. 

\subsection{First approach: comparing mixtures under LDP}
\label{subsec:reduction-one-dimensional-mixtures}

Recall that, without privacy constraints, the goal is to estimate $F_\varphi(\boldsymbol{p})$ in \eqref{additive-fucntional} directly from the observations $X_1,\ldots,X_n$. To derive a lower bound for this non-private problem, one can compare the two laws of $X_{1:n}$ obtained by drawing $\boldsymbol{p}$ from two different priors. But the difficulty in comparing these mixtures is that the law of each $X_i$ depends on all the probabilities $p_1,\ldots,p_S$, via $\mathbb P(X_i=k)=p_k$ for $k\in[S]$. By contrast, the comparison is easier in a model having a product structure, where each observation depends on just one of these probabilities. Functional estimation problems are therefore often reduced to this type of model, as in the estimation of entropy and power sums~\cite{WY,JiaoVenkatHanWeissman2015}.

A standard way to obtain such a product model is to use the following  Poissonization trick~\cite{MitzenmacherUpfal2017}.  It replaces the fixed sample size $n$ by a random size $N\sim\operatorname{Poi}(n)$, independent of $(X_i)_{i\geq1}$, so that the category counts $N_k:=\sum_{i\in[N]}\mathbf 1_{\{X_i=k\}}$ are independent with laws $N_k\sim\operatorname{Poi}(np_k)$. The joint law is therefore $P_{\boldsymbol{p}}^{N_{1:S}}=\bigotimes_{k=1}^S\operatorname{Poi}(np_k)$, and this product structure provides one independent observation $N_k$ for each parameter $p_k$.

\textit{Mixture factorization in classical models.}
To see how the above structure helps compare mixtures,
consider the canonical  model with this structure: $P_{\boldsymbol\theta}^{Y_{1:S}}=\bigotimes_{j=1}^S P_{\theta_j}^{Y_j}$, where the observations $Y_1,\ldots,Y_S$ are independent and the law of each $Y_j$ depends only on the parameter $\theta_j$. This model covers the above Poissonized model, as well as  Gaussian sequence models $Y_j\overset{\mathrm{ind.}}{\sim}\mathcal N(\theta_j,\sigma^2)$~\cite{CCT}.
For a product prior $\pi=\bigotimes_{j=1}^S\pi_j$ on the parameter vector $\boldsymbol\theta=(\theta_1,\ldots,\theta_S)$,  the induced mixture distribution satisfies
\begin{equation}
\label{eq:classical-mixture-definitions}
M_\pi^{Y_{1:S}}
= \bigotimes_{j=1}^S M_{\pi_j}^{Y_j},
\qquad \text{where }\;
M_\pi^{Y_{1:S}}
:=\int P_{\boldsymbol\theta}^{Y_{1:S}}\,\pi(d\boldsymbol\theta),
\qquad
M_{\pi_j}^{Y_j}
:=\int P_{\theta_j}^{Y_j}\,\pi_j(d\theta_j),
\end{equation}
 meaning that the observations $Y_1,\ldots,Y_S$ remain independent after mixing by $\pi$, with laws $Y_j\sim M_{\pi_j}^{Y_j}$. In particular, for two product priors $\pi=\bigotimes_{j=1}^S\pi_j$ and $\pi'=\bigotimes_{j=1}^S\pi'_j$, the tensorization property of the Kullback--Leibler divergence~\cite{PolyanskiyWu2025} gives
\begin{equation}
\label{eq:intro-mixture-kl-decomposition}
\mathrm{KL}\!\left(M_\pi^{Y_{1:S}}\,\middle\|\,M_{\pi'}^{Y_{1:S}}\right)
\ =\ \sum_{j\in[S]}
\mathrm{KL}\!\left(M_{\pi_j}^{Y_j}\,\middle\|\,M_{\pi'_j}^{Y_j}\right).
\end{equation}
Controlling the divergence between the full mixture distributions therefore reduces to controlling divergences between one-dimensional mixtures. Such one-dimensional comparisons can be handled using classical moment-matching arguments~\cite{WY,JiaoVenkatHanWeissman2015}, which we discuss in Section~\ref{subsec:mixture-approach-limitations}.

\textit{Failure of mixture factorization under LDP.}
For privatized observations, the factorization in~\eqref{eq:classical-mixture-definitions} is generally unavailable. Indeed, although $Z_1,\ldots,Z_n$ are independent for each fixed $\boldsymbol{p}$, with joint law $P_{\boldsymbol{p},K_{1:n}}^{Z_{1:n}}=\bigotimes_{i=1}^n P_{\boldsymbol{p},K_{1:n}}^{Z_i}$, the law of each observation $Z_i$ can depend on all parameters $p_1,\ldots,p_S$, through
\[
P_{\boldsymbol{p},K_{1:n}}^{Z_i}(A_i)
\ = \ 
\sum_{k\in[S]}p_kK_i(A_i\mid k).
\]
Consequently, the observations $Z_{1:n}$ do not generally remain independent after mixing over $\boldsymbol{p}$. More precisely, for a prior $\pi$ on $\mathcal P_S$, the mixture distributions satisfy, in general,
\begin{equation}
\label{eq:privatized-mixture-definitions}
M_{\pi,K_{1:n}}^{Z_{1:n}}
\neq \bigotimes_{i=1}^{n}M_{\pi,K_{1:n}}^{Z_i},
\qquad \text{where }\;
M_{\pi,K_{1:n}}^{Z_{1:n}}
:=\int P_{\boldsymbol p,K_{1:n}}^{Z_{1:n}}\,\pi(d\boldsymbol p),
\qquad
M_{\pi,K_{1:n}}^{Z_i}
:=\int P_{\boldsymbol p,K_{1:n}}^{Z_i}\,\pi(d\boldsymbol p).
\end{equation}
So we can no longer use~\eqref{eq:intro-mixture-kl-decomposition} to reduce the comparison to one-dimensional mixtures.

Overcoming this difficulty is a key contribution of this paper and leads to our first minimax lower bound for additive functionals $F_\varphi$ under LDP, presented in Section~\ref{sec2:general-lb-thms}. To prove this result, we show that, for our choice of priors $\pi$ and $\pi'$, the privatized mixtures $M_{\pi,K_{1:n}}^{Z_{1:n}}$ and $M_{\pi',K_{1:n}}^{Z_{1:n}}$ remain close in total variation, uniformly over all mechanisms $K_{1:n}\in\mathcal M_\alpha^n$. For a precise statement of this total variation bound, readers can turn directly to Theorem~\ref{thm:private-mixture-comparison} in Section~\ref{sec:proof-general-lower-bound}. We devote Section~\ref{sec:proof-private-mixture-comparison} to explaining and proving this bound.

\subsection{A limitation of the mixture approach and a complementary lower bound}
\label{subsec:mixture-approach-limitations}

To see where the mixture approach reaches its limits, we return to the one-dimensional mixtures in~\eqref{eq:classical-mixture-definitions} and recall how they are compared in the classical setting.

\textit{Moment matching and polynomial approximation.}
A classical way to compare one-dimensional mixtures is to use the link between moment matching and polynomial approximation~\cite{WY,JiaoVenkatHanWeissman2015,CCT}. Moment matching consists in choosing priors $\pi_j,\pi'_j$ on $[0,1]$ with the same first $L$ moments, meaning that 
\[\int x^\ell\,\pi_j(dx)=\int x^\ell\,\pi'_j(dx), \qquad \ell\in[L].\]
The number $L$ of matching moments has two competing effects. On the one hand, taking $L$ sufficiently large can make the KL divergences in~\eqref{eq:intro-mixture-kl-decomposition} small, so the mixtures are difficult to distinguish from the observations. On the other hand, the functional separation $\left|\int\varphi\,d\pi_j-\int\varphi\,d\pi'_j\right|$ cannot exceed twice the best approximation error of $\varphi$ by polynomials of degree at most $L$:
\begin{equation}
\label{eq:moment-matching-separation-bound}
\left|\int\varphi\,d\pi_j-\int\varphi\,d\pi'_j\right|
\ \leq\ 
2\inf_{\deg P\leq L}\sup_{x\in[0,1]}|\varphi(x)-P(x)|.
\end{equation}
This directly follows by subtracting any polynomial $P$ of degree at most $L$ and applying the triangle inequality, since moment matching gives $\int P\,d\pi_j=\int P\,d\pi'_j$. The non-trivial fact is that, for a continuous $\varphi$, the  bound~\eqref{eq:moment-matching-separation-bound} is actually attained by some pairs of priors $\pi_j,\pi'_j$.  (See Lemma~\ref{lem:approximation-moment-matching} for a formal statement of this classical result.)
Ultimately, this separation is what gives us a lower bound on the minimax quadratic risk. We therefore want it to be as large as possible, while keeping the two mixtures close (e.g., in KL divergence).

However, when $\varphi$ is well approximated by a polynomial (of degree at most $L$), the bound~\eqref{eq:moment-matching-separation-bound} shows that the separation $\left|\int\varphi\,d\pi_j-\int\varphi\,d\pi'_j\right|$ is small. Thus, this moment-matching method  provides too little separation to capture the difficulty of estimating $F_\varphi$. This motivates a second, complementary lower bound on the minimax quadratic risk.

\textit{A complementary lower bound.} In the usual application of Le Cam's two-point method under LDP, one fixes a pair $(\boldsymbol p,\boldsymbol p')$ with a large functional separation $\lvert F_\varphi(\boldsymbol p)-F_\varphi(\boldsymbol p')\rvert$, then shows that the KL divergence between the two laws of $Z_{1:n}$ is small uniformly over all privacy mechanisms $K_{1:n}\in\mathcal M_\alpha^n$ (see, e.g.,~\cite{DuchiJordanWainwright2018,AZ,IssartelRoueff}). In the present work, however, we do not require a single pair $(\boldsymbol p,\boldsymbol p')$ to remain difficult to distinguish under all privacy mechanisms. Instead, for each mechanism $K_{1:n}\in\mathcal M_\alpha^n$, we choose a pair $(\boldsymbol p,\boldsymbol p')_{K_{1:n}}$ whose privatized laws are close in KL divergence. This flexibility lets us keep the KL divergence small without sacrificing functional separation. We present this second lower bound in Section~\ref{sec2:general-lb-thms}.

\section{Two general lower bounds for additive functionals}
\label{sec2:general-lb-thms}

We now state the two minimax lower bounds motivated in Section~\ref{sec:intro-methodology}. We first introduce a common family of probability vectors (Section~\ref{subsec:common-construction}), then present the mixture-based bound (Section~\ref{subsec:polynomial-lower-bound}) and the complementary lower bound (Section~\ref{subsec:two-point-lower-bound}).

\subsection{A common construction based on pairwise perturbations}
\label{subsec:common-construction}
Both lower bounds use the following class of probability vectors $\boldsymbol p$, based on $m$ pairs of probabilities with a common  scale $b$.

\begin{definition}
\label{def:pairwise-perturbations}
For an integer $1\leq m\leq\lfloor S/2\rfloor$
and a real number $b>0$ such that $2mb\leq1$, define
\[
\mathcal P_{m,b}^{\mathrm{pair}}
:=
\Bigl\{
\boldsymbol{p}\in\mathcal P_S\ :\ 
p_{2j-1}+p_{2j}=2b,
\ \text{ for every }j\in[m]
\Bigr\}.
\]
For $\boldsymbol p\in\mathcal P_{m,b}^{\mathrm{pair}}$ and $\boldsymbol t=(t_1,\ldots,t_m)\in[-1,1]^m$, define $\boldsymbol p(\boldsymbol t)$ by
\begin{equation}
\label{eq:paired-probabilities}
p_{2j-1}(\boldsymbol t)
\ = \ b(1+t_j),
\qquad  \qquad
p_{2j}(\boldsymbol t)
\ = \ b(1-t_j),
\quad \qquad  j\in[m],
\end{equation}
and $p_k(\boldsymbol t)=p_k$ for $2m<k\leq S$. The subfamily $\{\boldsymbol p(\boldsymbol t)\}_{\boldsymbol t\in[-1,1]^m}$ contains all vectors in $\mathcal P_{m,b}^{\mathrm{pair}}$ that agree with $\boldsymbol p$ on the last $S-2m$ entries.
\end{definition}

The vector $\boldsymbol p(\boldsymbol t)$ can be viewed as a pairwise perturbation of $\boldsymbol p$, with $m$ pairs $(p_{2j-1},p_{2j})$ parametrized by~$\boldsymbol t$. This construction allows the parameters $t_1,\ldots,t_m$ to vary independently while preserving the constraint $\sum_{k=1}^S p_k(\boldsymbol t)=1$. Each pair has total probability mass $2b$, and $t_j$ controls its imbalance: $t_j=0$ gives equal probabilities, whereas $|t_j|=1$ assigns the entire mass to one category. Changing the sign of $t_j$ swaps the two probabilities. One benefit of this construction, for example in mixture-based lower bounds, is that we can put a product prior $\mu^{\otimes m}$ on $\boldsymbol t$: then, drawing $\boldsymbol t\sim\mu^{\otimes m}$ always produces a probability vector $\boldsymbol p(\boldsymbol t)\in\mathcal P_S$.

Recall that $\varphi:[0,1]\to\mathbb R$ defines the additive functional $F_\varphi$, and that the minimax risk $R_{n,\alpha}(F_\varphi)$ over $\mathcal P_S$ is defined in~\eqref{eq:minimax-risk}. The lower bounds we prove, apply more precisely to the risk restricted to the subclass $\mathcal P_{m,b}^{\mathrm{pair}}$,  which is 
\[
R_{n,\alpha}(F_\varphi;\mathcal P_{m,b}^{\mathrm{pair}})
\ :=\ 
\inf_{K_{1:n}\in\mathcal M_\alpha^n}
\inf_{\widehat F}
\sup_{\boldsymbol p\in\mathcal P_{m,b}^{\mathrm{pair}}}
\mathbb E_{\boldsymbol p,K_{1:n}}
\left[\bigl(\widehat F(Z_{1:n})-F_\varphi(\boldsymbol p)\bigr)^2\right],
\]
where the infima are over the same privacy mechanisms and estimators as in the unrestricted risk $R_{n,\alpha}(F_\varphi)$. Since $\mathcal P_{m,b}^{\mathrm{pair}}\subseteq\mathcal P_S$, we  have $R_{n,\alpha}(F_\varphi)\geq R_{n,\alpha}(F_\varphi;\mathcal P_{m,b}^{\mathrm{pair}})$.

For a probability vector $\boldsymbol p(\boldsymbol t)$ of the form~\eqref{eq:paired-probabilities}, each pair of probabilities contributes $\varphi[b(1+t_j)]+\varphi[b(1-t_j)]$ to the functional $F_\varphi$. This contribution is described by the pair profile $\varphi_b$, defined below.

\begin{definition}[Pair profile]
\label{def:pair-profile}
For $b\in(0,1/2]$, define the
\emph{pair profile} of $\varphi$ by
\begin{equation}
\label{eq:paired-functional}
\varphi_b(t)
\ := \
\varphi\bigl[b(1+t)\bigr]+\varphi\bigl[b(1-t)\bigr],
\qquad t\in[0,1].
\end{equation}
\end{definition}
The pair profile $\varphi_b$ entirely determines how the functional $F_\varphi(\boldsymbol p(\boldsymbol t))$ varies with $\boldsymbol t$. Indeed, it is not difficult to see that changing a coordinate $t_j$ changes  $F_\varphi(\boldsymbol p(\boldsymbol t))$ only through $\varphi_b(|t_j|)$. As we will see in the next two subsections, this pair profile plays a central role in both lower bounds.


\subsection{A lower bound via polynomial approximation}
\label{subsec:polynomial-lower-bound}

We use two quantities: the error $E_L(f)$ of approximating a function $f$ by polynomials of degree at most $L$, and its range width $V(f)$, which are formally defined by 
\[
E_L(f)
:=
\inf_{\deg P\leq L}\sup_{x\in[0,1]}|f(x)-P(x)|,
\qquad\qquad
V(f)
:=
\max_{x\in[0,1]}f(x)-\min_{x\in[0,1]}f(x),
\]
for any continuous function $f:[0,1]\to\mathbb R$ and any integer $L\geq0$.

Our first lower bound relates the difficulty of estimating $F_\varphi$ to the best polynomial approximation error of its pair profile $\varphi_b$ (from Definition~\ref{def:pair-profile}).

\begin{theorem}
\label{thm:general-additive-lower-bound}
\label{thm:general-additive-lower-bound-pairwise}
Let $\varphi:[0,1]\to\mathbb R$ be a continuous function.
Let $S\geq3$, $n\geq1$, and $\alpha\in(0,1)$.
Let $m,b$ define a nonempty class $\mathcal P_{m,b}^{\mathrm{pair}}$ as in Definition~\ref{def:pairwise-perturbations}.
If an integer $L\geq0$ satisfies
\begin{equation}
\label{eq:general-lb-parameters}
L\ \geq\ \left\lceil\log_2(2m)\right\rceil,
\qquad\qquad
b\ \leq\ C_0\sqrt{\frac{L+1}{n\alpha^2}},
\end{equation}
with $C_0=(2\sqrt e\,(e-1))^{-1}$, then we have
\begin{equation}
\label{eq:general-additive-lower-bound}
R_{n,\alpha}(F_\varphi;\mathcal P_{m,b}^{\mathrm{pair}})
\ \geq\ 
\frac{m^2}{4}\,
\left[
 c_0 E_L^2(\varphi_b)
 -\frac{V^2(\varphi_b)}{m}
\right]_+,
\end{equation}
where $c_0=3/8$ and $[x]_+:=\max\{x,0\}$.
\end{theorem}

The bound above has a bias--variance type form: a squared approximation error $E_L^2(\varphi_b)$ minus a variance-type term $V^2(\varphi_b)/m$. Theorem~\ref{thm:general-additive-lower-bound} thus identifies a source of estimation difficulty: when the pair profile $\varphi_b$ is poorly approximated by polynomials of degree at most $L$, and its range width $V(\varphi_b)$ remains moderate, the minimax risk is large.

In applications to specific functionals, such as entropy or power sums, we optimize the lower bound~\eqref{eq:general-additive-lower-bound} over the polynomial degree $L$ and the parameters of $\mathcal P_{m,b}^{\mathrm{pair}}$: the number $m$ of pairwise perturbations and the probability scale $b$. Looking at~\eqref{eq:general-additive-lower-bound} suggests a simple way to maximize the lower bound: take $m$ large, of order $S$ if possible, while keeping $L$ small so that the polynomial approximation error $E_L(\varphi_b)$ remains large. However, condition~\eqref{eq:general-lb-parameters} restricts our choice of $(L,m,b)$, as we explain in the Discussion at the end of this subsection.

To our knowledge, Theorem~\ref{thm:general-additive-lower-bound} is the first minimax lower bound under LDP obtained through the mixture approach described in Section~\ref{sec:intro-methodology}, a method that has proved effective for non-private functional estimation~\cite{WY,CCT}. To use this approach under LDP, we must show that the mixtures remain hard to distinguish after privatization. As explained in Section~\ref{subsec:reduction-one-dimensional-mixtures}, the classical argument does not carry over to LDP, so we establish a new bound on the TV distance between the privatized mixtures, valid uniformly over all mechanisms in $\mathcal M_\alpha^n$. This bound, presented in Section~\ref{subsec:controlling-privatized-mixtures}, is the main new ingredient in the proof of Theorem~\ref{thm:general-additive-lower-bound}. The proof also uses two classical tools, the method of two fuzzy hypotheses and the link between moment matching and polynomial approximation, both recalled in Section~\ref{subsec:classical-tools}. We give the proof of Theorem~\ref{thm:general-additive-lower-bound} in Section~\ref{proof_thm1}.

\paragraph{Discussion.}
We discuss the origin of the terms in the lower bound and the effects of the parameters.
\begin{itemize}[label=--, leftmargin=0.6em, topsep=0pt, itemsep=3pt, parsep=0pt, partopsep=0pt]
\item \textit{The polynomial approximation term.} The term $E_L(\varphi_b)$ comes from the classical connection between moment matching and polynomial approximation recalled in~\eqref{eq:moment-matching-separation-bound}, with a formal statement in Section~\ref{subsec:classical-tools}. Applied to $\varphi_b$, this result gives priors $\pi_j,\pi'_j$ on $[0,1]$, with the same first $L$ moments, such that $\int\varphi_b\,d\pi_j-\int\varphi_b\,d\pi'_j=2E_L(\varphi_b)$. We take the priors $\pi:=\bigotimes_{j=1}^m\pi_j$ and $\pi':=\bigotimes_{j=1}^m\pi'_j$ for $\boldsymbol p(\boldsymbol t)$ by drawing the coordinates $t_j$ independently from $\pi_j$ and $\pi'_j$, respectively. The contributions of the $m$ pairwise perturbations in $\mathcal P_{m,b}^{\mathrm{pair}}$ add up, giving a separation $2mE_L(\varphi_b)$ between the expectations of $F_\varphi(\boldsymbol p(\boldsymbol t))$ under $\pi$ and $\pi'$. Squaring this separation explains the term $m^2E_L^2(\varphi_b)$ in the lower bound~\eqref{eq:general-additive-lower-bound}.
\item \textit{The variance-type term.} Under the priors $\pi$ and $\pi'$ on $\boldsymbol t$, the functional $F_\varphi(\boldsymbol p(\boldsymbol t))$ fluctuates around its expectation. Its random part is a sum of $m$ independent contributions, each with variance controlled by the squared range width $V^2(\varphi_b)$. The total variance is therefore controlled by $mV^2(\varphi_b)$. Factoring out $m^2$ in the lower bound~\eqref{eq:general-additive-lower-bound} gives the variance-type term $V^2(\varphi_b)/m$.
\item \textit{Polynomial degree $L$.} As in~\eqref{eq:privatized-mixture-definitions}, let $M_{\pi,K_{1:n}}^{Z_{1:n}}$ and $M_{\pi',K_{1:n}}^{Z_{1:n}}$ be the privatized mixtures corresponding to the priors $\pi$ and $\pi'$, whose moments match up to order $L$. Increasing the polynomial degree $L$ therefore imposes more matching moments on these priors, which helps keep the mixtures close. However, increasing $L$ may also reduce the polynomial approximation error $E_L(\varphi_b)$ in the lower bound~\eqref{eq:general-additive-lower-bound}.
\item \textit{Probability scale $b$ and privacy level $\alpha$.} The typical noise scale for estimating one coordinate $p_k$ of the probability vector $\boldsymbol p$ under LDP is $1/\sqrt{n\alpha^2}$. By comparison, the condition $b\lesssim\sqrt{(L+1)/(n\alpha^2)}$ in~\eqref{eq:general-lb-parameters} allows the probability scale $b$ to be larger by a factor of order $\sqrt{L+1}$. This is made possible by matching the first $L$ moments of the priors $\pi$ and $\pi'$, which keeps the corresponding privatized mixtures close even at this larger scale $b$. This gives us extra room in choosing $b$ to optimize the lower bound~\eqref{eq:general-additive-lower-bound}. 
\item \textit{Number of pairwise perturbations $m$.} Increasing $m$ introduces more perturbed coordinates in the model $\mathcal P_{m,b}^{\mathrm{pair}}$, whose contributions to the distance between the privatized mixtures $M_{\pi,K_{1:n}}^{Z_{1:n}}$ and $M_{\pi',K_{1:n}}^{Z_{1:n}}$ accumulate. To keep these mixtures close, the priors $\pi$ and $\pi'$ must therefore match more moments, as expressed by the condition $L\geq\lceil\log_2(2m)\rceil$ in~\eqref{eq:general-lb-parameters}.
\end{itemize}

\subsection{A complementary lower bound via local variation}
\label{subsec:two-point-lower-bound}

As discussed in Section~\ref{subsec:mixture-approach-limitations}, moment matching can give a functional separation that is too small when the functional is well approximated by a polynomial. In Theorem~\ref{thm:general-additive-lower-bound}, this limitation appears through the term $E_L(\varphi_b)$, which can even vanish. For instance, when $\varphi(x)=x^2$, we obtain the quadratic functional $F_\varphi(\boldsymbol p)=\sum_{j\in[S]}p_j^2$, whose pair profile is
\begin{equation}
\label{eq:quadratic-paired-functional}
\varphi_b(t)
\ = \
\bigl[b(1+t)\bigr]^2+\bigl[b(1-t)\bigr]^2
\ = \
2b^2(1+t^2),
\end{equation}
so $E_L(\varphi_b)=0$ for any $L\geq2$, and Theorem~\ref{thm:general-additive-lower-bound} gives only the trivial lower bound zero, thus failing to capture the estimation difficulty of the functional.

This limitation motivates our second, complementary lower bound, which depends on the local variations of the pair profile $\varphi_b$ from Definition~\ref{def:pair-profile}.

\begin{theorem}
\label{thm:orientation-lower-bound}
\label{thm:orientation-lower-bound-pairwise}
Let $\varphi:[0,1]\to\mathbb R$ be a function.
Let $S\geq3$, $n\geq1$, and $\alpha\in(0,1)$.
Let $m,b$ define a nonempty class $\mathcal P_{m,b}^{\mathrm{pair}}$ as in Definition~\ref{def:pairwise-perturbations}.
Define $\delta_{m,b}:=1\wedge\bigl(c_1mb^2n\alpha^2\bigr)^{-1/2}$, where $c_1=16e(e-1)^2$. Then, we have
\begin{equation}
\label{eq:orientation-lower-bound}
R_{n,\alpha}(F_\varphi;\mathcal P_{m,b}^{\mathrm{pair}})
\ \geq \
\frac{m^2}{16}\,
\sup_{\substack{t,t'\in[0,1]\\ |t-t'|\leq\delta_{m,b}}}
\left|\varphi_b(t)-\varphi_b(t')\right|^2.
\end{equation}
\end{theorem}

The lower bound~\eqref{eq:orientation-lower-bound} depends on how much the pair profile $\varphi_b$ varies between nearby parameter values $t,t'\in[0,1]$, at distance at most $|t-t'|\leq\delta_{m,b}$. The larger these local variations of $\varphi_b$, the larger the lower bound. The parameters in $\mathcal P_{m,b}^{\mathrm{pair}}$ ---the number $m$ of pairwise perturbations and the probability scale $b$--- should be chosen to maximize this lower bound. A large $m$, of order $S$ if possible, increases the factor $m^2$ in~\eqref{eq:orientation-lower-bound}, while a larger $b$ can enlarge the variations of the pair profile $\varphi_b$. However, increasing either $m$ or $b$ reduces the distance $|t-t'|\leq\delta_{m,b}$ allowed in the supremum. The choice of $(m,b)$ therefore involves a trade-off through its effects on the three quantities $m^2$, $\varphi_b$, and $\delta_{m,b}$ appearing in the lower bound~\eqref{eq:orientation-lower-bound}. In our applications to entropy and power sums, we use this trade-off to choose $m$ and $b$ and obtain lower bounds for each functional.

Privacy plays a key role in the distance $\delta_{m,b}$ allowed in the supremum. Under $\alpha$-LDP, this distance is $\delta_{m,b}\asymp 1\wedge(b\sqrt{mn\alpha^2})^{-1}$, whereas in the non-private setting it would typically be of order $1\wedge(nmb)^{-1/2}$. This corresponds to an effective sample size reduction from $n$ to $\alpha^2bn$, which, since $b<1$, is stronger than the usual reduction to $\alpha^2n$ in many scalar estimation problems. For example, when $b\asymp1/m$, the effective sample size is of order $\alpha^2n/m$, giving an extra reduction factor of $1/m$. This larger range $\delta_{m,b}$ is what allows us to obtain an interesting larger lower bound on the risk $R_{n,\alpha}$. In Section~\ref{subsec:orientation-kl-bound}, we explain this reduction from $n$ to $\alpha^2bn$ by comparing the KL divergences with and without privacy.

\paragraph{What the second lower bound adds.}
Unlike Theorem~\ref{thm:general-additive-lower-bound}, the lower bound~\eqref{eq:orientation-lower-bound} can be positive even when $\varphi$ is a polynomial of small degree. An example is the quadratic functional introduced just before Theorem~\ref{thm:orientation-lower-bound}. For this functional, Theorem~\ref{thm:orientation-lower-bound} gives a lower bound of order $1\wedge(n\alpha^2)^{-1}$. To see this, take $(m,b)=(1,1/2)$ and $(t,t')=(1,1-\delta_{1,1/2})$, so that, using $\varphi_{1/2}(t)=(1+t^2)/2$ from~\eqref{eq:quadratic-paired-functional},
\[
\varphi_{1/2}(1)-\varphi_{1/2}(1-\delta_{1,1/2})
\ = \
\frac{(1+1^2)-\bigl[1+(1-\delta_{1,1/2})^2\bigr]}{2}
\ = \
\delta_{1,1/2}-\frac{\delta_{1,1/2}^2}{2}
\ \geq\ \frac{\delta_{1,1/2}}2;
\]
then Theorem~\ref{thm:orientation-lower-bound}  gives $R_{n,\alpha}(F_\varphi;\mathcal P_{1,1/2}^{\mathrm{pair}})\geq\delta_{1,1/2}^2/64\geq c(1\wedge(n\alpha^2)^{-1})$, with $c>0$ a universal constant.

Interestingly, this  lower bound of order $1\wedge(n\alpha^2)^{-1}$ on the global minimax risk $R_{n,\alpha}(F_\varphi)$, matches known upper bounds for this problem~\cite{ButuceaIssartel2021} (up to log factors in the non-interactive setting, and up to constants in the sequentially interactive setting).
Theorem~\ref{thm:orientation-lower-bound} can therefore capture the estimation difficulty of a problem where Theorem~\ref{thm:general-additive-lower-bound} fails, illustrating the complementarity of these two lower-bound tools. To explore their power beyond the above  example, we apply them to entropy (Section~\ref{sec:entropy}) and power sums (Section~\ref{sec-power-sum}) ---two widely studied functionals.


\paragraph{Proof idea for Theorem~\ref{thm:orientation-lower-bound}.}
The proof relies on the classical Le Cam two-point method. Usually, this approach simply consists in fixing two observation laws $P$ and $P'$ in advance, with separated values of the functional, and proving that their divergence $\operatorname{KL}(P\|P')$ is small, making them hard to distinguish. In our proof, instead of fixing two distributions in advance, we work with two families of distributions in $\mathcal P_{m,b}^{\mathrm{pair}}$, denoted here by $\mathcal F_1$ and $\mathcal F_2$. For any fixed privacy mechanism $K_{1:n}$, we aim to find a distribution $\boldsymbol p\in\mathcal F_1$ and a distribution $\boldsymbol p'\in\mathcal F_2$ for which $\operatorname{KL}\bigl(P_{\boldsymbol p,K_{1:n}}^{Z_{1:n}}\|P_{\boldsymbol p',K_{1:n}}^{Z_{1:n}}\bigr)$ is small. To show that such a pair exists, we use the probabilistic method~\cite{AlonSpencer2016}, which reduces this combinatorial problem to the simpler task of bounding the average $\mathbb E\bigl[\operatorname{KL}\bigl(P_{\boldsymbol p,K_{1:n}}^{Z_{1:n}}\|P_{\boldsymbol p',K_{1:n}}^{Z_{1:n}}\bigr)\bigr]$ when $(\boldsymbol p,\boldsymbol p')$ is drawn at random from $\mathcal F_1\times\mathcal F_2$. Note that the selected pair $(\boldsymbol p,\boldsymbol p')$ may depend on the privacy mechanism $K_{1:n}$. This contrasts with the usual approach to privacy lower bounds, where a pair $(\boldsymbol p,\boldsymbol p')$ is fixed in advance and the KL divergence is then bounded uniformly over all privacy mechanisms $K_{1:n}$.

It remains to ensure that the functional separation $\lvert F_\varphi(\boldsymbol p)-F_\varphi(\boldsymbol p')\rvert$ is sufficiently large. This could be delicate, since the argument above only guarantees the existence of a pair $(\boldsymbol p,\boldsymbol p')$ with small KL divergence, and this pair may depend on the privacy mechanism $K_{1:n}$. Fortunately, the construction of $\mathcal P_{m,b}^{\mathrm{pair}}$ in Definition~\ref{def:pairwise-perturbations} makes this step simple. The functional $F_\varphi$ is constant within each family $\mathcal F_1$ and $\mathcal F_2$, so every pair in $\mathcal F_1\times\mathcal F_2$ has the same functional separation. We can thus choose a pair $(\boldsymbol p,\boldsymbol p')$ to make the KL divergence small without sacrificing separation. See  Section~\ref{proof:orientation-lower-bound} for the detailed proof.





\section{Shannon entropy estimation under LDP}
\label{sec:entropy}

Recall that the Shannon entropy is $F_{\varphi_H}(\boldsymbol p)=\sum_{k\in[S]}\varphi_H(p_k)=-\sum_{k\in[S]}p_k\log p_k$. We study how its minimax quadratic risk $R_{n,\alpha}(F_{\varphi_H})$ depends on the sample size $n$, the privacy parameter $\alpha$, and the alphabet size $S$. We first present a simple estimator and its risk bound (Section~\ref{sec:entropy-plugin}). We then determine the minimax rate up to log factors (Section~\ref{section-entrpy-polynomial}), followed by the exact asymptotic rate (Section~\ref{sec:entropy-complementary}). We discuss these rates (Section~\ref{sec:entropy-discussion}) and conclude with open questions (Section~\ref{sec:entropy-open-questions}). For a summary of the bounds presented in this section, see Table~\ref{tab:intro-entropy-bounds}.


\subsection{A plug-in estimator.}
\label{sec:entropy-plugin}
A natural approach to entropy estimation is to first
estimate the probabilities $p_k$ and then plug these
estimates into the entropy functional,
$-\sum_{k\in[S]} \widehat p_k\log\widehat p_k$.
This approach was  proposed by
Bravo-Hermsdorff et al.~\cite{BH}\footnote{Theorem~5.1 of~\cite{BH} gives a sample-complexity bound for their estimator. However, its proof omits log factors, when converting the error in estimating the distribution into an error in estimating its entropy. These log factors are necessary, as we will see in Section~\ref{sec:entropy-complementary}.}, using private
distribution estimators designed to reduce
 communication costs.
Since we impose no communication constraints, we use a standard Laplace mechanism to implement this plug-in approach, as described below.

\noindent\textit{Estimator.} For each observation $X_i$, release
$Z_i=(Z_{ik})_{k\in[S]}$, where
\[
Z_{ik}
=
\mathbf{1}_{\{X_i=k\}}+\frac{2}{\alpha}W_{ik},
\qquad i\in[n],\quad k\in[S],
\]
and the $W_{ik}$ are i.i.d.\ centered standard Laplace
random variables, independent of the $X_i$'s.
Denote by $K_i^{\mathrm{Lap},\alpha}$ the corresponding privatization kernel, mapping $X_i$ to $Z_i$, and write $K_{1:n}^{\mathrm{Lap},\alpha}=(K_1^{\mathrm{Lap},\alpha},\ldots,K_n^{\mathrm{Lap},\alpha})$ for the mechanism obtained by applying these $n$ kernels independently to $X_1,\ldots,X_n$. It is straightforward to check that $K_i^{\mathrm{Lap},\alpha}\in\mathcal M_\alpha$ for every $i\in[n]$, and hence $K_{1:n}^{\mathrm{Lap},\alpha}\in\mathcal M_\alpha^n$.
Then, we estimate each probability $p_k$ by averaging
the privatized observations and projecting onto $[0,1]$,
then plug these estimates into the entropy functional:
\[
\widehat p_k
=
T_{[0,1]}\left(\frac1n\sum_{i\in[n]} Z_{ik}\right),
\, \,  k\in[S],
\qquad \quad 
\widehat F_{\varphi_H}
=
T_{[0,\log S]}
\left(-\sum_{k\in[S]}\widehat p_k\log\widehat p_k\right),
\]
where $T_{[a,b]}(x)=(x\wedge b)\vee a$ denotes the projection onto $[a,b]$.
The next proposition provides an upper bound
on the quadratic risk of the estimator $\widehat F_{\varphi_H}$.
The proof is given in Appendix~\ref{app:entropy-plug-in}.

\begin{proposition}
\label{prop:entropy-laplace-upper}
There exists a universal constant $C>0$ such that, for every $n\geq1$, $S\geq3$, and $\alpha\in(0,1)$, under the privatization mechanism $K_{1:n}^{\mathrm{Lap},\alpha}$, the estimator $\widehat F_{\varphi_H}$ satisfies
\begin{equation}
\label{eq:entropy-upper-comparison}
\sup_{\boldsymbol{p}\in\mathcal P_S}
\mathbb E_{\boldsymbol{p},K_{1:n}^{\mathrm{Lap},\alpha}}
\bigl[(\widehat F_{\varphi_H}-F_{\varphi_H}(\boldsymbol{p}))^2\bigr]
\ \leq \
C\,\pa{
(\log S)^2
\wedge
\frac{S^2\log^2(e+n\alpha^2)}{n\alpha^2}}.
\end{equation}
\end{proposition}

In what follows, we use this upper bound as a reference, and study its minimax optimality.

\subsection{Non-asymptotic rate and polynomial dependence on the alphabet size.}
\label{section-entrpy-polynomial} The quadratic dependence on $S$ in the upper bound~\eqref{eq:entropy-upper-comparison} raises a natural question: is this dependence unavoidable, or is it the price of using a perhaps naive plug-in approach? To our knowledge, the best available lower bound is due to Asoodeh and Zhang~\cite[Corollary~2]{AZ}, who show that
\begin{equation}
\label{eq:entropy-previous-lower}
R_{n,\alpha}(F_{\varphi_H})
\ \gtrsim \
(\log S)^2
\left(1\wedge\frac{1}{n\alpha^2}\right),
\end{equation}
for $S\geq3$, $n\geq1$, and $\alpha\in(0,1)$. This lower bound does not capture the factor $S^2$ in the upper bound~\eqref{eq:entropy-upper-comparison}, leaving the question open.
Our general lower bound in Theorem~\ref{thm:general-additive-lower-bound}, applied to $\varphi_H(t)=-t\log t$, gives the following result. Its proof is deferred to Appendix~\ref{section:proof_entropy_lb_appli_Thm1}.

\begin{corollary}
\label{thm:entropy-lower}
There exist universal constants $c>0$, $S_0\geq4$
such that, for all $n\geq1$, $S\geq S_0$, $\alpha\in(0,1)$, we have
\[
R_{n,\alpha}(F_{\varphi_H})
\ \geq \ 
c\left(
\frac{1}{(\log S)^4}
\wedge
\frac{S^2}{n\alpha^2(\log S)^3}
\right).
\]
\end{corollary}

Comparing this  bound with 
\eqref{eq:entropy-upper-comparison} characterizes
the minimax risk up to log factors:
\begin{equation}
\label{eq:entropy-minimax-up-to-logs}
R_{n,\alpha}(F_{\varphi_H})
\ \  \widetilde{\asymp} \ \
1 \wedge \frac{S^2}{n\alpha^2},
\end{equation}
where $\widetilde{\asymp}$ suppresses
log factors.
This closes the polynomial gap in $S$ left by~\eqref{eq:entropy-previous-lower}, and also shows that the plug-in estimator has the optimal polynomial dependence on $(S,n,\alpha)$, among all LDP mechanisms and estimators. This is somewhat surprising, since adding Laplace noise to each observation strongly perturbs the empirical frequencies $\frac1n\sum_{i=1}^n\mathbf1_{\{X_i=k\}}$ of the categories $k\in[S]$, and may seem too naive to achieve good rates. It is therefore natural that other work has tried to improve on the simple plug-in approach (e.g.~\cite{ButuceaIssartel2021} for power sums), but Corollary~\ref{thm:entropy-lower} rules out any improvement in the polynomial dependence on $(S,n,\alpha)$ for entropy estimation. We discuss the non-asymptotic rate~\eqref{eq:entropy-minimax-up-to-logs} at the end of this section.

\subsection{Asymptotic  rate.}
\label{sec:entropy-complementary}

We now consider the asymptotic regime $n\alpha^2\to\infty$, with $S$ fixed. (The privacy parameter $\alpha$ may vary with $n$.) The plug-in upper bound~\eqref{eq:entropy-upper-comparison} is then of order $\log^2(n\alpha^2)/(n\alpha^2)$, while the existing lower bound~\eqref{eq:entropy-previous-lower} and our lower bound in Corollary~\ref{thm:entropy-lower} are both of order $1/(n\alpha^2)$. Is this factor $\log^2(n\alpha^2)$ necessary? Our general lower bound in Theorem~\ref{thm:orientation-lower-bound}, applied to $\varphi_H(t)=-t\log t$, gives the following result. Its proof is deferred to Appendix~\ref{appendix:cmplement_lower_bound_entropy}.

\begin{corollary}
\label{cor:entropy-complementary-lower}
For any  $S\geq3$, we have $R_{n,\alpha}(F_{\varphi_H})\gtrsim_S\log^2(n\alpha^2)/(n\alpha^2)$, as $n\alpha^2\to\infty$ with $\alpha\in(0,1)$. Here, $\gtrsim_S$ hides a positive constant depending only on $S$.
\end{corollary}

Together with Proposition~\ref{prop:entropy-laplace-upper}, this gives the asymptotic minimax rate
\begin{equation}
\label{eq:entropy-fixed-alphabet-minimax}
R_{n,\alpha}(F_{\varphi_H})
\asymp_S
\frac{\log^2(n\alpha^2)}{n\alpha^2}
\qquad\text{as }n\alpha^2\to\infty,
\end{equation}
where $\asymp_S$ hides constants that may depend on $S$. Thus, the Laplace plug-in estimator achieves the optimal dependence on $n$ and $\alpha$ for a fixed alphabet. Surprisingly, the factor $\log^2(n\alpha^2)$ cannot be removed by changing the estimator or the privacy mechanism. We discuss this minimax rate below.


\subsection{Discussion of the rates.}
\label{sec:entropy-discussion}
 In many scalar estimation problems under LDP, privacy simply reduces the effective sample size from $n$ to $n\alpha^2$. For entropy estimation, however, its effect is more complex, as explained below.

 \textit{Comparing the non-asymptotic rate~\eqref{eq:entropy-minimax-up-to-logs} with the minimax classical rate.}
 Without privacy constraints, the results of Wu and Yang~\cite{WY} give the minimax quadratic risk $1\wedge\left(\frac{S^2}{n^2}+\frac1n\right)$, up to log factors. Comparing this with our private rate~\eqref{eq:entropy-minimax-up-to-logs} reveals two effects of privacy. (i) For small sample sizes $n$, both rates are constant, corresponding to the term $1$. The sample size required to leave this constant-risk regime increases from order $S$ without privacy to order $S^2/\alpha^2$ with privacy. (ii) For large sample sizes, when $n\alpha^2\geq S^2$, the non-private rate is $\frac1n$, whereas the private rate is $\frac{S^2}{n\alpha^2}$. Thus, the private rate preserves a dependence on the alphabet size even when the non-private rate no longer has one. These comparisons are understood up to log factors.

\textit{Interpreting~\eqref{eq:entropy-minimax-up-to-logs} through small probabilities.} The difference between the private and non-private rates can be understood through the difficulty of estimating small probabilities. Without privacy, the empirical frequency of a category $k\in[S]$, with probability $p_k$, has standard deviation of order $\sqrt{p_k/n}$ for small $p_k$. With the Laplace mechanism, the added noise contributes an error of order $1/\sqrt{n\alpha^2}$, which does not decrease with $p_k$. When the estimation error is as large as $p_k$, estimating $p_k$ becomes difficult; this happens at $p_k\asymp1/n$ without privacy, and at $p_k\asymp1/\sqrt{n\alpha^2}$ with the Laplace mechanism. We simply write $\varepsilon$ for these two critical scales. When $S\varepsilon\leq1$, a distribution can have a number of such probabilities proportional to $S$, so their aggregate entropy contribution is of order $S\varepsilon$, suggesting a squared error of order $S^2\varepsilon^2$. This recovers the terms $\frac{S^2}{n^2}$ without privacy and $\frac{S^2}{n\alpha^2}$ with privacy in the above rates. When $S\varepsilon>1$, a distribution can have $S$ probabilities of order $1/S$, all below the critical scale, so their aggregate entropy contribution is of order $1$, suggesting the constant term in the above rates. Together, these two cases recover our private rate~\eqref{eq:entropy-minimax-up-to-logs}, and the contribution $1\wedge\frac{S^2}{n^2}$ to the non-private rate of~\cite{WY}. Thus, the difference between the rates can be understood simply: privacy noise makes some small probabilities hard to estimate, although they could be estimated accurately without privacy. This simple, high-level explanation should not suggest that the proof is straightforward: proving the lower bounds requires much more delicate arguments, as detailed in Sections~\ref{sec:proof-general-lower-bound}--\ref{proof:orientation-lower-bound}. 

\textit{The asymptotic rate~\eqref{eq:entropy-fixed-alphabet-minimax} for fixed $S$.} The same reasoning also explains the log factor in our asymptotic minimax rate~\eqref{eq:entropy-fixed-alphabet-minimax}. A probability $\varepsilon$ contributes $\varepsilon\log(1/\varepsilon)$ to entropy. At the private scale $\varepsilon\asymp1/\sqrt{n\alpha^2}$, the square of this contribution is of order $\frac{\log^2(n\alpha^2)}{n\alpha^2}$, which gives~\eqref{eq:entropy-fixed-alphabet-minimax}. At the non-private scale $\varepsilon\asymp1/n$, the squared contribution is only of order $\frac{(\log n)^2}{n^2}$, which is negligible compared with the asymptotic non-private minimax risk $1/n$~\cite{WY}. Thus, the sensitivity of entropy to small probabilities determines the asymptotic rate under local privacy.

\begin{center}
\captionof{table}{Minimax risk bounds for entropy estimation under $\alpha$-LDP, with $\alpha\in(0,1)$.}
\label{tab:intro-entropy-bounds}
\begingroup
\setlength{\tabcolsep}{10pt}
\renewcommand{\arraystretch}{2.2}
\begin{tabular}{@{}lccc@{}}
\toprule
 & \shortstack{Previous\\lower bound} & Upper bound & Our lower bound \\
\midrule
\addlinespace[6pt]
\shortstack[l]{Non-asymptotic,\\up to log factors}
& $\displaystyle\frac{1}{\alpha^2 n}\wedge1$~\cite{AZ}
& $\displaystyle\frac{S^2}{\alpha^2 n}\wedge1$~[Prop.~\ref{prop:entropy-laplace-upper}]
& $\displaystyle\frac{S^2}{\alpha^2 n}\wedge1$~[Cor.~\ref{thm:entropy-lower}] \\
\addlinespace[6pt]
\midrule
\addlinespace[6pt]
\shortstack[l]{Asymptotic: \\ $n\alpha^2\to\infty$, fixed $S$}
& $\displaystyle\frac{1}{\alpha^2 n}$~\cite{AZ}
& $\displaystyle\frac{\log^2(n\alpha^2)}{n\alpha^2}$~[Prop.~\ref{prop:entropy-laplace-upper}]
& $\displaystyle\frac{\log^2(n\alpha^2)}{n\alpha^2}$~[Cor.~\ref{cor:entropy-complementary-lower}] \\
\addlinespace[6pt]
\bottomrule
\end{tabular}
\endgroup
\end{center}

\subsection{Open question.}
\label{sec:entropy-open-questions}
Our results, summarized in Table~\ref{tab:intro-entropy-bounds} above, characterize the asymptotic rate up to constants, but the non-asymptotic rate only up to log factors. While this improves our understanding of entropy estimation under LDP, interesting questions remain about the precise non-asymptotic rate. To examine these questions more closely, we first combine Theorems~\ref{thm:general-additive-lower-bound} and~\ref{thm:orientation-lower-bound} to give a lower bound with all log factors written explicitly. Its proof is given in Appendix~\ref{appendix:cmplement_lower_bound_entropy}.

\begin{corollary}
\label{cor:entropy-combined-lower}
There exist universal constants $c>0$, $S_0\geq4$ such that, for all $n\geq1$, $S\geq S_0$, and $\alpha\in(0,1)$ satisfying $n\alpha^2\geq S^2\log S$, we have
\[
R_{n,\alpha}(F_{\varphi_H})
\geq c\left\{
\frac{S^2}{n\alpha^2(\log S)^3}
+
\frac{S}{n\alpha^2}\log^2\!\left(\frac{e\,n\alpha^2}{S}\right)
\right\}.
\]
\end{corollary}

In this corollary, we have restricted attention to the regime $n\alpha^2\geq S^2\log S$ to simplify the lower bound. A natural question is whether an estimator can attain this bound up to a universal constant. Our upper bound from Section~\ref{sec:entropy-plugin} is larger, and we suspect that this gap reflects a limitation of the plug-in estimator itself, as in non-private entropy estimation~\cite{JiaoVenkatHanWeissman2017MLE}. One possible direction is to adapt the polynomial approximation methods used in~\cite{WY} for non-private entropy estimation. The idea is to approximate each entropy term $-p_k\log p_k$ by a polynomial in $p_k$, then estimate the powers $p_k^\ell$ appearing in that polynomial. Under LDP, these powers must be estimated from privatized observations, and the added noise may make high-degree terms difficult to estimate accurately. Analyzing this approach and determining whether it  attains our lower bound  raises further technical challenges, which lie beyond the scope of the present paper. We thus leave this  question open.

\section{Power-sum estimation under LDP}
\label{sec-power-sum}

Recall that the power-sum functional is $F_{\varphi_\gamma}(\boldsymbol p)=\sum_{k\in[S]}\varphi_\gamma(p_k)=\sum_{k\in[S]}p_k^\gamma$. For $\gamma\in(0,2)$, we study how its minimax quadratic risk $R_{n,\alpha}(F_{\varphi_\gamma})$ depends on the sample size $n$, the privacy parameter $\alpha$, and the alphabet size $S$. We exclude the trivial case $\gamma=1$, since $F_{\varphi_1}(\boldsymbol p)=1$ for every probability vector $\boldsymbol p$. We first present two estimators from~\cite{ButuceaIssartel2021} and their risk bounds (Section~\ref{subsec:power-upper-bounds}). We then apply our general lower bounds to study their minimax optimality (Section~\ref{subsec:power-minimax-optimality}). We discuss the resulting rates (Section~\ref{subsec:power-discussion-proposal}) and close with an open question (Section~\ref{subsec:power-open-question}). Table~\ref{tab:power-sum-bounds-summary} summarizes the bounds.


\subsection{Existing upper bounds.}
\label{subsec:power-upper-bounds}

The work of~\cite{ButuceaIssartel2021} provides two estimators for power sums under LDP: a Laplace plug-in estimator and a thresholded version. We recall their risk bounds, which will serve as a reference when studying minimax optimality.

\textit{Laplace plug-in estimator.} The first procedure of~\cite{ButuceaIssartel2021} uses the Laplace mechanism described in Section~\ref{sec:entropy-plugin} to produce privatized observations $Z_{1:n}$; it   estimates the probabilities $p_k$ from $Z_{1:n}$, and then plugs these estimates into the  functional to obtain $\sum_{k\in[S]}\widehat p_k^\gamma$. More precisely, their probability estimates are $\widehat p_k:=T_{[0,2]}(n^{-1}\sum_{i\in[n]}Z_{ik})$. Since $F_{\varphi_\gamma}(\boldsymbol p)\in[0,S^{1-\gamma}]$ for $\gamma\in(0,1)$, and $F_{\varphi_\gamma}(\boldsymbol p)\in[0,1]$ for $\gamma\in(1,2)$, we project the plug-in estimator $\sum_{k\in[S]}\widehat p_k^\gamma$ onto these intervals, and obtain the final estimator $\widehat F_{\varphi_\gamma}:=T_{[0,S^{(1-\gamma)\vee0}]}\bigl(\sum_{k\in[S]}\widehat p_k^\gamma\bigr)$.
The risk bound from~\cite[Theorem~2.1]{ButuceaIssartel2021}, combined with this projection, gives
\begin{equation}
\label{eq:power-plugin-upper-sublinear}
\sup_{\boldsymbol p\in\mathcal P_S}
\mathbb E_{\boldsymbol p,K_{1:n}^{\mathrm{Lap},\alpha}}
\!\left[
\bigl(\widehat F_{\varphi_\gamma}-F_{\varphi_\gamma}(\boldsymbol p)\bigr)^2
\right]
\lesssim_\gamma
\frac{S^2}{(\alpha^2 n)^\gamma}\wedge S^{2(1-\gamma)},
\qquad \text{for }\gamma\in(0,1),
\end{equation}
for $S\geq2$, $\alpha\in(0,1)$, and $\alpha^2 n\geq1$. Here, $\lesssim_\gamma$ hides  constants depending only on $\gamma$.
Informally, $S^2/(\alpha^2 n)^\gamma$ is a bias-type term that comes from estimating small probabilities, as explained in Section~\ref{subsec:power-discussion-proposal}. The variance term, of order $S/(\alpha^2 n)^\gamma$, is absorbed into this bound. The term $S^{2(1-\gamma)}$ is simply the worst-case squared error of the constant estimator $\widehat F=0$.

\textit{Thresholded estimator for powers $\gamma\in(1,2)$.} The above procedure tries to estimate even the small probabilities hidden by the privacy noise. But these probabilities contribute little to the functional $F_{\varphi_\gamma}(\boldsymbol p)=\sum_{k\in[S]}p_k^\gamma$ for $\gamma\in(1,2)$, and estimating them mainly adds noise to the estimator. A natural idea is therefore to discard them, and apply the plug-in only to probabilities detected above a threshold. This is the idea behind the thresholded estimator in~\cite[Section~2.4]{ButuceaIssartel2021}. It improves the polynomial  dependence on $S$ of the risk bound in the regime $\alpha^2 n<S^2$ of  moderate effective sample sizes  \cite[Theorem~2.3]{ButuceaIssartel2021}. We therefore use this thresholded estimator when $\alpha^2 n<S^2$, and the ordinary Laplace plug-in estimator otherwise. Denoting this estimator by $\widehat F_{\varphi_\gamma}^{\mathrm{th}}$, its quadratic risk satisfies \cite[Theorems~2.1 and~2.3]{ButuceaIssartel2021}
\begin{equation}
\label{eq:power-upper-reference}
\sup_{\boldsymbol p\in\mathcal P_S}
\mathbb E_{\boldsymbol p,K_{1:n}^{\mathrm{Lap},\alpha}}
\!\left[
\bigl(\widehat F_{\varphi_\gamma}^{\mathrm{th}}-F_{\varphi_\gamma}(\boldsymbol p)\bigr)^2
\right]
\begin{cases}
\displaystyle
\widetilde{\lesssim}_\gamma\;
\frac{S^2}{(\alpha^2 n)^\gamma}
\wedge\frac1{(\alpha^2 n)^{\gamma-1}},
& \text{for }\gamma\in(1,2),\quad \alpha^2 n\leq T_\gamma(S),\\[9pt]
\displaystyle
\lesssim_\gamma\;
\frac{S^{3-2\gamma}+1}{\alpha^2 n},
& \text{for }\gamma\in(1,2),\quad \alpha^2 n>T_\gamma(S),
\end{cases}
\end{equation}
for $S\geq2$, $\alpha\in(0,1)$, and $\alpha^2 n\geq1$, with $n\geq C_0\log S$ in the first line, for a universal constant $C_0$. Here, $\widetilde{\lesssim}_\gamma$ also suppresses log factors in $S$ and $n$, and $T_\gamma(S):=S^{(2\gamma-1)/(\gamma-1)}$ for $\gamma\in(1,3/2)$, while $T_\gamma(S):=S^{2/(\gamma-1)}$ for $\gamma\in[3/2,2)$. In what follows, we use the upper bounds in~\eqref{eq:power-plugin-upper-sublinear} and~\eqref{eq:power-upper-reference} as a reference, and study their minimax optimality.

\subsection{Minimax optimality.}
\label{subsec:power-minimax-optimality}

We begin with~\eqref{eq:power-plugin-upper-sublinear} and the first line of~\eqref{eq:power-upper-reference}. The existing lower bounds of~\cite{ButuceaIssartel2021} leave polynomial gaps in both the alphabet size $S$ and the effective sample size $\alpha^2 n$. Indeed, for $\gamma\in(0,1)\cup(1,2)$, their lower bounds have $S^{2-\gamma}/(\alpha^2 n)^\gamma$ instead of $S^2/(\alpha^2 n)^\gamma$, leaving a gap of a factor $S^\gamma$. For $\gamma\in(1,2)$, they also have $(\alpha^2 n)^{-2(\gamma-1)}$ instead of $(\alpha^2 n)^{-(\gamma-1)}$, leaving a gap of a factor $(\alpha^2 n)^{\gamma-1}$. To close these polynomial gaps, we apply Theorem~\ref{thm:general-additive-lower-bound} to $\varphi_\gamma(t)=t^\gamma$. This gives the following lower bounds, which match the upper bound in~\eqref{eq:power-plugin-upper-sublinear} and the first line of~\eqref{eq:power-upper-reference}, up to log factors.
The proof is deferred to Appendix~\ref{section:proof_power_sum_general_lb}.

\begin{corollary}
\label{cor:power-bias-lower-combined}
Let $\gamma\in(0,2)$, with $\gamma\neq1$. There exists a constant $S_\gamma\geq4$, depending only on $\gamma$, such that, for all $n\geq1$, $S\geq S_\gamma$, and $\alpha\in(0,1)$ satisfying $\alpha^2 n\geq S_\gamma^2$, we have
\[
R_{n,\alpha}(F_{\varphi_\gamma})
\widetilde{\gtrsim}_\gamma
\begin{cases}
\displaystyle
\frac{S^2}{(\alpha^2 n)^\gamma}
\wedge S^{2(1-\gamma)},
& \text{for }\gamma\in(0,1),\\[9pt]
\displaystyle
\frac{S^2}{(\alpha^2 n)^\gamma}
\wedge\frac1{(\alpha^2 n)^{\gamma-1}},
& \text{for }\gamma\in(1,2),
\end{cases}
\]
where $\widetilde{\gtrsim}_\gamma$ hides log factors in $(S, \alpha^2 n)$, and  positive constants depending only on $\gamma$.
\end{corollary}

We now focus on the upper bound in the second line of~\eqref{eq:power-upper-reference}. \cite{ButuceaIssartel2021} only provides a lower bound of order $1/(\alpha^2 n)$, thus leaving a gap of a factor $S^{3-2\gamma}$ unresolved. Applying our Theorem~\ref{thm:orientation-lower-bound} to $\varphi_\gamma(t)=t^\gamma$ closes this gap, with the following matching lower bound. Its proof is deferred to Appendix~\ref{appendix:power-variance-lower}.

\begin{corollary}
\label{cor:power-variance-lower}
Let $\gamma\in(1,2)$. For all $S\geq3$, $n\geq1$, and $\alpha\in(0,1)$ satisfying $\alpha^2 n>T_\gamma(S)$, we have $R_{n,\alpha}(F_{\varphi_\gamma})\gtrsim_\gamma (S^{3-2\gamma}+1)/(\alpha^2 n)$.
\end{corollary}

Together with the upper bounds~(\ref{eq:power-plugin-upper-sublinear}-\ref{eq:power-upper-reference}), Corollaries~\ref{cor:power-bias-lower-combined} and~\ref{cor:power-variance-lower} characterize the minimax rates up to log factors, and even up to constants in several regimes. Surprisingly, these optimal rates are attained by the plug-in estimator and its thresholded version presented in Section~\ref{subsec:power-upper-bounds}. We discuss these rates below in Section~\ref{subsec:power-discussion-proposal}.

\begin{center}
\captionof{table}{Minimax risk bounds for power-sum  estimation,  with $\alpha\in(0,1)$. Write $n_\alpha:=\alpha^2 n$.}
\label{tab:power-sum-bounds-summary}
\begingroup
\setlength{\tabcolsep}{8pt}
\renewcommand{\arraystretch}{2.2}
\begin{tabular*}{\textwidth}{@{\extracolsep{\fill}}lccc@{}}
\toprule
Regime
& \shortstack{Lower bound\\\cite{ButuceaIssartel2021}}
& \shortstack{Upper bound\\\cite{ButuceaIssartel2021}}
& Our lower bound \\
\midrule
\addlinespace[6pt]
$\gamma\in(0,1)$
& $\displaystyle \frac{S^{2-\gamma}}{n_\alpha^\gamma}\wedge S^{2(1-\gamma)}$
& $\displaystyle \frac{S^2}{n_\alpha^\gamma}\wedge S^{2(1-\gamma)}$
& $\displaystyle \frac{S^2}{n_\alpha^\gamma}\wedge S^{2(1-\gamma)}$~[Cor.~\ref{cor:power-bias-lower-combined}] \\
\addlinespace[6pt]
\midrule
\addlinespace[6pt]
\shortstack[l]{$\gamma\in(1,2)$\\[4pt] $n_\alpha\leq T_\gamma(S)$}
& $\displaystyle\Bigl(\frac{S^{2-\gamma}}{n_\alpha^\gamma}\wedge\frac{1}{n_\alpha^{2(\gamma-1)}}\Bigr)+\frac{1}{n_\alpha}$
& $\displaystyle\frac{S^2}{n_\alpha^\gamma}\wedge\frac{1}{n_\alpha^{\gamma-1}}$
& $\displaystyle\frac{S^2}{n_\alpha^\gamma}\wedge\frac{1}{n_\alpha^{\gamma-1}}$~[Cor.~\ref{cor:power-bias-lower-combined}] \\
\addlinespace[6pt]
\midrule
\addlinespace[6pt]
\shortstack[l]{$\gamma\in(1,2)$\\[4pt] $n_\alpha>T_\gamma(S)$}
& $\displaystyle\frac{1}{n_\alpha}$
& $\displaystyle\frac{S^{3-2\gamma}+1}{n_\alpha}$
& $\displaystyle\frac{S^{3-2\gamma}+1}{n_\alpha}$~[Cor.~\ref{cor:power-variance-lower}] \\
\addlinespace[6pt]
\bottomrule
\end{tabular*}
\endgroup
\end{center}

\subsection{Discussion of the rates.}
\label{subsec:power-discussion-proposal}

As for entropy in Section~\ref{sec:entropy-discussion}, the difficulty of estimating small probabilities helps explain the power-sum rates. The difference is that their contribution now depends on $\gamma$: the function $x^\gamma$ is concave for $\gamma\in(0,1)$, and convex for $\gamma\in(1,2)$. Write $\varepsilon:=1/\sqrt{\alpha^2 n}$ for the probability estimation error introduced by a standard Laplace mechanism, such as the one described in Section~\ref{sec:entropy-plugin}. Probabilities $p_k\lesssim\varepsilon$, comparable to this error or smaller, are difficult to estimate.

\textit{Small probabilities for $\gamma\in(0,1)$.} The two terms in the minimax rate for $\gamma\in(0,1)$ in Table~\ref{tab:power-sum-bounds-summary} correspond to two configurations. (i) When $S\varepsilon\leq1$, a distribution $\boldsymbol p$ can have a number $\asymp S$ of probabilities of size $p_k\asymp\varepsilon$. Their aggregate contribution to $F_{\varphi_\gamma}$ is of order $S\varepsilon^\gamma$, suggesting a squared error of order $S^2\varepsilon^{2\gamma}=S^2/(\alpha^2 n)^\gamma$. (ii) When $S\varepsilon>1$, a distribution can have $S$ probabilities of order $1/S$, all below the noise level. Their aggregate contribution can be of order $S^{1-\gamma}$, suggesting the risk term $S^{2(1-\gamma)}$. Together, these two cases recover the rate in Table~\ref{tab:power-sum-bounds-summary} for $\gamma\in(0,1)$.   

\textit{What changes for $\gamma\in(1,2)$.} When $S\varepsilon\leq1$, the same reasoning explains the term $S^2/(\alpha^2 n)^\gamma$ in the second line of Table~\ref{tab:power-sum-bounds-summary}. The difference appears when $S\varepsilon>1$. Since $x\mapsto x^\gamma$ is convex, concentrating the mass on fewer probabilities near the noise level gives a larger contribution to the functional than spreading it over $S$ probabilities of size $1/S$. We therefore consider about $1/\varepsilon<S$ categories, each having probability of order $\varepsilon$. Their aggregate contribution to $F_{\varphi_\gamma}$ is of order $(1/\varepsilon)\varepsilon^\gamma=\varepsilon^{\gamma-1}$, whose square gives $(\alpha^2 n)^{-(\gamma-1)}$. Together, these two cases recover the rate in the second line of Table~\ref{tab:power-sum-bounds-summary}. 

\textit{Large effective sample sizes $\alpha^2 n>T_\gamma(S)$.} The rate in the third line of Table~\ref{tab:power-sum-bounds-summary} comes from a different source of difficulty, related to variance rather than bias. Indeed, even probabilities above the noise level are estimated with error, and these fluctuations affect the functional. Specifically, for $p_k\gg\varepsilon$, a first-order expansion gives $\widehat p_k^\gamma-p_k^\gamma\approx\gamma p_k^{\gamma-1}(\widehat p_k-p_k)$, and since the Laplace noise introduces an error of typical size $|\widehat p_k-p_k|\asymp\varepsilon$, this expansion gives fluctuations of size $p_k^{\gamma-1}\varepsilon$, and hence a variance of order $p_k^{2\gamma-2}\varepsilon^2$. Since the added noises are independent across categories, these variances add up, giving a contribution of order $\varepsilon^2\sum_{k\in[S]}p_k^{2\gamma-2}$. (i) For $\gamma\in(1,3/2)$, concavity of $x\mapsto x^{2\gamma-2}$ shows that this sum is maximized at the uniform distribution $p_k=1/S$, giving the risk term $\varepsilon^2S^{3-2\gamma}=S^{3-2\gamma}/(\alpha^2 n)$. Thus, this term comes from fluctuations across many categories, even though their probabilities are above the  privacy noise (they satisfy $p_k=1/S>\varepsilon$, since we are in the regime $\alpha^2 n>T_\gamma(S)$ and $T_\gamma(S)\geq S^2$).  (ii) For $\gamma\in[3/2,2)$, the inequalities $p_k^{2\gamma-2}\leq p_k$, valid because $p_k\in[0,1]$ and $2\gamma-2\geq1$, give the variance contribution $\varepsilon^2\sum_{k\in[S]}p_k^{2\gamma-2}\leq\varepsilon^2\sum_{k\in[S]}p_k=\varepsilon^2$, which is therefore at most of order $\varepsilon^2=1/(\alpha^2 n)$. An error of this order occurs, for example, when $p_1=p_2=1/2$ and all other probabilities are zero. Together, cases (i) and (ii) recover the rate in the third line of Table~\ref{tab:power-sum-bounds-summary}.

Although these arguments give a simple interpretation of the rates, proving the lower bounds requires a more delicate analysis, detailed in Sections~\ref{sec:proof-general-lower-bound}--\ref{proof:orientation-lower-bound}.

\subsection{Open question.}
\label{subsec:power-open-question}
The plug-in estimator and its thresholded version (Section~\ref{subsec:power-upper-bounds}) achieve the optimal polynomial dependence on $(S,n,\alpha)$, but can we improve the log factors in the two unresolved regimes (the first two rows of Table~\ref{tab:power-sum-bounds-summary})? One possibility is to adapt the polynomial approximation method used for non-private power-sum estimation in~\cite{JiaoVenkatHanWeissman2015}, by approximating $x^\gamma$ with a low-degree polynomial and estimating the resulting polynomial from privatized observations. As discussed for entropy in Section~\ref{sec:entropy-open-questions}, analyzing this approach would require a separate study. We therefore leave this question open.



\FloatBarrier
\section{Proof of Theorem \ref{thm:general-additive-lower-bound}}
\label{sec:proof-general-lower-bound}

The proof combines a new bound on the total variation distance between privatized mixtures (Section~\ref{subsec:controlling-privatized-mixtures}) with  classical tools: the duality between polynomial approximation and moment matching, the two fuzzy hypotheses method, and an elementary concentration bound (Section~\ref{subsec:classical-tools}). We then combine these ingredients to prove Theorem~\ref{thm:general-additive-lower-bound} (Section~\ref{proof_thm1}).

\textit{Notation.} We write $\|P_1-P_0\|_{\mathrm{TV}}:=\sup_{A\text{ measurable}}|P_1(A)-P_0(A)|$ for the total variation distance between two probability measures $P_0$ and $P_1$ on the same measurable space.

\subsection{Controlling privatized mixtures}
\label{subsec:controlling-privatized-mixtures}

The next theorem addresses the difficulty described in Section~\ref{subsec:reduction-one-dimensional-mixtures} by bounding the total variation distance between joint privatized mixtures, without requiring these mixture laws to factorize as products of mixtures over individual parameters. The notation for the priors $\bigotimes_j\pi_j$ and $\bigotimes_j\pi'_j$ used in Sections~\ref{sec:intro-methodology} and~\ref{sec2:general-lb-thms} will be replaced below by $\mu_0^{\otimes m}$ and $\mu_1^{\otimes m}$, respectively.

Fix $\boldsymbol p\in\mathcal P_{m,b}^{\mathrm{pair}}$ and consider its associated family $\{\boldsymbol p(\boldsymbol t):\boldsymbol t\in[-1,1]^m\}$ from Definition~\ref{def:pairwise-perturbations}. Fix also a privacy mechanism $K_{1:n}\in\mathcal M_\alpha^n$. For each fixed $\boldsymbol t$, the privatized observations are independent, with
\[
P_{\boldsymbol p(\boldsymbol t),K_{1:n}}^{Z_{1:n}}
=\bigotimes_{i=1}^n P_{\boldsymbol p(\boldsymbol t),K_{1:n}}^{Z_i},
\qquad \qquad
P_{\boldsymbol p(\boldsymbol t),K_{1:n}}^{Z_i}(A_i)
=\sum_{k\in[S]}p_k(\boldsymbol t)K_i(A_i\mid k),
\]
for measurable $A_i\subseteq\mathcal Z_i$.
Let $\mu_0,\mu_1$ be two probability measures on $[0,1]$. For $r\in\{0,1\}$, drawing $\boldsymbol t\sim\mu_r^{\otimes m}$ gives the mixture law
\begin{equation}
\label{eq:privatized-mixture-laws}
M_{r,\boldsymbol p,K_{1:n}}^{Z_{1:n}}
:=
\mathbb E_{\boldsymbol t\sim\mu_r^{\otimes m}}
\!\left[P_{\boldsymbol p(\boldsymbol t),K_{1:n}}^{Z_{1:n}}\right]
=
\int_{[0,1]^m}
P_{\boldsymbol p(\boldsymbol t),K_{1:n}}^{Z_{1:n}}\,
\mu_r^{\otimes m}(d\boldsymbol t).
\end{equation}
The following theorem bounds the total variation distance between these two mixture laws, uniformly over $\boldsymbol p\in\mathcal P_{m,b}^{\mathrm{pair}}$ and $K_{1:n}\in\mathcal M_\alpha^n$. The bound is expressed in terms of the moment differences
\begin{equation}
\label{eq:moment-differences}
\delta_\ell
:=
\int_{[0,1]}t^\ell\,\mu_1(dt)
-
\int_{[0,1]}t^\ell\,\mu_0(dt),
\qquad \ell\geq0.
\end{equation}

\begin{theorem}
\label{thm:private-mixture-comparison}
Let $S\geq3$, $n\geq1$, and $\alpha>0$. Let  $m\geq1$ and $b>0$ define a nonempty class $\mathcal P_{m,b}^{\mathrm{pair}}$ as in Definition~\ref{def:pairwise-perturbations}. For any two probability measures $\mu_0,\mu_1$ on $[0,1]$, the mixture laws  in~\eqref{eq:privatized-mixture-laws} satisfy
\begin{equation}
\label{eq:private-mixture-comparison}
\sup_{\boldsymbol p\in\mathcal P_{m,b}^{\mathrm{pair}}}
\sup_{K_{1:n}\in\mathcal M_\alpha^n}
\left\|
M_{1,\boldsymbol p,K_{1:n}}^{Z_{1:n}}
-
M_{0,\boldsymbol p,K_{1:n}}^{Z_{1:n}}
\right\|_{\mathrm{TV}}
\leq
\frac m2
\sqrt{
\sum_{\ell=1}^{n}
\binom n\ell
\delta_\ell^2[b(e^\alpha-1)]^{2\ell}
},
\end{equation}
where $\delta_\ell$ denotes the difference between the moments of order $\ell$ of $\mu_1$ and $\mu_0$, as defined in~\eqref{eq:moment-differences}. 
\end{theorem}

If $\mu_0$ and $\mu_1$ have the same first $L$ moments, then $\delta_\ell=0$ for $\ell\leq L$, and the first $L$ terms in the sum in~\eqref{eq:private-mixture-comparison} vanish; in particular, the bound is zero if $n\leq L$. The bound therefore depends only on the remaining moment differences $\delta_\ell$, for $\ell>L$. For $0<\alpha\leq1$, we have $e^\alpha-1\leq2\alpha$, so this sum is bounded by $\sum_{\ell=L+1}^{n}\binom n\ell\delta_\ell^2(4b^2\alpha^2)^\ell$. When $b\alpha\ll1$, we can then choose how many moments $L$ to match to make the total variation distance small. The bound holds uniformly over all $\boldsymbol p\in\mathcal P_{m,b}^{\mathrm{pair}}$ and $K_{1:n}\in\mathcal M_\alpha^n$.

Theorem~\ref{thm:private-mixture-comparison} is a key contribution of our work, and its proof is given in Section~\ref{sec:proof-private-mixture-comparison}. The proof exploits the affine dependence of the privatized observation laws $P_{\boldsymbol p(\boldsymbol t),K_{1:n}}^{Z_i}$ on the perturbation parameters $\boldsymbol t$. By independence for fixed $\boldsymbol t$, the joint law $P_{\boldsymbol p(\boldsymbol t),K_{1:n}}^{Z_{1:n}}=\bigotimes_{i=1}^n P_{\boldsymbol p(\boldsymbol t),K_{1:n}}^{Z_i}$ is a product of these affine laws and therefore depends polynomially on $\boldsymbol t$. Matching moments, $\delta_\ell=0$ for $\ell\leq L$, then cancels the contributions of polynomial terms of degree at most $L$. To bound the remaining contributions, we use the $\alpha$-LDP condition, which limits how much the privatized laws $P_{\boldsymbol p(\boldsymbol t),K_{1:n}}^{Z_i}$ can change as $\boldsymbol t$ varies.

\subsection{Classical tools}
\label{subsec:classical-tools}

We recall two classical tools commonly used to prove minimax lower bounds in estimation problems without privacy constraints, together with an elementary concentration bound.

\textit{Polynomial approximation and moment matching.}\label{subsec:approximation-moment-matching} We recall the duality between polynomial approximation and moment matching discussed in Section~\ref{subsec:mixture-approach-limitations}; see \cite{WY,CCT} for details and proofs. Recall that $E_L(f)$, defined in Section~\ref{subsec:polynomial-lower-bound}, denotes the best uniform approximation error of $f$ on $[0,1]$ by polynomials of degree at most $L$. We use $\delta_\ell$, defined in~\eqref{eq:moment-differences}, for the difference between the moments of order $\ell$ of two probability measures $\mu_1$ and $\mu_0$ on $[0,1]$.

\begin{lemma}
\label{lem:approximation-moment-matching}
Let $f\in C([0,1])$ and let $L\geq0$ be an integer. Then, we have
\begin{equation}
\label{eq:approximation-matching-separation}
\max_{\substack{
\mu_0,\mu_1\in\mathcal P([0,1])\\
\delta_\ell=0,\ \ell=1,\ldots,L
}}
\left|
\int_{[0,1]}f(t)\,\mu_1(dt)
-
\int_{[0,1]}f(t)\,\mu_0(dt)
\right|
=2E_L(f),
\end{equation}
where the maximum is over pairs of probability measures on $[0,1]$ with the same first $L$ moments.
\end{lemma}


\textit{Two fuzzy hypotheses.}\label{subsec:two-fuzzy-hypotheses} Instead of comparing two fixed parameters, this method compares two priors $\pi_0,\pi_1$. The functional $F$ must take separated values with high probability under these priors, while the observation mixtures $M_0,M_1$ induced by $\pi_0,\pi_1$, respectively, must remain close in TV distance. We recall the following squared-error version, adapted from \cite{Tsybakov2009}.

\begin{lemma}
\label{lem:two-fuzzy-hypotheses}
Let $\Theta$ and $\mathcal Z$ be measurable spaces, let $Q$ be a Markov kernel from $\Theta$ to $\mathcal Z$, and let $F:\Theta\to\mathbb R$ be measurable. For two probability measures $\pi_0,\pi_1$ on $\Theta$, define
\[
M_r:=\int_\Theta Q(\cdot\mid\theta)\,\pi_r(d\theta),
\qquad r\in\{0,1\}.
\]
Suppose that there exist a threshold $c\in\mathbb R$, a separation parameter $\Delta>0$, and $\eta_0,\eta_1\in[0,1]$ such that
\[
\pi_0\bigl(\{\theta:F(\theta)\leq c-\Delta\}\bigr)\geq1-\eta_0,
\qquad \qquad
\pi_1\bigl(\{\theta:F(\theta)\geq c+\Delta\}\bigr)\geq1-\eta_1.
\]
Then, every measurable estimator $\widehat F:\mathcal Z\to\mathbb R$ satisfies
\[
\sup_{\theta\in\Theta}
\int
\bigl(\widehat F(z)-F(\theta)\bigr)^2\,Q(dz\mid\theta)
\ \geq \ 
\frac{\Delta^2}{2}
\left(1-\|M_1-M_0\|_{\mathrm{TV}}-\eta_0-\eta_1\right).
\]
\end{lemma}



\textit{Concentration of sums of functions of independent random variables.} The following elementary lemma controls how far a sum $\sum_{j=1}^m F_j(U_j)$ of functions of independent random variables can deviate from its mean, in terms of the oscillations $V(F_j)$. Recall from Section~\ref{subsec:polynomial-lower-bound} that $V(f)=\sup_{u\in[0,1]}f(u)-\inf_{u\in[0,1]}f(u)$, for any bounded function $f:[0,1]\to\mathbb R$.
The proof is  in Appendix~\ref{sec:proofs-main-tools}.

\begin{lemma}
\label{lem:product-measure-concentration}
Let $\mu$ be a probability measure on $[0,1]$, let $m\geq1$,
and  let $U=(U_1,\ldots,U_m)$ be a random vector with
distribution $\mu^{\otimes m}$.
Let $F_1,\ldots,F_m:[0,1]\to\mathbb R$ be bounded measurable functions, and define
$F(u)=\sum_{j\in[m]} F_j(u_j)$, for
$u=(u_1,\ldots,u_m)\in[0,1]^m$.
Then, for any $\Delta>0$, 
\[
\mathbb P_{\mu^{\otimes m}}\left(
\left|F(U)-\mathbb E_{\mu^{\otimes m}}[F(U)]\right|
\geq\frac{\Delta}{4}
\right)
\ \leq \
\frac{4\sum_{j\in[m]}V(F_j)^2}
{\Delta^{\,2}}.
\]
\end{lemma}

\subsection{Proof of Theorem \ref{thm:general-additive-lower-bound}}
\label{proof_thm1}

The proof of Theorem~\ref{thm:general-additive-lower-bound} proceeds in four steps. Step~1 uses the duality between polynomial approximation and moment matching (Lemma~\ref{lem:approximation-moment-matching}) to construct priors with separated functional means. Step~2 applies the two fuzzy hypotheses lemma (Lemma~\ref{lem:two-fuzzy-hypotheses}) to obtain a risk lower bound from this separation. This leaves us with two quantities to control: the fluctuations of the functional around its mean under each prior, and the TV distance between the privatized mixture laws. Step~3 controls these fluctuations using the concentration bound of Lemma~\ref{lem:product-measure-concentration}. Step~4 controls this TV distance by applying the moment-based TV bound (Theorem~\ref{thm:private-mixture-comparison}).

\begin{proof}[Proof of Theorem~\ref{thm:general-additive-lower-bound}]
If $E_L(\varphi_b)=0$, the lower bound is zero and the conclusion is immediate. Assume henceforth that $E_L(\varphi_b)>0$. Fix $\boldsymbol p\in\mathcal P_{m,b}^{\mathrm{pair}}$ and consider the family $\boldsymbol p(\boldsymbol t)$ from Definition~\ref{def:pairwise-perturbations}. The last $S-2m$ probabilities remain fixed, so, for $\boldsymbol t\in[0,1]^m$, the definition of the pair profile gives
\begin{equation}
\label{eq:general-lb-decomposition}
F_\varphi(\boldsymbol p(\boldsymbol t))
\ = \
\sum_{k=2m+1}^{S}\varphi(p_k)
+\sum_{j=1}^{m}\varphi_b(t_j).
\end{equation}

\textit{Step 1: Construction of the priors, and separation of the means.}
Recall the notation $\delta_\ell$ in~\eqref{eq:moment-differences} for the difference between moments of order $\ell$. By Lemma~\ref{lem:approximation-moment-matching}, applied to $\varphi_b$, there exist two probability measures $\mu_0,\mu_1$ on $[0,1]$ with the same first $L$ moments, that is, $\delta_\ell=0$ for $\ell=1,\ldots,L$, and satisfying the separation between the expectations of $\varphi_b$: $\int_{[0,1]}\varphi_b(t)\,\mu_1(dt)-\int_{[0,1]}\varphi_b(t)\,\mu_0(dt)\ =\ 2E_L(\varphi_b)$.
We use the product measures $\mu_0^{\otimes m}$ and $\mu_1^{\otimes m}$ as priors on $\boldsymbol t\in[0,1]^m$. Define the difference between the functional means under these priors by $\Delta\ :=\ \mathbb E_{\boldsymbol t\sim\mu_1^{\otimes m}}\bigl[F_\varphi(\boldsymbol p(\boldsymbol t))\bigr]-\mathbb E_{\boldsymbol t\sim\mu_0^{\otimes m}}\bigl[F_\varphi(\boldsymbol p(\boldsymbol t))\bigr]$.
By \eqref{eq:general-lb-decomposition} and the separation $2E_L(\varphi_b)$ obtained above, we obtain
\begin{equation}
\label{eq:general-lb-separation}
\Delta
\ = \
m\left(
\int \varphi_b(t)\,\mu_1(dt)
-
\int \varphi_b(t)\,\mu_0(dt)
\right)
\ = \
2mE_L(\varphi_b).
\end{equation}

\textit{Step 2: Applying the two fuzzy hypotheses method.}
Fix a mechanism $K_{1:n}\in\mathcal M_\alpha^n$. For $r\in\{0,1\}$, let $M_{r,\boldsymbol p,K_{1:n}}^{Z_{1:n}}:=\mathbb E_{\boldsymbol t\sim\mu_r^{\otimes m}}\bigl[P_{\boldsymbol p(\boldsymbol t),K_{1:n}}^{Z_{1:n}}\bigr]$ be the observation mixture under the prior $\mu_r^{\otimes m}$ constructed in Step~1.
For $r\in\{0,1\}$, set
\[
\eta_r
\ := \
\mu_r^{\otimes m}\left(
\left\{
\boldsymbol t:
\left|
F_\varphi(\boldsymbol p(\boldsymbol t))
-
\mathbb E_{\boldsymbol v\sim\mu_r^{\otimes m}}
\bigl[F_\varphi(\boldsymbol p(\boldsymbol v))\bigr]
\right|
\geq\frac{\Delta}{4}
\right\}
\right).
\]
Take the threshold
\[
c
\ := \
\frac12\left(
\mathbb E_{\boldsymbol t\sim\mu_0^{\otimes m}}
\bigl[F_\varphi(\boldsymbol p(\boldsymbol t))\bigr]
+
\mathbb E_{\boldsymbol t\sim\mu_1^{\otimes m}}
\bigl[F_\varphi(\boldsymbol p(\boldsymbol t))\bigr]
\right).
\]
By the definition of $\Delta$ in Step~1 and the definition of $c$, the means of $F_\varphi(\boldsymbol p(\boldsymbol t))$ under $\mu_0^{\otimes m}$ and $\mu_1^{\otimes m}$ are $c-\Delta/2$ and $c+\Delta/2$, respectively. A deviation of less than $\Delta/4$ from the first mean implies $F_\varphi(\boldsymbol p(\boldsymbol t))\leq c-\Delta/4$, and a deviation of less than $\Delta/4$ from the second implies $F_\varphi(\boldsymbol p(\boldsymbol t))\geq c+\Delta/4$. Thus, by the definition of $\eta_r$,
\[
\begin{aligned}
\mu_0^{\otimes m}\bigl(\{\boldsymbol t:F_\varphi(\boldsymbol p(\boldsymbol t))\leq c-\Delta/4\}\bigr)
&\ \geq \ 1-\eta_0,\\
\mu_1^{\otimes m}\bigl(\{\boldsymbol t:F_\varphi(\boldsymbol p(\boldsymbol t))\geq c+\Delta/4\}\bigr)
&\ \geq \ 1-\eta_1.
\end{aligned}
\]
We apply the two fuzzy hypotheses lemma (Lemma~\ref{lem:two-fuzzy-hypotheses}) to the functional $\boldsymbol t\mapsto F_\varphi(\boldsymbol p(\boldsymbol t))$, the Markov kernel $\boldsymbol t\mapsto P_{\boldsymbol p(\boldsymbol t),K_{1:n}}^{Z_{1:n}}$, and the priors $\mu_0^{\otimes m},\mu_1^{\otimes m}$, with threshold $c$ and separation parameter $\Delta/4$. Since $\boldsymbol p(\boldsymbol t)\in\mathcal P_{m,b}^{\mathrm{pair}}$, this gives, for every measurable estimator $\widehat F$,
\begin{align}
\sup_{\boldsymbol q\in\mathcal P_{m,b}^{\mathrm{pair}}}
\mathbb E_{\boldsymbol q,K_{1:n}}
\bigl[(\widehat F(Z_{1:n})-F_\varphi(\boldsymbol q))^2\bigr]
&\ \geq \
\sup_{\boldsymbol t\in[0,1]^m}
\int
\bigl(\widehat F(z_{1:n})-F_\varphi(\boldsymbol p(\boldsymbol t))\bigr)^2
\,P_{\boldsymbol p(\boldsymbol t),K_{1:n}}^{Z_{1:n}}(dz_{1:n})\notag\\
&\ \geq \
\frac{\Delta^2}{32}
\left(
1-\left\|M_{1,\boldsymbol p,K_{1:n}}^{Z_{1:n}}
-M_{0,\boldsymbol p,K_{1:n}}^{Z_{1:n}}\right\|_{\mathrm{TV}}
-\eta_0-\eta_1
\right).\label{eq:general-lb-fuzzy-reduction}
\end{align}

\textit{Step 3: Concentration under the priors.}
By \eqref{eq:general-lb-decomposition}, $F_\varphi(\boldsymbol p(\boldsymbol t))$ differs from $\sum_{j\in[m]}\varphi_b(t_j)$ by a constant, which does not affect deviations from the mean. Hence, applying the concentration bound of Lemma~\ref{lem:product-measure-concentration} with $F_j=\varphi_b$ for every $j\in[m]$ gives
\begin{equation}
\label{eq:general-lb-prior-concentration}
\eta_r
\ \leq \
\frac{4mV(\varphi_b)^2}{\Delta^2}
\ = \
\frac{V(\varphi_b)^2}{mE_L(\varphi_b)^2},
\qquad \quad r\in\{0,1\},
\end{equation}
where the equality follows from \eqref{eq:general-lb-separation}.

\textit{Step 4: Controlling the TV distance.}
By the moment matching condition in Step~1, $\delta_\ell=0$ for $1\leq\ell\leq L$. Moreover, $|\delta_\ell|\leq1$, since $0\leq t^\ell\leq1$ on $[0,1]$. The moment-based TV bound (Theorem~\ref{thm:private-mixture-comparison}) therefore gives
\[
\left\|M_{1,\boldsymbol p,K_{1:n}}^{Z_{1:n}}
-M_{0,\boldsymbol p,K_{1:n}}^{Z_{1:n}}\right\|_{\mathrm{TV}}
\ \leq \
\frac m2
\left\{
\sum_{\ell=L+1}^{n}
\binom n\ell
[b(e^\alpha-1)]^{2\ell}
\right\}^{1/2},
\]
where the sum is zero when $n\leq L$.
Since $0<\alpha<1$, we have $e^\alpha-1\leq(e-1)\alpha$. The condition $b\leq C_0\sqrt{(L+1)/(n\alpha^2)}$ in~\eqref{eq:general-lb-parameters}, with $C_0=(2\sqrt e\,(e-1))^{-1}$, then gives $enb^2(e^\alpha-1)^2\ \leq\ enb^2(e-1)^2\alpha^2\ \leq\ (L+1)/4$.
Using $\binom n\ell\leq(en/\ell)^\ell$, we obtain, for
$L+1\leq\ell\leq n$, 
\[
\binom n\ell[b(e^\alpha-1)]^{2\ell}
\ \leq \
\left(\frac{enb^2(e^\alpha-1)^2}{\ell}\right)^\ell
\ \leq \
\left(\frac{L+1}{4\ell}\right)^\ell
\ \leq \
4^{-\ell}.
\]
Therefore,
\begin{equation}
\label{eq:general-lb-mixture-tv}
\left\|M_{1,\boldsymbol p,K_{1:n}}^{Z_{1:n}}
-M_{0,\boldsymbol p,K_{1:n}}^{Z_{1:n}}\right\|_{\mathrm{TV}}
\ \leq \
\frac m2
\left(\sum_{\ell=L+1}^{\infty}4^{-\ell}\right)^{1/2}
\ = \
\frac{m2^{-L}}{2\sqrt3}
\ \leq \
\frac1{4\sqrt3}
\ < \
\frac14,
\end{equation}
where the second inequality follows from $m2^{-L}\leq1/2$, since $L\geq\lceil\log_2(2m)\rceil$ by~\eqref{eq:general-lb-parameters}.


\textit{Putting everything together.}
Substituting \eqref{eq:general-lb-separation}, \eqref{eq:general-lb-prior-concentration} and
\eqref{eq:general-lb-mixture-tv} into
\eqref{eq:general-lb-fuzzy-reduction} yields
\[
\begin{aligned}
\sup_{\boldsymbol q\in\mathcal P_{m,b}^{\mathrm{pair}}}
\mathbb E_{\boldsymbol q,K_{1:n}}
\bigl[(\widehat F(Z_{1:n})-F_\varphi(\boldsymbol q))^2\bigr]
&\ \geq \
\frac{m^2E_L(\varphi_b)^2}{8}
\left(
1-\frac14
-\frac{2V(\varphi_b)^2}{mE_L(\varphi_b)^2}
\right)\\
&\ = \
\frac{3m^2}{32}E_L(\varphi_b)^2
-\frac m4V(\varphi_b)^2.
\end{aligned}
\]
Taking the infimum over $K_{1:n}\in\mathcal M_\alpha^n$ and all measurable estimators $\widehat F$, and using the nonnegativity of the risk, gives the lower bound of Theorem~\ref{thm:general-additive-lower-bound} and concludes the proof.
\end{proof}

\section{Proof of the TV bound for privatized mixtures}
\label{sec:proof-private-mixture-comparison}

This section is devoted to the proof of Theorem~\ref{thm:private-mixture-comparison}. We first develop general tools (in Section~\ref{sec:main-technical-tools}) for bounding the TV distance between mixtures, and for bounding the $\chi^2$ divergence between privatized laws. These tools may be of independent interest, beyond the statistical problem considered in this paper. We then combine them to prove Theorem~\ref{thm:private-mixture-comparison} (in Section~\ref{proo:general_tv_bound}).

\textit{Notation.} For two probability measures $\mu,\nu$ on the same measurable
space, define their $\chi^2$ divergence:
$\chi^2(\nu\|\mu):=\int\bigl(\frac{d\nu}{d\mu}-1\bigr)^2\,d\mu$, if $\nu\ll\mu$,
and $\chi^2(\nu\|\mu):=+\infty$ otherwise.
For a Markov kernel $Q$ and a probability measure $\mu$ on its input space,
we define the mixture law $Q\circ\mu$ by
$(Q\circ\mu)(A):=\int Q(A\mid u)\,\mu(du)$,
for every measurable set $A$ in the output space of $Q$.

\subsection{Main tools}
\label{sec:main-technical-tools}

The three tools below are each stated in their own general setting.
The first reduces the problem of bounding
the TV distance between mixtures for multidimensional kernels to one-dimensional kernels.
The second gives a moment-based TV bound for mixtures associated with one-dimensional kernels.
The third establishes a $\chi^2$--$\mathrm{TV}^2$ Lipschitz bound for privacy kernels.

\subsubsection{Reduction to one-dimensional kernels}

Let $m\geq1$ and consider the product probability measures $\mu:=\bigotimes_{j\in[m]}\mu_j$ and $\nu:=\bigotimes_{j\in[m]}\nu_j$, where $\mu_j$ and $\nu_j$ are probability measures on a common measurable space $\mathcal T_j$, for $j\in[m]$.
Let $Q$ be a Markov kernel from $\prod_{j=1}^m\mathcal T_j$
to an arbitrary measurable space.
Comparing the $m$-dimensional mixtures $Q\circ\mu$ and $Q\circ\nu$
can be reduced to comparing the following $1$-dimensional mixtures with respect to $\mu_j$ and $\nu_j$.
For $\boldsymbol t=(t_1,\ldots,t_m)\in\prod_{k=1}^m\mathcal T_k$ and $j\in[m]$, let $\boldsymbol t_{-j}=(t_k)_{k\neq j}$ denote the coordinates other than $t_j$. Holding these coordinates fixed, we define the $1$-dimensional Markov kernel
\[
Q^{\boldsymbol t_{-j}}(\,\cdot\mid t_j)
\ :=
\
Q(\,\cdot\mid t_1,\ldots,t_m).
\]
The resulting laws $Q^{\boldsymbol t_{-j}}\circ\mu_j$ and
$Q^{\boldsymbol t_{-j}}\circ\nu_j$ are thus $1$-dimensional mixtures over $\mathcal T_j$.
As shown in the next lemma, the TV distances between these $1$-dimensional mixtures control that between the $m$-dimensional mixtures $Q\circ\mu$ and $Q\circ\nu$.

\begin{lemma}
\label{lem:coordinatewise-mixture-comparison}
Consider the above setting. Suppose that, for each $j\in[m]$, there exists
$\varepsilon_j\geq0$ such that
\[
\sup_{\boldsymbol t_{-j}\in\prod_{k\neq j}\mathcal T_k}
\left\|
Q^{\boldsymbol t_{-j}}\circ\nu_j-Q^{\boldsymbol t_{-j}}\circ\mu_j
\right\|_{\mathrm{TV}}
\ \leq \
\varepsilon_j\,.
\]
Then, we have
$\|Q\circ\nu-Q\circ\mu\|_{\mathrm{TV}}
\ \leq \
\sum_{j\in[m]}\varepsilon_j$.
\end{lemma}

Note that this reduction holds for any Markov kernel $Q$. It can be viewed as a triangle inequality along the coordinates $j\in[m]$. This becomes clear in the proof in Appendix~\ref{sec:proofs-main-tools}.


\subsubsection{A moment-based TV bound for affine one-dimensional kernels}

We now focus on $1$-dimensional mixtures; for the purposes of this paper, it suffices to consider mixing distributions supported on $[0,1]$.
Formally, let $\mu,\nu$ be two probability measures on $[0,1]$. Fix $n\geq1$, let $\mathcal Z_i$ be a measurable space for each $i\in[n]$, and let $Q$ be a Markov kernel from $[0,1]$ to $\prod_{i=1}^n\mathcal Z_i$, equipped with the product $\sigma$-algebra.
The next proposition bounds the TV distance between the mixtures $Q\circ\mu$ and $Q\circ\nu$ under the following assumptions.
Assume that $Q(\,\cdot\mid t)=\bigotimes_{i=1}^nQ_i(\,\cdot\mid t)$ for every $t\in[0,1]$, where each $Q_i$ is a Markov kernel from $[0,1]$ to $\mathcal Z_i$ that is affine in $t$:
\[
Q_i(\,\cdot\mid t)\ =\ (1-t)Q_i(\,\cdot\mid0)+tQ_i(\,\cdot\mid1),
\qquad t\in[0,1],\quad i\in[n]\,.
\]
Here, write $\delta_\ell:=\int_{[0,1]}t^\ell\,\nu(dt)-\int_{[0,1]}t^\ell\,\mu(dt)$ for the difference between the moments of order $\ell$ of $\nu$ and $\mu$, for $\ell\geq0$.
Proposition~\ref{lem:mixtures-product-measures} bounds the TV distance between $Q\circ\mu$ and $Q\circ\nu$ in terms of these moment differences $\delta_\ell$, and the $\chi^2$ divergences between the endpoint distributions $Q_i(\cdot\mid1)$ and $Q_i(\cdot\mid0)$.

\begin{proposition}
\label{lem:mixtures-product-measures}
In the above setting, we have
\[
\|Q\circ\nu-Q\circ\mu\|_{\mathrm{TV}}^2
\ \leq \
\frac14
\sum_{I\subseteq[n],\,I\neq\varnothing}
\delta_{|I|}^2
\prod_{i\in I}\chi^2\bigl(Q_i(\,\cdot\mid 1)\|Q_i(\,\cdot\mid 0)\bigr)\,.\]
\end{proposition}

Here, $|I|$ denotes the cardinality of the subset $I$ of $[n]$, and we use the convention $0\cdot(+\infty)=0$, as the $\chi^2$ divergences may be infinite. If $\mu$ and $\nu$ have matching moments through order $L$, then $\delta_\ell=0$ for $\ell\leq L$, so only subsets $I$ with $|I|>L$ contribute to the bound. As outlined in the last paragraph of Section~\ref{subsec:controlling-privatized-mixtures}, the proof relies on the affine dependence of each $Q_i(\cdot\mid t)$ on $t$, which makes the product $Q(\cdot\mid t)=\bigotimes_{i=1}^n Q_i(\cdot\mid t)$ polynomial in $t$. Taking the difference between the mixtures under $\nu$ and $\mu$ then replaces each power $t^\ell$ by the moment difference $\delta_\ell$. See the proof in Appendix~\ref{app:proof-mixtures-product-measures} for further details.


\subsubsection{$\chi^2$--$\mathrm{TV}^2$ Lipschitz property for 
$\alpha$-LDP
Markov kernels}

We first introduce a Lipschitz property for Markov kernels in the $\chi^2$--$\mathrm{TV}^2$ sense.

\begin{definition}
\label{def:tv-chi2-contracting-kernel}
For $\beta>0$, we say that a Markov kernel $K$ is
$\beta$-Lipschitz in the $\chi^2$--$\mathrm{TV}^2$ sense if, for any two
probability measures $P_0,P_1$ on its input space,
\[
\chi^2\left(K\circ P_1\,\middle\|\,K\circ P_0\right)
\ \leq \
\beta\|P_1-P_0\|_{\mathrm{TV}}^2.
\]
\end{definition}

The next lemma shows that local differential privacy implies the Lipschitz property of Definition~\ref{def:tv-chi2-contracting-kernel}.
To state it, let $\mathcal M_\alpha(\mathcal X)$ denote the class of $\alpha$-LDP Markov kernels from a measurable space $\mathcal X$ to an arbitrary measurable space, for $\alpha>0$. As in Section~\ref{subsec:ldp-observation-model}, this means that $K(A\mid x)\leq e^\alpha K(A\mid x')$ for all $x,x'\in\mathcal X$ and every measurable set $A$ in the output space of $K$.

\begin{lemma}
\label{lem:chi2-contraction-ldp}
Let $\mathcal X$ be a measurable space, let $\alpha>0$, and set $\beta_\alpha:=(e^\alpha-1)^2$.
Any $K\in\mathcal M_\alpha(\mathcal X)$ is $\beta_\alpha$-Lipschitz in the $\chi^2$--$\mathrm{TV}^2$ sense.
\end{lemma}

The Lipschitz constant $\beta_\alpha$ is of order $\alpha^2$ as $\alpha\downarrow0$. When $\alpha\in(0,1)$, this gives an indication of how privacy reduces the information available for distinguishing the input laws $P_0$ and $P_1$. Although it is only an upper bound on the optimal $\chi^2$--$\mathrm{TV}^2$ Lipschitz constant, it is sufficient for our statistical analysis.
The proof is given in Appendix~\ref{sec:proofs-main-tools}.

\begin{remark}
The lemma also yields a $\chi^2$ contraction bound with constant $\min\{1,\beta_\alpha/4\}$:
\[
\chi^2(K\circ P_1\|K\circ P_0)
\ \leq \ (1 \wedge \beta_\alpha/4)\,\chi^2(P_1\|P_0),
\]
where the constant $1$ comes from the data processing inequality, and $\beta_\alpha/4$ follows from Lemma~\ref{lem:chi2-contraction-ldp} and the standard inequality $\|P_1-P_0\|_{\mathrm{TV}}^2\leq\chi^2(P_1\|P_0)/4$.
However, this bound can be much weaker than the $\chi^2$--$\mathrm{TV}^2$ Lipschitz bound in Definition~\ref{def:tv-chi2-contracting-kernel}. For example, take $P_0=(1-\varepsilon,\varepsilon)$ and $P_1=(1,0)$ with $\varepsilon\in(0,1)$; then $\chi^2(P_1\|P_0)=\varepsilon/(1-\varepsilon)$, while $\|P_1-P_0\|_{\mathrm{TV}}^2=\varepsilon^2$.
\end{remark}

\subsection{Proof of Theorem~\ref{thm:private-mixture-comparison}}
\label{proo:general_tv_bound}

The proof proceeds in three steps. We first reduce the comparison to one-dimensional mixtures (Lemma~\ref{lem:coordinatewise-mixture-comparison}).
We then apply the moment-based TV bound for one-dimensional kernels (Proposition~\ref{lem:mixtures-product-measures}).
We control the resulting $\chi^2$ divergences using the $\chi^2$--$\mathrm{TV}^2$ Lipschitz bound under LDP (Lemma~\ref{lem:chi2-contraction-ldp}) and conclude the proof.

\begin{proof}[\textbf{Proof.}]
Fix a probability vector $\boldsymbol p\in\mathcal P_{m,b}^{\mathrm{pair}}$, a privacy mechanism $K_{1:n}\in\mathcal M_\alpha^n$, and two probability measures $\mu_0,\mu_1$ on $[0,1]$. For convenience, we introduce the following notation. For fixed $\boldsymbol t\in[0,1]^m$, write $Q(\cdot\mid\boldsymbol t):=P_{\boldsymbol p(\boldsymbol t),K_{1:n}}^{Z_{1:n}}$ for the joint distribution of $Z_{1:n}$, and $Q_i(\cdot\mid\boldsymbol t):=P_{\boldsymbol p(\boldsymbol t),K_{1:n}}^{Z_i}$ for the marginal distribution of $Z_i$, for $i\in[n]$. For $r\in\{0,1\}$, the mixture law can then be written as
\begin{equation}
\label{eq:private-mixture-kernel-notation}
M_{r,\boldsymbol p,K_{1:n}}^{Z_{1:n}}
:=\mathbb E_{\boldsymbol t\sim\mu_r^{\otimes m}}
\!\left[P_{\boldsymbol p(\boldsymbol t),K_{1:n}}^{Z_{1:n}}\right]
=Q\circ\mu_r^{\otimes m}.
\end{equation}

For  fixed $\boldsymbol t$, the  $Z_i$ are independent, so $Q(\cdot\mid\boldsymbol t)=\bigotimes_{i=1}^n Q_i(\cdot\mid\boldsymbol t)$, with $Q_i(\cdot\mid\boldsymbol t)=K_i\circ\boldsymbol p(\boldsymbol t)$.
These identities also hold for the induced one-dimensional kernels, which we define below.
For $j\in[m]$, denoting the fixed coordinates by $\boldsymbol t_{-j}=(t_k)_{k\ne j}\in[0,1]^{m-1}$,  write
\[
Q^{\boldsymbol t_{-j}}(\cdot\mid v)
 := Q(\cdot\mid t_1,\ldots,t_{j-1},v,t_{j+1},\ldots,t_m),
\qquad
Q_i^{\boldsymbol t_{-j}}(\cdot\mid v)
 := Q_i(\cdot\mid t_1,\ldots,t_{j-1},v,t_{j+1},\ldots,t_m)
\]
for $v\in[0,1]$ and $i\in[n]$. Similarly, write $\boldsymbol p^{\boldsymbol t_{-j}}(v):=\boldsymbol p(t_1,\ldots,t_{j-1},v,t_{j+1},\ldots,t_m)$ for the law of $X_i$ when $t_j=v$. Then, these one-dimensional kernels satisfy
\begin{equation}
\label{eq:private-mixture-output-law}
Q^{\boldsymbol t_{-j}}(\cdot\mid v)
\ =\ \bigotimes_{i=1}^n Q_i^{\boldsymbol t_{-j}}(\cdot\mid v),
\qquad \qquad
Q_i^{\boldsymbol t_{-j}}(\cdot\mid v)
\ =\ K_i\circ\boldsymbol p^{\boldsymbol t_{-j}}(v).
\end{equation}

\textit{Step 1: Reduction to one-dimensional kernels.}
Applying Lemma~\ref{lem:coordinatewise-mixture-comparison} to the kernel $Q$ and the product measures $\mu_0^{\otimes m}$ and $\mu_1^{\otimes m}$, and using~\eqref{eq:private-mixture-kernel-notation}, gives
\begin{equation}
\label{eq:private-mixture-coordinate-reduction}
\left\|M_{1,\boldsymbol p,K_{1:n}}^{Z_{1:n}}-M_{0,\boldsymbol p,K_{1:n}}^{Z_{1:n}}\right\|_{\mathrm{TV}}
\ \leq\ m\, \max_{j\in[m]} \ \sup_{\boldsymbol t_{-j}\in[0,1]^{m-1}}
\left\|Q^{\boldsymbol t_{-j}}\circ\mu_1-Q^{\boldsymbol t_{-j}}\circ\mu_0\right\|_{\mathrm{TV}}.
\end{equation}

\textit{Step 2: Moment-based TV bound.}
Fix $j\in[m]$ and $\boldsymbol t_{-j}\in[0,1]^{m-1}$.
Assume for the moment that
\begin{equation}
\label{eq:coordinatewise-affine-p}
\boldsymbol{p}^{\boldsymbol t_{-j}}(v)
=(1-v)\boldsymbol{p}^{\boldsymbol t_{-j}}(0)+v\,\boldsymbol{p}^{\boldsymbol t_{-j}}(1),
\qquad v\in[0,1].
\end{equation}
Then by~\eqref{eq:private-mixture-output-law} and the linearity of the composition operation $\circ$ with $K_i$, we obtain
\[
Q_i^{\boldsymbol t_{-j}}(\cdot\mid v)
=(1-v)Q_i^{\boldsymbol t_{-j}}(\cdot\mid0)
+vQ_i^{\boldsymbol t_{-j}}(\cdot\mid1).
\]
Thus, Proposition~\ref{lem:mixtures-product-measures} applies to $Q^{\boldsymbol t_{-j}}$, and yields
\begin{equation}
\label{eq:private-mixture-scalar-chi2}
\left\|
Q^{\boldsymbol t_{-j}}\circ\mu_1
-Q^{\boldsymbol t_{-j}}\circ\mu_0
\right\|^2_{\mathrm{TV}}
\ \leq \ \frac14
\sum_{I\subseteq[n],\,I\neq\varnothing}
\delta_{|I|}^2\prod_{i\in I}
\chi^2\Bigl(Q_i^{\boldsymbol t_{-j}}(\cdot\mid1)
\,\Big\|\,Q_i^{\boldsymbol t_{-j}}(\cdot\mid0)\Bigr)
,
\end{equation}
with $\delta_\ell:=\int_{[0,1]}t^\ell\,\mu_1(dt)-\int_{[0,1]}t^\ell\,\mu_0(dt)$ for $\ell\geq0$, the difference between the moments of  $\mu_1$ and $\mu_0$.

\noindent\textit{Proof of~\eqref{eq:coordinatewise-affine-p}.} Recall that, for the family $(\boldsymbol p(\boldsymbol t))_{\boldsymbol t\in[0,1]^m}$ associated with $\boldsymbol{p}\in\mathcal P_{m,b}^{\mathrm{pair}}$ in Definition~\ref{def:pairwise-perturbations},
\begin{equation*}
p_{2j-1}(\boldsymbol t)\ =\ b(1+t_j),
\qquad \qquad 
p_{2j}(\boldsymbol t)\ =\ b(1-t_j),
\end{equation*}
while $p_k(\boldsymbol t)$ does not depend on $t_j$ for $k\in[S]\setminus\{2j-1,2j\}$. Since
\[
\bigl(b(1+t_j),b(1-t_j)\bigr)
=(1-t_j)(b,b)+t_j(2b,0),
\]
we obtain
$\boldsymbol p(\boldsymbol t)
=(1-t_j)\boldsymbol{p}^{\boldsymbol t_{-j}}(0)
+t_j\boldsymbol{p}^{\boldsymbol t_{-j}}(1)$.
Taking $t_j=v$ proves~\eqref{eq:coordinatewise-affine-p}.

 \textit{Step 3: $\chi^2$--$\mathrm{TV}^2$ Lipschitz bound and conclusion.}
For any $j\in[m]$, $\boldsymbol t_{-j}\in[0,1]^{m-1}$, $i\in[n]$,
and  $v\in[0,1]$,
the distribution of $Z_i$ under the input law $\boldsymbol p^{\boldsymbol t_{-j}}(v)$ is $Q_i^{\boldsymbol t_{-j}}(\cdot\mid v)=K_i\circ\boldsymbol p^{\boldsymbol t_{-j}}(v)$. Since $K_i$ is $\alpha$-LDP, applying Lemma~\ref{lem:chi2-contraction-ldp} with $(K,P_0,P_1)=(K_i,\boldsymbol p^{\boldsymbol t_{-j}}(0),\boldsymbol p^{\boldsymbol t_{-j}}(1))$ gives
\[
\chi^2\Bigl(Q_i^{\boldsymbol t_{-j}}(\cdot\mid1)
\,\Big\|\,Q_i^{\boldsymbol t_{-j}}(\cdot\mid0)\Bigr)
\ \leq \ (e^\alpha-1)^2\left\|\boldsymbol p^{\boldsymbol t_{-j}}(1)-\boldsymbol p^{\boldsymbol t_{-j}}(0)\right\|_{\mathrm{TV}}^2.
\]
As seen in the proof of~\eqref{eq:coordinatewise-affine-p}, the laws $\boldsymbol p^{\boldsymbol t_{-j}}(0)$ and $\boldsymbol p^{\boldsymbol t_{-j}}(1)$ agree in all but coordinates $2j-1$ and $2j$, where their entries are $(b,b)$ and $(2b,0)$, respectively. Hence 
$\|\boldsymbol p^{\boldsymbol t_{-j}}(1)-\boldsymbol p^{\boldsymbol t_{-j}}(0)\|_{\mathrm{TV}}=\tfrac12(|2b-b|+|0-b|)=b$.
Combining this identity with the above $\chi^2$--$\mathrm{TV}^2$ bound, and substituting into~\eqref{eq:private-mixture-scalar-chi2}, we obtain
\begin{equation*}
\begin{aligned}
\left\|Q^{\boldsymbol t_{-j}}\circ\mu_1-Q^{\boldsymbol t_{-j}}\circ\mu_0\right\|_{\mathrm{TV}}^2
\ \leq \ \frac14
\sum_{I\subseteq[n],\,I\neq\varnothing}
\delta_{|I|}^2[b(e^\alpha-1)]^{2|I|}
\ =\ \frac14
\sum_{\ell\in[n]}\binom n\ell
\delta_\ell^2[b(e^\alpha-1)]^{2\ell},
\end{aligned}
\end{equation*}
where the equality follows by grouping the terms according to the cardinality $|I|$.
Applying~\eqref{eq:private-mixture-coordinate-reduction} and taking the supremum over $\boldsymbol{p}\in\mathcal P_{m,b}^{\mathrm{pair}}$ and $K_{1:n}\in\mathcal M_\alpha^n$ concludes the proof of Theorem~\ref{thm:private-mixture-comparison}.
\end{proof}


\section{Proof of Theorem~\ref{thm:orientation-lower-bound}}
\label{proof:orientation-lower-bound}

We first present the construction and functional separation (Section~\ref{subsec:orientation-construction-separation}), then state the KL bound, which is the key ingredient of the proof (Section~\ref{subsec:orientation-kl-bound}). We combine these results with Le Cam's inequality to prove Theorem~\ref{thm:orientation-lower-bound} (Section~\ref{subsec:orientation-theorem-proof}). The proof of the KL bound is deferred to Section~\ref{sec:proof:lemma-oreinted-KL}.

\subsection{Construction and functional separation}
\label{subsec:orientation-construction-separation}

Fix an integer $1\leq m\leq\lfloor S/2\rfloor$ and a real number $b\in(0,1/2]$ such that $\mathcal P_{m,b}^{\mathrm{pair}}\neq\varnothing$.
Fix a probability vector $\boldsymbol p\in\mathcal P_{m,b}^{\mathrm{pair}}$ and consider the induced family $\boldsymbol p(\boldsymbol t)$, indexed by $\boldsymbol t\in[-1,1]^m$, from Definition~\ref{def:pairwise-perturbations}. We choose a vector $\boldsymbol t$ whose coordinates have a constant amplitude $t\in[0,1]$; that is, we set $t_j=\omega_jt$ for all $j\in[m]$, where $\omega=(\omega_1,\ldots,\omega_m)\in\{-1,1\}^m$. By Definition~\ref{def:pairwise-perturbations},
\begin{equation}
\label{eq:constant-amplitude-pairs}
p_{2j-1}(t\omega)=b(1+t\omega_j),
\qquad
p_{2j}(t\omega)=b(1-t\omega_j),
\qquad j\in[m],
\end{equation}
while $p_k(t\omega)=p_k$ for $k>2m$.

As proved in the next lemma, the variations of $F_\varphi$ between the two probability vectors $\boldsymbol p(t\omega)$ and $\boldsymbol p(t'\omega)$ reduce to those of the pair profile, as expressed in~\eqref{eq:orientation-functional-separation}.
Crucially, this functional separation is the same for every $\omega$. Indeed, changing $\omega$ only swaps the two probabilities within each pair in~\eqref{eq:constant-amplitude-pairs}, leaving the functional $F_\varphi(\boldsymbol p(t\omega))$ unchanged. We can therefore focus on finding a sign vector $\omega$ for which the corresponding privatized observation laws are close (in KL divergence), without changing this separation; this is the purpose of the next subsection.

\begin{lemma}
\label{lem:orientation-functional-separation}
For every $t,t'\in[0,1]$ and every $\omega\in\{-1,1\}^m$,
we have
\begin{equation}
\label{eq:orientation-functional-separation}
\Delta
:=\left|F_\varphi(\boldsymbol p(t\omega))-F_\varphi(\boldsymbol p(t'\omega))\right|
=m\left|\varphi_b(t)-\varphi_b(t')\right|.
\end{equation}
In particular, the separation $\Delta$ does not depend on the sign vector $\omega$.
\end{lemma}

\begin{proof}
Recall that $F_\varphi(\boldsymbol p)=\sum_{k\in[S]}\varphi(p_k)$ for every $\boldsymbol p\in\mathcal P_S$. Then, by~\eqref{eq:constant-amplitude-pairs}, for every $\omega\in\{-1,1\}^m$, 
\[
F_\varphi(\boldsymbol p(t\omega))
=m\left\{\varphi\bigl(b(1+t)\bigr)+\varphi\bigl(b(1-t)\bigr)\right\}
+\sum_{k>2m}\varphi(p_k)
=m\varphi_b(t)+\sum_{k>2m}\varphi(p_k),
\]
using the definition $\varphi_b(t)=\varphi(b(1+t))+\varphi(b(1-t))$. Using the same identity with $t'$ in place of $t$, we subtract $F_\varphi(\boldsymbol p(t'\omega))$ from $F_\varphi(\boldsymbol p(t\omega))$ and easily obtain~\eqref{eq:orientation-functional-separation}.
\end{proof}

\subsection{KL divergence between privatized observation laws}
\label{subsec:orientation-kl-bound}

Fix  $\boldsymbol p\in\mathcal P_{m,b}^{\mathrm{pair}}$ and consider the family $\boldsymbol p(t\omega)$ constructed in Section~\ref{subsec:orientation-construction-separation}. Fix a privacy mechanism $K_{1:n}\in\mathcal M_\alpha^n$, and write $P_{\omega,t}^{Z_{1:n}}:=P_{\boldsymbol p(t\omega),K_{1:n}}^{Z_{1:n}}$ for the joint law of the privatized observations $Z_{1:n}$ under the input distribution $\boldsymbol p(t\omega)$ and the mechanism $K_{1:n}$. The next proposition shows that, for at least one orientation $\omega$, the two privatized observation laws have a small KL divergence.

\begin{proposition}
\label{prop:orientation-kl-bound-general}
For any $K_{1:n}\in\mathcal M_\alpha^n$ and any $t,t'\in[0,1]$, there exists $\omega\in\{-1,1\}^m$ such that
\[
\operatorname{KL}\!\left(
P_{\omega,t}^{Z_{1:n}}
\,\middle\|\,
P_{\omega,t'}^{Z_{1:n}}
\right)
\leq ne^\alpha(e^\alpha-1)^2mb^2(t-t')^2.
\]
\end{proposition}

This result is the key ingredient in the proof of Theorem~\ref{thm:orientation-lower-bound}. Interestingly, we do not require a single orientation $\omega$ to yield a small KL divergence uniformly over all privacy mechanisms $K_{1:n}$. Instead, the orientation $\omega$ may depend on the privacy mechanism $K_{1:n}$ and the amplitudes $t,t'$.

In the proof of Proposition~\ref{prop:orientation-kl-bound-general}, the trick is to use a probabilistic argument that reduces the problem to controlling the average $\mathbb E_\omega\!\left[\operatorname{KL}\!\left(P_{\omega,t}^{Z_{1:n}}\middle\|P_{\omega,t'}^{Z_{1:n}}\right)\right]$, where $\omega\overset{\mathrm{unif.}}{\sim}\{-1,1\}^m$. This uniform average is easier to work with to obtain sharp bounds. More generally, the probabilistic method is widely used in combinatorics to establish the existence of objects with desired properties by studying random constructions~\cite{AlonSpencer2016}. See Section~\ref{sec:proof:lemma-oreinted-KL} for the detailed proof.

To see the role of privacy in Proposition~\ref{prop:orientation-kl-bound-general}, let us compare with the non-private setting. For two input distributions $(P,P')$, independence in the non-private setting gives $\operatorname{KL}(P^{\otimes n}\|P'^{\otimes n})=n\operatorname{KL}(P\|P')$. Taking $P=\boldsymbol p(t\omega)$ and $P'=\boldsymbol p(t'\omega)$, this divergence typically scales as $nmb(t-t')^2$, regardless of the orientation $\omega$ (see Remark~\ref{rmk:non-private-scale-copy}). Under local privacy, however, Proposition~\ref{prop:orientation-kl-bound-general} guarantees that some $\omega$ satisfies $\operatorname{KL}(P_{\omega,t}^{Z_{1:n}}\|P_{\omega,t'}^{Z_{1:n}})\lesssim \alpha^2nmb^2(t-t')^2$. This reflects the loss of information due to privacy, through a reduction in effective sample size from $n$ to $\alpha^2bn$, in the regime  $\alpha<1$ of strong privacy.

\begin{remark}
\label{rmk:non-private-scale-copy}
\textit{Non-private scale.} Recall that, for $n$ independent observations drawn from $\boldsymbol p$, a small perturbation $\Delta_k$ of a probability $p_k$ contributes a term of order $n\Delta_k^2/p_k$ to the non-private KL divergence. In our construction $\mathcal P_{m,b}^{\mathrm{pair}}$, each of the $m$ pairs consists of the probabilities $b(1+t)$ and $b(1-t)$, which vary by $\pm b(t'-t)$ when $t$ is replaced by $t'$. When $t,t'$ stay bounded away from $1$, these probabilities are of order $b$, so the total non-private KL divergence is of order $nmb(t-t')^2$. 
\end{remark}

\subsection{Proof of Theorem~\ref{thm:orientation-lower-bound}}
\label{subsec:orientation-theorem-proof}

The proof of Theorem~\ref{thm:orientation-lower-bound}
uses the following version of Le Cam's two-point inequality
for quadratic loss, which follows from
\cite[Theorem~2.2]{Tsybakov2009}.

\begin{lemma}[Le Cam's two-point inequality]
\label{lem:le-cam-two-point}
Let $P_0$ and $P_1$ be two probability measures
on the same measurable space, and let
$\theta_0,\theta_1\in\mathbb R$.
For every real-valued estimator $\widehat\theta$,
\[
\max_{r\in\{0,1\}}
\mathbb E_{P_r}\bigl[(\widehat\theta-\theta_r)^2\bigr]
\ \geq\
\frac{(\theta_1-\theta_0)^2}{8}
\left(1-\|P_1-P_0\|_{\mathrm{TV}}\right).
\]
\end{lemma}

To prove Theorem~\ref{thm:orientation-lower-bound}, we now apply this Le Cam inequality together with the above functional separation (Lemma~\ref{lem:orientation-functional-separation}) and KL bound (Proposition~\ref{prop:orientation-kl-bound-general}).

\begin{proof}[Proof of Theorem~\ref{thm:orientation-lower-bound}]
Fix $\boldsymbol p\in\mathcal P_{m,b}^{\mathrm{pair}}$ and consider the family $\boldsymbol p(t\omega)$ constructed in Section~\ref{subsec:orientation-construction-separation}, where $t\in[0,1]$ and $\omega\in\{-1,1\}^m$. Let $\alpha\in(0,1)$ and fix a privacy mechanism $K_{1:n}\in\mathcal M_\alpha^n$. Recall that $P_{\omega,t}^{Z_{1:n}}$ denotes the joint law of the privatized observations $Z_{1:n}$ under the input distribution $\boldsymbol p(t\omega)$ and the mechanism $K_{1:n}$. For any $t,t'\in[0,1]$ satisfying $|t-t'|\leq\delta_{m,b}:=1\wedge(16e(e-1)^2mb^2n\alpha^2)^{-1/2}$, Proposition~\ref{prop:orientation-kl-bound-general} gives an orientation $\omega\in\{-1,1\}^m$ such that
\[
\operatorname{KL}\!\left(P_{\omega,t}^{Z_{1:n}}\middle\|P_{\omega,t'}^{Z_{1:n}}\right)
\leq ne^\alpha(e^\alpha-1)^2mb^2(t-t')^2
\leq e(e-1)^2n\alpha^2mb^2\delta_{m,b}^2
\leq\frac1{16},
\]
where we used $\alpha<1$ and the definition of $\delta_{m,b}$. 
Now, fix $(t,t',\omega)$ satisfying the above bound.
Pinsker's inequality then gives
\[
\left\|P_{\omega,t}^{Z_{1:n}}-P_{\omega,t'}^{Z_{1:n}}\right\|_{\mathrm{TV}}
\leq
\sqrt{\frac12\operatorname{KL}\!\left(P_{\omega,t}^{Z_{1:n}}\middle\|P_{\omega,t'}^{Z_{1:n}}\right)}
\leq\frac12.
\]
Let $\Delta:=|F_\varphi(\boldsymbol p(t\omega))-F_\varphi(\boldsymbol p(t'\omega))|$ denote the functional separation. Let $\widehat F$ be any estimator based on the privatized observations $Z_{1:n}$. Since $\boldsymbol p(t\omega),\boldsymbol p(t'\omega)\in\mathcal P_{m,b}^{\mathrm{pair}}$, Lemma~\ref{lem:le-cam-two-point} yields
\[
\sup_{\boldsymbol q\in\mathcal P_{m,b}^{\mathrm{pair}}}
\mathbb E_{\boldsymbol q,K_{1:n}}\!\left[(\widehat F-F_\varphi(\boldsymbol q))^2\right]
\geq
\frac{\Delta^2}{8}
\left(1-\left\|P_{\omega,t}^{Z_{1:n}}-P_{\omega,t'}^{Z_{1:n}}\right\|_{\mathrm{TV}}\right)
\geq\frac{\Delta^2}{16}
=\frac{m^2}{16}|\varphi_b(t)-\varphi_b(t')|^2,
\]
where the equality follows from Lemma~\ref{lem:orientation-functional-separation}. The lower bound in Theorem~\ref{thm:orientation-lower-bound} follows by taking the supremum over $t,t'\in[0,1]$ satisfying $|t-t'|\leq\delta_{m,b}$, and then the infimum over all mechanisms $K_{1:n}\in\mathcal M_\alpha^n$ and all estimators $\widehat F$.
\end{proof}


\subsection{Proof of Proposition~\ref{prop:orientation-kl-bound-general}}
\label{sec:proof:lemma-oreinted-KL}

The proof of Proposition~\ref{prop:orientation-kl-bound-general} uses the following lemma, which gives expressions and bounds for the densities of the privatized observations.
Its proof is deferred to  Appendix~\ref{app:proof-orientation-densities}. Fix a probability vector $\boldsymbol p=(p_x)_{x\in[S]}\in\mathcal P_{m,b}^{\mathrm{pair}}$, and
write $\delta_x$ for the Dirac probability measure at $x$.

\begin{lemma}
\label{lem:orientation-densities}
Let $K\in\mathcal M_\alpha([S])$ and fix $t\in[0,1]$. Consider the family $\{\boldsymbol p(t\omega)\}_{\omega\in\{-1,1\}^m}$ constructed in Section~\ref{subsec:orientation-construction-separation}. For every $x\in[S]$, we have $K\circ\delta_x\ll K\circ\boldsymbol p(\boldsymbol0)$ and
\begin{equation}
\label{eq:orientation-density-difference}
\frac{d(K\circ\delta_x)}{d(K\circ\boldsymbol p(\boldsymbol0))}
\ \geq \ e^{-\alpha},
\qquad \qquad
\left|\frac{d(K\circ\delta_x-K\circ\delta_{x'})}{d(K\circ\boldsymbol p(\boldsymbol0))}\right|
\ \leq \ e^\alpha-1,
\qquad x,x'\in[S].
\end{equation}
For every $\omega\in\{-1,1\}^m$, we have $K\circ\boldsymbol p(t\omega)\ll K\circ\boldsymbol p(\boldsymbol0)$, with
\begin{equation}
\label{eq:orientation-output-density}
\frac{d(K\circ\boldsymbol p(t\omega))}{d(K\circ\boldsymbol p(\boldsymbol0))}
\ = \ 1+bt\sum_{j\in[m]}\omega_j
\frac{d(K\circ\delta_{2j-1}-K\circ\delta_{2j})}{d(K\circ\boldsymbol p(\boldsymbol0))}
\ \geq \ e^{-\alpha}.
\end{equation}
All these relations hold simultaneously for all $x,x'\in[S]$, $\omega\in\{-1,1\}^m$, outside a  $K\circ\boldsymbol p(\boldsymbol0)$-null set.
\end{lemma}

\begin{remark}\label{remark-finiteness-input-set-contraction-privacy}
The finiteness of the input space $[S]$ avoids a technical difficulty with conditional densities. For each $x\in[S]$, let $q_x:=\frac{d(K\circ\delta_x)}{d(K\circ\boldsymbol p(\boldsymbol 0))}$ denote the conditional density of the privatized observation given input $x$, with respect to $K\circ\boldsymbol p(\boldsymbol 0)$. For each pair $(x,x')\in[S]^2$, the LDP condition gives $q_x\leq e^\alpha q_{x'}$ outside a $K\circ\boldsymbol p(\boldsymbol 0)$-null set. Since there are finitely many pairs $(x,x')\in[S]^2$, these inequalities hold simultaneously outside a single $K\circ\boldsymbol p(\boldsymbol 0)$-null set. The same finite-union argument applies to \eqref{eq:orientation-density-difference} over all pairs $(x,x')\in[S]^2$, and to \eqref{eq:orientation-output-density} over all orientations $\omega\in\{-1,1\}^m$, with $t$ fixed.

\noindent Although this simultaneous validity plays no key role in our proof of Theorem~\ref{thm:orientation-lower-bound}, it can be important in other problems involving conditional densities on uncountable input spaces. Indeed, on such spaces, the exceptional sets cannot in general be combined into a single null set. As discussed in \cite[Section 4.2.1]{IssartelRoueff}, the kernel privacy condition does not automatically give a conditional density inequality valid simultaneously for all input values.
\end{remark}

\begin{proof}[\textbf{Proof of Proposition~\ref{prop:orientation-kl-bound-general}}]
Fix $K_{1:n}\in\mathcal M_\alpha^n$ and $t,t'\in[0,1]$.
    Recall that $P_{\omega,t}^{Z_{1:n}}:=P_{\boldsymbol p(t\omega),K_{1:n}}^{Z_{1:n}}$ denotes
    the joint law of the privatized observations $Z_{1:n}=(Z_1,\ldots,Z_n)$. Let
    $P_{\omega,t}^{Z_i}$ denote its $i$th marginal, and use the same notation
    with $t'$ in place of $t$.
    For each fixed $\omega\in\{-1,1\}^m$, independence gives
    \begin{equation}
    \label{eq:orientation-kl-additivity}
    \operatorname{KL}
    \bigl(P_{\omega,t}^{Z_{1:n}}\,\|\,P_{\omega,t'}^{Z_{1:n}}\bigr)
    \ = \
    \sum_{i\in[n]}
    \operatorname{KL}
    \bigl(P_{\omega,t}^{Z_i}\,\|\,P_{\omega,t'}^{Z_i}\bigr).
    \end{equation}
For each $i\in[n]$, we have $P_{\omega,t}^{Z_i}=K_i\circ\boldsymbol p(t\omega)$, and likewise with $t'$ in place of $t$. Denote by $\rho_{i,\omega,t}$ and $\rho_{i,\omega,t'}$ the densities of these two laws with respect to $K_i\circ\boldsymbol p(\boldsymbol 0)$, which exist by Lemma~\ref{lem:orientation-densities}. We bound each term in \eqref{eq:orientation-kl-additivity} using the standard inequality $\operatorname{KL}\leq\chi^2$ (see, e.g., the proof of Lemma~3 in \cite{IssartelRoueff}):
\[
\begin{aligned}
\operatorname{KL}\bigl(P_{\omega,t}^{Z_i}\|P_{\omega,t'}^{Z_i}\bigr)
\leq\chi^2\bigl(P_{\omega,t}^{Z_i}\|P_{\omega,t'}^{Z_i}\bigr)
&=\int\frac{(\rho_{i,\omega,t}-\rho_{i,\omega,t'})^2}
{\rho_{i,\omega,t'}}\,d(K_i\circ\boldsymbol p(\boldsymbol 0))\\
&\leq e^\alpha b^2(t-t')^2
\int\biggl(\sum_{j\in[m]}\omega_j
\frac{d(K_i\circ\delta_{2j-1}-K_i\circ\delta_{2j})}{d(K_i\circ\boldsymbol p(\boldsymbol 0))}\biggr)^2
\,d(K_i\circ\boldsymbol p(\boldsymbol 0)),
\end{aligned}
\]
where the second inequality follows from \eqref{eq:orientation-output-density}, applied with $K=K_i$, which gives the expression for $\rho_{i,\omega,t}-\rho_{i,\omega,t'}$ and the lower bound $\rho_{i,\omega,t'}\geq e^{-\alpha}$.
    We now average this KL bound over $\omega\overset{\mathrm{unif.}}{\sim}\{-1,1\}^m$.
    Writing $\mathbb E_\omega$ for this average, we have
    $\mathbb E_\omega[\omega_j\omega_k]=0$ for $j\neq k$
    and $\mathbb E_\omega[\omega_j^2]=1$.
    Consequently, writing $q_{i,k}:=\frac{d(K_i\circ\delta_k)}{d(K_i\circ\boldsymbol p(\boldsymbol 0))}$ for each $i\in[n]$ and $k\in[S]$, we obtain
\[\mathbb E_\omega\,
    \bigg{(}\sum_{j\in[m]}
    \omega_j(q_{i,2j-1}-q_{i,2j})\bigg{)}^2
    \ = \
    \sum_{j\in[m]}(q_{i,2j-1}-q_{i,2j})^2.
    \]
    Taking $\mathbb E_\omega$ in
    \eqref{eq:orientation-kl-additivity} and using the last two
    displays yields
    \[
    \mathbb E_\omega\!
    \left[
    \operatorname{KL}
    \bigl(P_{\omega,t}^{Z_{1:n}}\,\|\,P_{\omega,t'}^{Z_{1:n}}\bigr)
    \right]
    \ \leq \
    e^\alpha b^2(t-t')^2
    \sum_{i\in[n]}
    \int \sum_{j\in[m]}
    (q_{i,2j-1}-q_{i,2j})^2\,\mathrm d(K_i\circ\boldsymbol p(\boldsymbol 0)) \, .
    \]
    Using \eqref{eq:orientation-density-difference}, which gives $\lvert q_{i,2j-1}-q_{i,2j}\rvert\leq e^\alpha-1$ almost everywhere with respect to $K_i\circ\boldsymbol p(\boldsymbol 0)$, we obtain
    \[
    \mathbb E_\omega\!
    \left[
    \operatorname{KL}
    \bigl(P_{\omega,t}^{Z_{1:n}}\,\|\,P_{\omega,t'}^{Z_{1:n}}\bigr)
    \right]
    \ \leq \
    ne^\alpha(e^\alpha-1)^2mb^2(t-t')^2,
    \]
    Hence, there exists an orientation $\omega\in\{-1,1\}^m$ for which this KL divergence satisfies the same bound. This proves Proposition~\ref{prop:orientation-kl-bound-general}.
\end{proof}


\paragraph{AI use statement.} 
We started this project and obtained the results well before we began using generative AI tools. 
We used these tools to streamline the writing, review the proofs, and check the calculations and notation. We take full responsibility for the final manuscript.


{\footnotesize
\let\originalbibitem\bibitem
\renewcommand{\bibitem}[2][]{%
  \def\currentbibkey{#2}\def\shortbibkey{BH}%
  \ifx\currentbibkey\shortbibkey
    \originalbibitem[BBG+22]{#2}%
  \else
    \originalbibitem[#1]{#2}%
  \fi
}
\bibliography{biblio}
\bibliographystyle{alpha}
}

\appendix
\numberwithin{theorem}{section}
\numberwithin{proposition}{section}
\numberwithin{lemma}{section}
\numberwithin{definition}{section}
\numberwithin{corollary}{section}
\numberwithin{example}{section}
\numberwithin{remark}{section}
\addtocontents{toc}{\protect\setcounter{tocdepth}{1}}

\section{Proofs of Lemmas~\ref{lem:product-measure-concentration}, \ref{lem:coordinatewise-mixture-comparison}, \ref{lem:chi2-contraction-ldp} and~\ref{lem:orientation-densities}}
\label{sec:proofs-main-tools}

\subsection{Proof of Lemma~\ref{lem:product-measure-concentration}}

For any square-integrable random variable $X$ and any constant $c$,
\[
\mathbb E[(X-c)^2]
=\operatorname{Var}(X)+(\mathbb EX-c)^2
\geq\operatorname{Var}(X).
\]
For each $j$, take $c_j=(\sup_{u\in[0,1]}F_j(u)+\inf_{u\in[0,1]}F_j(u))/2$.
Since $|F_j(U_j)-c_j|\leq V(F_j)/2$, we obtain
\[
\operatorname{Var}(F_j(U_j))
\leq\mathbb E[(F_j(U_j)-c_j)^2]
\leq\frac{V(F_j)^2}{4}.
\]
Independence gives $\operatorname{Var}(F(U))=\sum_{j\in[m]}\operatorname{Var}(F_j(U_j))$.
Chebyshev's inequality then yields
\[
\mathbb P\left(|F(U)-\mathbb EF(U)|\geq\Delta/4\right)
\leq\frac{16\operatorname{Var}(F(U))}{\Delta^2}
\leq\frac{4\sum_{j\in[m]}V(F_j)^2}{\Delta^2},
\]
which concludes the proof of Lemma~\ref{lem:product-measure-concentration}.


\subsection{Proof of Lemma~\ref{lem:coordinatewise-mixture-comparison}}

To compare $Q\circ\mu$ and $Q\circ\nu$, we replace the marginals $\mu_j$ by $\nu_j$ one at a time.
For $j=0,\ldots,m$, define the mixture law
\[
M^{(j)}
\ :=\ Q\circ
\left[
\left(\bigotimes_{k=1}^{j}\nu_k\right)
\otimes
\left(\bigotimes_{k=j+1}^{m}\mu_k\right)
\right].
\]
Thus $M^{(0)}=Q\circ\mu$ and $M^{(m)}=Q\circ\nu$.
Fix $j\in[m]$. The mixing measure of $M^{(j)}$ is obtained from that of $M^{(j-1)}$ by replacing $\mu_j$ with $\nu_j$, and keeping all other marginals fixed.
Thus, for every measurable set $A$ in the output space of $Q$, Fubini's theorem gives
\[
M^{(j)}(A)-M^{(j-1)}(A)
\ =\ \int
\Bigl\{
(Q^{\boldsymbol t_{-j}}\circ\nu_j)(A)
-(Q^{\boldsymbol t_{-j}}\circ\mu_j)(A)
\Bigr\}
\prod_{k<j}\nu_k(dt_k)
\prod_{k>j}\mu_k(dt_k).
\]
For every fixed $\boldsymbol t_{-j}$, the integrand is bounded in absolute value by $\|Q^{\boldsymbol t_{-j}}\circ\nu_j-Q^{\boldsymbol t_{-j}}\circ\mu_j\|_{\mathrm{TV}}$, which by assumption is at most $\varepsilon_j$.
Integrating over the remaining coordinates therefore gives $\lvert M^{(j)}(A)-M^{(j-1)}(A)\rvert\leq\varepsilon_j$.
Taking the supremum over measurable sets $A$ gives
$\left\|M^{(j)}-M^{(j-1)}\right\|_{\mathrm{TV}}
\leq\varepsilon_j$.
Finally, the triangle inequality implies
\[
\left\|M^{(m)}-M^{(0)}\right\|_{\mathrm{TV}}
\leq
\sum_{j=1}^m
\left\|M^{(j)}-M^{(j-1)}\right\|_{\mathrm{TV}}
\leq
\sum_{j=1}^m\varepsilon_j,
\]
which proves Lemma~\ref{lem:coordinatewise-mixture-comparison}.


\subsection{Proof of Lemma~\ref{lem:chi2-contraction-ldp}}

Before proving Lemma~\ref{lem:chi2-contraction-ldp}, we establish the following auxiliary lemma. It gives a uniform bound on the difference between two densities of the privatized observation.

\begin{lemma}
\label{lem:ldp-output-domination}
Let $K\in\mathcal M_\alpha(\mathcal X)$, and denote its output space by $\mathcal Z$. For any probability measures $\mu_0$, $\mu_1$, and $\mu_2$ on $\mathcal X$, we have $K\circ\mu_i\leq e^\alpha(K\circ\mu_j)$ for all $i,j\in\{0,1,2\}$. In particular, we have $K\circ\mu_i\ll K\circ\mu_0$ for $i\in[2]$, and
\[
\left|\frac{d(K\circ\mu_1-K\circ\mu_2)}{d(K\circ\mu_0)}\right|
\leq e^\alpha-1,
\qquad K\circ\mu_0\text{-almost everywhere}.
\]
\end{lemma}

\begin{proof}
Fix a measurable set $A\subseteq\mathcal Z$. Since $K$ is $\alpha$-LDP, we have $K(A\mid x)\leq e^\alpha K(A\mid x')$ for all $x,x'\in\mathcal X$. Integrating over $x'$ with respect to $\mu_j$, and over $x$ with respect to $\mu_i$, gives $(K\circ\mu_i)(A)\leq e^\alpha(K\circ\mu_j)(A)$ for all $i,j\in\{0,1,2\}$. In particular, $K\circ\mu_i\ll K\circ\mu_0$,  so  the Radon-Nikodym derivative $d(K\circ\mu_i)/d(K\circ\mu_0)$ exists for $i\in[2]$.

Setting $\nu_i:=K\circ\mu_i$ for $i\in\{0,1,2\}$, and using  $\nu_i(A)\leq e^\alpha \nu_j(A)$ for all $i,j\in\{0,1,2\}$, we obtain
\[
\nu_1(A)-\nu_2(A)
\leq (1-e^{-\alpha})\nu_1(A)
\leq (e^\alpha-1)\nu_0(A).
\]
Exchanging $\nu_1$ and $\nu_2$ yields $\lvert(\nu_1-\nu_2)(A)\rvert\leq(e^\alpha-1)\nu_0(A)$. Since this holds for every measurable set $A$, the Radon-Nikodym derivative $d(\nu_1-\nu_2)/d\nu_0$ satisfies the bound stated in the lemma.
\end{proof}

\begin{proof}[\textbf{Proof of Lemma~\ref{lem:chi2-contraction-ldp}.}]
Fix $K\in\mathcal M_\alpha(\mathcal X)$ and two probability measures $P_0$ and $P_1$ on $\mathcal X$. By Lemma~\ref{lem:ldp-output-domination}, $K\circ P_1\ll K\circ P_0$, so the Radon--Nikodym derivative of $K\circ P_1$ with respect to $K\circ P_0$ exists. Thus
\begin{equation}
\label{eq:ldp-chi2-integral}
\chi^2(K\circ P_1\|K\circ P_0)
=\int\left(\frac{d(K\circ P_1-K\circ P_0)}{d(K\circ P_0)}\right)^2\,d(K\circ P_0).
\end{equation}
Applying the bound in Lemma~\ref{lem:ldp-output-domination} with $(\mu_0,\mu_1,\mu_2)=(P_0,P_1,P_0)$ would give $\chi^2(K\circ P_1\|K\circ P_0)\leq(e^\alpha-1)^2$, which is insufficient to prove Lemma~\ref{lem:chi2-contraction-ldp}.

Set $\delta:=\|P_1-P_0\|_{\mathrm{TV}}$. If $\delta=0$, then $P_1=P_0$ and $\chi^2(K\circ P_1\|K\circ P_0)=0$, so the inequality in Lemma~\ref{lem:chi2-contraction-ldp} holds. Assume henceforth that $\delta>0$. Set $\mu:=(P_1-P_0)/\delta$, and let $\mu=\mu_+-\mu_-$ be its Hahn-Jordan decomposition. Since $\mu(\mathcal X)=0$, we have $\mu_+(\mathcal X)=\mu_-(\mathcal X)=\sup_{A\text{ measurable}}|\mu(A)|=1$. Thus $\mu_+$, $\mu_-$ are probability measures, so Lemma~\ref{lem:ldp-output-domination} applies with $(\mu_0,\mu_1,\mu_2)=(P_0,\mu_+,\mu_-)$ and gives
\[
\left|\frac{d(K\circ P_1-K\circ P_0)}{d(K\circ P_0)}\right|
=\delta\left|\frac{d(K\circ\mu_+-K\circ\mu_-)}{d(K\circ P_0)}\right|
\leq\delta(e^\alpha-1),
\qquad K\circ P_0\text{-almost everywhere}.
\]
Substituting this bound into \eqref{eq:ldp-chi2-integral} gives $\chi^2(K\circ P_1\|K\circ P_0)\leq(e^\alpha-1)^2\delta^2$, which proves Lemma~\ref{lem:chi2-contraction-ldp}.
\end{proof}

\subsection{Proof of Lemma~\ref{lem:orientation-densities}}
\label{app:proof-orientation-densities}

For any $x,x'\in[S]$, applying Lemma~\ref{lem:ldp-output-domination} with $(\mu_0,\mu_1,\mu_2)=(\boldsymbol p(\boldsymbol 0),\delta_x,\delta_{x'})$ gives $K\circ\delta_x\ll K\circ\boldsymbol p(\boldsymbol 0)$ and
\[
\left|\frac{d(K\circ\delta_x-K\circ\delta_{x'})}{d(K\circ\boldsymbol p(\boldsymbol 0))}\right|
\leq e^\alpha-1.
\]
For any $\omega\in\{-1,1\}^m$, the construction in Section~\ref{subsec:orientation-construction-separation} gives $\boldsymbol p(t\omega)=\boldsymbol p(\boldsymbol 0)+bt\sum_{j\in[m]}\omega_j(\delta_{2j-1}-\delta_{2j})$. Composing with $K$, we obtain
\[
K\circ\boldsymbol p(t\omega)
=K\circ\boldsymbol p(\boldsymbol 0)
+bt\sum_{j\in[m]}\omega_j
\bigl(K\circ\delta_{2j-1}-K\circ\delta_{2j}\bigr).
\]
By Lemma~\ref{lem:ldp-output-domination}, we have $K\circ\boldsymbol p(t\omega)\ll K\circ\boldsymbol p(\boldsymbol 0)$, so we can take Radon--Nikodym derivatives with respect to $K\circ\boldsymbol p(\boldsymbol 0)$ to obtain
\[
\frac{d(K\circ\boldsymbol p(t\omega))}{d(K\circ\boldsymbol p(\boldsymbol 0))}
=1+bt\sum_{j\in[m]}\omega_j
\frac{d(K\circ\delta_{2j-1}-K\circ\delta_{2j})}{d(K\circ\boldsymbol p(\boldsymbol 0))}.
\]
Finally, Lemma~\ref{lem:ldp-output-domination} gives $K\circ\boldsymbol p(\boldsymbol 0)\leq e^\alpha(K\circ\mu)$ for any probability measure $\mu$ on $[S]$. Applying this inequality with $\mu=\delta_x$ and with $\mu=\boldsymbol p(t\omega)$ yields
\[
\frac{d(K\circ\delta_x)}{d(K\circ\boldsymbol p(\boldsymbol 0))}
\geq e^{-\alpha},
\qquad
\frac{d(K\circ\boldsymbol p(t\omega))}{d(K\circ\boldsymbol p(\boldsymbol 0))}
\geq e^{-\alpha}.
\]
Since $t$ is fixed and the sets $[S]$ and $\{-1,1\}^m$ are finite, all these identities and inequalities hold simultaneously outside a single $K\circ\boldsymbol p(\boldsymbol 0)$-null set. This proves Lemma~\ref{lem:orientation-densities}.


\section{Proof of Proposition~\ref{lem:mixtures-product-measures}}
\label{app:proof-mixtures-product-measures}

We will use the  following three  elementary lemmas.
We omit the proof of the first lemma.

\begin{lemma}[Product expansion]
\label{lem:product-expansion}
For any real $t$ and any real coefficients $a_i$ and $b_i$ indexed by $i\in[n]$, we have $\prod_{i=1}^n(a_i+t b_i)=\sum_{I\subseteq[n]}t^{|I|}\prod_{i\in I}b_i\prod_{i\notin I}a_i$.
\end{lemma}

Recall that the total variation distance between two probability measures $P_1$ and $P_2$ is defined by $\|P_1-P_2\|_{\mathrm{TV}}:=\sup_{A\text{ measurable}}|P_1(A)-P_2(A)|$.

\begin{lemma}[Total variation in terms of densities]
\label{lem:total-variation-densities}
Let $P_1$, $P_2$ and $R$ be probability measures on the same measurable space. Suppose that $P_1-P_2\ll R$. Then
\[
\|P_1-P_2\|_{\mathrm{TV}}
\ =\ \frac12\int\left|\frac{d(P_1-P_2)}{dR}\right|\,dR
\ \leq\ \frac12\left(\int\left(\frac{d(P_1-P_2)}{dR}\right)^2dR\right)^{1/2}.
\]
\end{lemma}

\begin{proof}
Write $h=d(P_1-P_2)/dR$. Since $\int h\,dR=1-1=0$, the positive and negative parts of $h$ have the same integral. Hence $\|P_1-P_2\|_{\mathrm{TV}}=\frac12\int|h|\,dR$. The inequality follows from Cauchy--Schwarz.
\end{proof}

Recall that the $\chi^2$ divergence between two probability measures $P$ and $R$ is defined by $\chi^2(P\|R):=\int\left(\frac{dP}{dR}-1\right)^2dR$ if $P\ll R$, and by $\chi^2(P\|R):=+\infty$ otherwise.

\begin{lemma}[Orthogonality of products of centered densities]
\label{lem:centered-density-products}
For each $i\in[n]$, let $P_i$ and $R_i$ be probability measures on the same measurable space, and suppose that $\chi^2(P_i\|R_i)<\infty$. Then, for any real coefficients $c_I$ indexed by nonempty subsets $I\subseteq[n]$,
\[
\int\left(\sum_{\varnothing\neq I\subseteq[n]}c_I\prod_{i\in I}\left(\frac{dP_i}{dR_i}(z_i)-1\right)\right)^2 \bigotimes_{i=1}^n R_i(dz_i)
\ =\ \sum_{\varnothing\neq I\subseteq[n]}c_I^2\prod_{i\in I}\chi^2(P_i\|R_i).
\]
\end{lemma}

\begin{proof}
Write $h_i=dP_i/dR_i$. We have $\int(h_i-1)\,dR_i=0$ and $\int(h_i-1)^2\,dR_i=\chi^2(P_i\|R_i)<\infty$. Expanding the square gives a sum over pairs of nonempty subsets $I,I'\subseteq[n]$:
\[
\begin{aligned}
\int\left(\prod_{i\in I}(h_i(z_i)-1)\right)\left(\prod_{j\in I'}(h_j(z_j)-1)\right)\bigotimes_{i=1}^n R_i(dz_i)
&\ =\ \prod_{i\in I\cap I'}\int(h_i-1)^2\,dR_i
       \prod_{j\in I\triangle I'}\int(h_j-1)\,dR_j\\
&\ =\ \begin{cases}
0, & \text{if } I\neq I',\\
\displaystyle\prod_{i\in I}\chi^2(P_i\|R_i), & \text{if } I=I',
\end{cases}
\end{aligned}
\]
where $I\triangle I':=(I\setminus I')\cup(I'\setminus I)$ denotes the symmetric difference of $I$ and $I'$.
Thus, only the diagonal terms remain, each with coefficient $c_I^2$, which gives the stated identity.
\end{proof}

\subsection{Proof of the proposition}

Recall the setting of Proposition~\ref{lem:mixtures-product-measures}. Let $\mu$ and $\nu$ be probability measures on $[0,1]$.
 For $\ell\geq0$, write $\delta_\ell:=\int_{[0,1]}t^\ell\,(\nu-\mu)(dt)$.
 For every $t\in[0,1]$, the kernel $Q$ satisfies $Q(\cdot\mid t)=\bigotimes_{i=1}^n Q_i(\cdot\mid t)$, where each $Q_i$ is a Markov kernel from $[0,1]$ to a measurable space $\mathcal Z_i$ and satisfies $Q_i(\cdot\mid t)=(1-t)Q_i(\cdot\mid0)+tQ_i(\cdot\mid1)$.
Let $J$ denote the set of indices at which the two endpoint laws differ,
\[
J\ :=\ \bigl\{i\in[n]:Q_i(\cdot\mid1)\neq Q_i(\cdot\mid0)\bigr\}.
\]
The next lemma expresses the difference between the two mixtures in terms of the moment differences $\delta_\ell$ and the differences between the endpoint laws.

\begin{lemma}
\label{lem:mixture-difference-expansion}
For every measurable rectangle $A=\prod_{i=1}^n A_i$, we have
\begin{equation}
\label{eq:mixture-difference-expansion}
(Q\circ\nu)(A)-(Q\circ\mu)(A)
\ =\ \sum_{\varnothing\neq I\subseteq J}\delta_{|I|}
\prod_{i\in I}\bigl[Q_i(A_i\mid1)-Q_i(A_i\mid0)\bigr]
\prod_{i\in[n]\setminus I}Q_i(A_i\mid0).
\end{equation}
In particular, we have $Q\circ\mu=Q\circ\nu$ whenever $\delta_\ell=0$ for all $1\leq\ell\leq |J|$.
\end{lemma}

\begin{proof}
Using the product form of $Q(\cdot\mid t)$ and the affine dependence of each $Q_i(\cdot\mid t)$ on $t$,
we obtain, for every measurable rectangle $A=\prod_{i=1}^n A_i$,
\[
\begin{aligned}
(Q\circ\nu)(A)-(Q\circ\mu)(A)
&\ =\ \int_{[0,1]}\prod_{i=1}^n Q_i(A_i\mid t)\,(\nu-\mu)(dt)\\
&\ =\ \int_{[0,1]}\prod_{i=1}^n
\Bigl[Q_i(A_i\mid0)+t\bigl(Q_i(A_i\mid1)-Q_i(A_i\mid0)\bigr)\Bigr]\,(\nu-\mu)(dt)\\
&\ =\ \sum_{I\subseteq[n]}
\left(\int_{[0,1]}t^{|I|}\,(\nu-\mu)(dt)\right)
\prod_{i\in I}\bigl(Q_i(A_i\mid1)-Q_i(A_i\mid0)\bigr)
\prod_{i\notin I}Q_i(A_i\mid0),
\end{aligned}
\]
where the last line follows from Lemma~\ref{lem:product-expansion}.
The integral in the last line equals $\delta_{|I|}$. Since $\delta_0=\nu([0,1])-\mu([0,1])=0$, the summand indexed by $I=\varnothing$ vanishes. For each $I\not\subseteq J$, there exists an index $i\in I\setminus J$ for which $Q_i(A_i\mid1)-Q_i(A_i\mid0)=0$ by definition of $J$, so the summand indexed by $I$ also vanishes. This proves~\eqref{eq:mixture-difference-expansion}.
In particular, if $\delta_\ell=0$ for all $1\leq\ell\leq|J|$, then~\eqref{eq:mixture-difference-expansion} gives equality of the two mixtures on all measurable rectangles. Since these rectangles form a $\pi$-system generating the product $\sigma$-algebra, it follows that $Q\circ\nu=Q\circ\mu$.
\end{proof}

When $Q_i(\cdot\mid1)$ has a density with respect to $Q_i(\cdot\mid0)$, Lemma~\ref{lem:mixture-difference-expansion} gives the following formula.

\begin{lemma}
\label{lem:mixture-difference-density}
Suppose that $Q_i(\cdot\mid1)$ has a density $q_i$ with respect to $Q_i(\cdot\mid0)$ for all $i\in[n]$. Then $Q\circ\nu-Q\circ\mu$ has the following Radon--Nikodym derivative with respect to $Q(\cdot\mid0)$,
\[
\frac{d(Q\circ\nu-Q\circ\mu)}{dQ(\cdot\mid0)}(z)
\ =\ \sum_{\varnothing\neq I\subseteq[n]}\delta_{|I|}\prod_{i\in I}(q_i(z_i)-1) ,
\qquad Q(\cdot\mid0)\text{-almost everywhere}.
\]
\end{lemma}

\begin{proof}
In the formula of Lemma~\ref{lem:mixture-difference-expansion}, we can replace $\sum_{\varnothing\neq I\subseteq J}$ by $\sum_{\varnothing\neq I\subseteq[n]}$, since the summand is zero whenever $I\not\subseteq J$. 
For any measurable $A_i\subseteq\mathcal Z_i$, we have
\[
Q_i(A_i\mid1)-Q_i(A_i\mid0)
\ =\ \int_{A_i}q_i(z_i)\,Q_i(dz_i\mid0)-\int_{A_i}1\,Q_i(dz_i\mid0)
\ =\ \int_{A_i}(q_i(z_i)-1)\,Q_i(dz_i\mid0).
\]
Substituting this identity into the formula of Lemma~\ref{lem:mixture-difference-expansion}, we obtain, for any measurable  $A=\prod_{i=1}^n A_i$,
\[
\begin{aligned}
(Q\circ\nu-Q\circ\mu)(A)
&\ =\ \sum_{\varnothing\neq I\subseteq[n]}\delta_{|I|}
\prod_{i\in I}\left\{\int_{A_i}(q_i(z_i)-1)\,Q_i(dz_i\mid0)\right\}
\prod_{i\notin I}\int_{A_i}1\,Q_i(dz_i\mid0)\\
&\ =\ \sum_{\varnothing\neq I\subseteq[n]}\delta_{|I|}
\int_A\prod_{i\in I}(q_i(z_i)-1)\,Q(dz\mid0)\\
&\ =\ \int_A\sum_{\varnothing\neq I\subseteq[n]}\delta_{|I|}
\prod_{i\in I}(q_i(z_i)-1)\,Q(dz\mid0),
\end{aligned}
\]
where the second equality follows from Fubini's theorem and $Q(\cdot\mid0)=\bigotimes_{i=1}^n Q_i(\cdot\mid0)$.
Both sides define bounded signed measures, and we have shown their equality on measurable rectangles $\prod_{i=1}^n A_i$.
Since measurable rectangles form a $\pi$-system generating the product $\sigma$-algebra, this equality extends to all measurable sets.
This proves Lemma~\ref{lem:mixture-difference-density}.
\end{proof}

\begin{proof}[\textbf{Proof of Proposition~\ref{lem:mixtures-product-measures}.}]
If $Q\circ\nu=Q\circ\mu$, the inequality in Proposition~\ref{lem:mixtures-product-measures} holds trivially. Suppose henceforth that $Q\circ\nu\neq Q\circ\mu$. By Lemma~\ref{lem:mixture-difference-expansion}, there exists $1\leq\ell\leq|J|$ such that $\delta_\ell\neq0$.

Suppose first that $\chi^2\bigl(Q_j(\cdot\mid1)\|Q_j(\cdot\mid0)\bigr)=+\infty$ for some $j\in[n]$. Then $j\in J$. Let $I\subseteq J$ satisfy $|I|=\ell$ and $j\in I$. We have
\[
\delta_{|I|}^{\,2}\prod_{i\in I}\chi^2\bigl(Q_i(\cdot\mid1)\|Q_i(\cdot\mid0)\bigr)
\ =\ +\infty,
\]
since $\delta_{|I|}^2=\delta_\ell^2>0$, and all these   $\chi^2$ factors  indexed by $i\in I\subseteq J$ are strictly positive, and the factor indexed by $j$ is infinite. Thus the sum on the right-hand side of Proposition~\ref{lem:mixtures-product-measures} is infinite, so the inequality holds.

Suppose now that $\chi^2\bigl(Q_i(\cdot\mid1)\|Q_i(\cdot\mid0)\bigr)<\infty$ for each $i\in[n]$. Then $Q_i(\cdot\mid1)\ll Q_i(\cdot\mid0)$, and we write $q_i:=dQ_i(\cdot\mid1)/dQ_i(\cdot\mid0)$. By Lemma~\ref{lem:mixture-difference-density},
\[
\frac{d(Q\circ\nu-Q\circ\mu)}{dQ(\cdot\mid0)}(z)
\ =\ \sum_{\varnothing\neq I\subseteq[n]}\delta_{|I|}\prod_{i\in I}(q_i(z_i)-1),
\qquad Q(\cdot\mid0)\text{-almost everywhere}.
\]
Applying Lemma~\ref{lem:total-variation-densities} with $(P_1,P_2)=\left(Q\circ\nu,Q\circ\mu\right)$ and $R= Q(\cdot\mid0)$, we obtain
\[
\|Q\circ\nu-Q\circ\mu\|_{\mathrm{TV}}^2
\ \leq\ \frac14\int\left(\sum_{\varnothing\neq I\subseteq[n]}\delta_{|I|}
\prod_{i\in I}(q_i(z_i)-1)\right)^2Q(dz\mid0).
\]
Applying Lemma~\ref{lem:centered-density-products} with $(P_i,R_i)=\left(Q_i(\cdot\mid1),Q_i(\cdot\mid0)\right)$ and $c_I=\delta_{|I|}$, we obtain
\[
\int\left(\sum_{\varnothing\neq I\subseteq[n]}\delta_{|I|}
\prod_{i\in I}(q_i(z_i)-1)\right)^2Q(dz\mid0)
\ =\ \sum_{\varnothing\neq I\subseteq[n]}\delta_{|I|}^{\,2}
\prod_{i\in I}\chi^2\bigl(Q_i(\cdot\mid1)\|Q_i(\cdot\mid0)\bigr).
\]
Combining the last two displays proves Proposition~\ref{lem:mixtures-product-measures}.
\end{proof}

\section{Proof of Proposition~\ref{prop:entropy-laplace-upper}}
\label{app:entropy-plug-in}

In this section, we prove the risk bound for the plug-in estimator given in Proposition~\ref{prop:entropy-laplace-upper}. Recall that $T_{[a,b]}(x)=(x\wedge b)\vee a$ denotes projection onto $[a,b]$, and that $\varphi_H(t)=-t\log t$, with $\varphi_H(0)=0$.

We first recall an elementary  bound for the entropy function. Its proof is  at the end of this section.

\begin{lemma}
\label{lem:entropy-coordinate-continuity}
For any $x,y\in[0,1]$, we have $|\varphi_H(x)-\varphi_H(y)|\leq |x-y|\log\!\left(\frac{e^2}{|x-y|}\right)$, where the right-hand side is interpreted as zero when $x=y$.
\end{lemma}

\begin{proof}[Proof of Proposition~\ref{prop:entropy-laplace-upper}]
Fix $k\in[S]$. Set $\overline Z_k=n^{-1}\sum_{i\in[n]} Z_{ik}$.
Since the Laplace variables are centered and have variance $2$,
\[
\mathbb E_{\boldsymbol{p},K_{1:n}^{\mathrm{Lap},\alpha}}[\overline Z_k]=p_k,
\qquad
\operatorname{Var}_{\boldsymbol{p},K_{1:n}^{\mathrm{Lap},\alpha}}(\overline Z_k)
=
\frac{p_k(1-p_k)}{n}+\frac{8}{n\alpha^2}.
\]
Since $\widehat p_k=T_{[0,1]}(\overline Z_k)$ and $p_k\in[0,1]$, and $T_{[0,1]}$ is $1$-Lipschitz, we have $|\widehat p_k-p_k|\leq|\overline Z_k-p_k|$.
Hence 
\begin{align*}
\mathbb E_{\boldsymbol{p},K_{1:n}^{\mathrm{Lap},\alpha}}
\bigl[(\widehat p_k-p_k)^2\bigr]
\leq
\mathbb E_{\boldsymbol{p},K_{1:n}^{\mathrm{Lap},\alpha}}
\bigl[(\overline Z_k-p_k)^2\bigr]
&=
\operatorname{Var}_{\boldsymbol{p},K_{1:n}^{\mathrm{Lap},\alpha}}(\overline Z_k)
+
\left(
\mathbb E_{\boldsymbol{p},K_{1:n}^{\mathrm{Lap},\alpha}}[\overline Z_k]-p_k
\right)^2.
\end{align*}
Combining the last displays, we obtain
\begin{equation}
\label{eq:laplace-probability-quadratic-error}
\mathbb E_{\boldsymbol{p},K_{1:n}^{\mathrm{Lap},\alpha}}
\bigl[(\widehat p_k-p_k)^2\bigr]
\leq
\frac{p_k(1-p_k)}{n}+\frac{8}{n\alpha^2}
\ \leq\ \frac{9}{n\alpha^2},
\end{equation}
where the last inequality uses $\alpha\in(0,1)$.

Applying Lemma~\ref{lem:entropy-coordinate-continuity} with
$x=\widehat p_k$ and $y=p_k$ gives
\[
(\varphi_H(\widehat p_k)-\varphi_H(p_k))^2
\leq
(\widehat p_k-p_k)^2
\log^2\!\left(\frac{e^2}{|\widehat p_k-p_k|}\right) \ := \  f((\widehat p_k-p_k)^2),
\]
where
$f(u)=u\log^2(e^2/\sqrt{u})$ for $0<u\leq1$, with $f(0)=0$.
Observe that $f$ is
increasing and concave on $[0,1]$.
Then, 
Jensen's inequality gives,
for $n\alpha^2\geq9$,
\begin{align*}
\mathbb E_{\boldsymbol{p},K_{1:n}^{\mathrm{Lap},\alpha}}
\bigl[(\varphi_H(\widehat p_k)-\varphi_H(p_k))^2\bigr]
\leq
\mathbb E_{\boldsymbol{p},K_{1:n}^{\mathrm{Lap},\alpha}}
\bigl[f((\widehat p_k-p_k)^2)\bigr]
&\leq
f\left(
\mathbb E_{\boldsymbol{p},K_{1:n}^{\mathrm{Lap},\alpha}}
\bigl[(\widehat p_k-p_k)^2\bigr]
\right).
\end{align*}
Combining this bound with \eqref{eq:laplace-probability-quadratic-error} and the monotonicity of $f$, we obtain
\begin{align*}
\mathbb E_{\boldsymbol{p},K_{1:n}^{\mathrm{Lap},\alpha}}
\bigl[(\varphi_H(\widehat p_k)-\varphi_H(p_k))^2\bigr]
\leq f\left(\frac{9}{n\alpha^2}\right)
=
\frac{9}{n\alpha^2}
\log^2\!\left(\frac{e^2\sqrt{n\alpha^2}}{3}\right)
\leq
C\frac{\log^2(e+n\alpha^2)}{n\alpha^2},
\end{align*}
under the condition $n\alpha^2\geq9$, which ensures that
$9/(n\alpha^2)$ belongs to the interval $[0,1]$
on which $f$ is increasing. Here, $C$ is a universal constant.

Since
$\widehat F_{\varphi_H}=T_{[0,\log S]}(\sum_{k\in[S]}\varphi_H(\widehat p_k))$
and $F_{\varphi_H}(\boldsymbol{p})\in[0,\log S]$, and $T_{[0,\log S]}$ is $1$-Lipschitz, we have
\[
|\widehat F_{\varphi_H}-F_{\varphi_H}(\boldsymbol{p})|
\leq
\left|\sum_{k\in[S]}
\big{[}\varphi_H(\widehat p_k)-\varphi_H(p_k)\big{]}\right|.
\]
For $n\alpha^2\geq9$, using the last two bounds, and then Cauchy--Schwarz, we obtain
\begin{align*}
\mathbb E_{\boldsymbol{p},K_{1:n}^{\mathrm{Lap},\alpha}}
\bigl[(\widehat F_{\varphi_H}-F_{\varphi_H}(\boldsymbol{p}))^2\bigr]
&\leq
\mathbb E_{\boldsymbol{p},K_{1:n}^{\mathrm{Lap},\alpha}}
\left[
\left(
\sum_{k\in[S]}
\big{[}\varphi_H(\widehat p_k)-\varphi_H(p_k)\big{]}
\right)^2
\right]\\
&\leq
S\sum_{k\in[S]}
\mathbb E_{\boldsymbol{p},K_{1:n}^{\mathrm{Lap},\alpha}}
\bigl[(\varphi_H(\widehat p_k)-\varphi_H(p_k))^2\bigr]\leq
C\frac{S^2\log^2(e+n\alpha^2)}{n\alpha^2}.
\end{align*}

On the other hand, since both $\widehat F_{\varphi_H}$ and $F_{\varphi_H}(\boldsymbol{p})$ belong
to $[0,\log S]$, we also have, for any $n\geq1$,
\[
\mathbb E_{\boldsymbol{p},K_{1:n}^{\mathrm{Lap},\alpha}}
\bigl[(\widehat F_{\varphi_H}-F_{\varphi_H}(\boldsymbol{p}))^2\bigr]
\leq(\log S)^2.
\]
In particular, when $n\alpha^2<9$, this yields
\[
\mathbb E_{\boldsymbol{p},K_{1:n}^{\mathrm{Lap},\alpha}}
\bigl[(\widehat F_{\varphi_H}-F_{\varphi_H}(\boldsymbol{p}))^2\bigr]
\leq(\log S)^2
\leq S^2
\leq
9\frac{S^2\log^2(e+n\alpha^2)}{n\alpha^2}.
\]

Thus, combining the risk bounds for $n\alpha^2\geq9$ and $n\alpha^2<9$, we obtain, for all $n\geq1$ and $\alpha\in(0,1)$,
\[
\mathbb E_{\boldsymbol{p},K_{1:n}^{\mathrm{Lap},\alpha}}
\bigl[(\widehat F_{\varphi_H}-F_{\varphi_H}(\boldsymbol{p}))^2\bigr]
\leq
C'\left(
(\log S)^2
\wedge
\frac{S^2\log^2(e+n\alpha^2)}{n\alpha^2}
\right),
\]
for some universal constant $C'>0$.
Since this bound holds for every $\boldsymbol{p}\in\mathcal P_S$,
taking the supremum over $\boldsymbol{p}$ proves Proposition~\ref{prop:entropy-laplace-upper}.
\end{proof}

\begin{proof}[Proof of Lemma~\ref{lem:entropy-coordinate-continuity}]
Let $x,y\in[0,1]$. Recall that $\varphi_H(t)=-t\log t$ for $t\in(0,1]$, with $\varphi_H(0)=0$. The case $x=y$ is immediate. Otherwise, by symmetry, assume $x<y$. Since $\varphi_H'(t)=-1-\log t$ for $t>0$, we have
\[
\begin{aligned}
|\varphi_H(y)-\varphi_H(x)|
\leq \int_x^y |1+\log t|\,dt
\leq \int_0^{y-x} (1-\log t)\,dt
&= \bigl[2t-t\log t\bigr]_0^{y-x}\\
&= 2(y-x)-(y-x)\log(y-x)\\
&= (y-x)\log\!\left(\frac{e^2}{y-x}\right),
\end{aligned}
\]
where the second inequality uses $|1+\log t|\leq1-\log t$ on $(0,1]$, and the fact that $t\mapsto1-\log t$ is decreasing. 
This proves Lemma~\ref{lem:entropy-coordinate-continuity}.
\end{proof}


\section{Proof of Corollary~\ref{thm:entropy-lower}}
\label{section:proof_entropy_lb_appli_Thm1}

We prove Corollary~\ref{thm:entropy-lower} by applying Theorem~\ref{thm:general-additive-lower-bound} with $\varphi=\varphi_H$, where $\varphi_H(t)=-t\log t$. For $b\in (0, 1/2]$, we denote by $\varphi_{H,b}$ the pair profile associated with $\varphi_H$, namely
\[
\varphi_{H,b}(t)
:=
-b(1+t)\log\!\bigl(b(1+t)\bigr)
-b(1-t)\log\!\bigl(b(1-t)\bigr),
\qquad t\in[0,1].
\]
We first bound its best polynomial approximation error $E_L(\varphi_{H,b})$ from below and compute its range width $V(\varphi_{H,b})$ (Section~\ref{subsection-approx-error-entropy}). We then apply Theorem~\ref{thm:general-additive-lower-bound} to prove Corollary~\ref{thm:entropy-lower} (Section~\ref{proof-LB-entropy-coro}).


\subsection{Polynomial approximation error and range width for entropy} 
\label{subsection-approx-error-entropy}

The following lemma gives an exact expression for $V(\varphi_{H,b})$ and a lower bound on $E_L(\varphi_{H,b})$.

\begin{lemma}
\label{lem:entropy-approximation}
There exist universal constants $c_0>0$ and $L_0\geq1$ such that $V(\varphi_{H,b})=2b\log2$ and $E_L(\varphi_{H,b})\geq c_0b/L^2$, for all $b\in(0,1/2]$ and all  $L\geq L_0$.
\end{lemma}

\begin{proof}[Proof of Lemma~\ref{lem:entropy-approximation}]
Let $L\geq1$ be an integer and $b\in(0,1/2]$. Define
\[
g(t):=-(1+t)\log(1+t)-(1-t)\log(1-t),
\qquad t\in[0,1].
\]
Expanding the logarithms in the definition of $\varphi_{H,b}$ gives
\begin{align}
\varphi_{H,b}(t)
&=-b(1+t)\bigl[\log b+\log(1+t)\bigr]
-b(1-t)\bigl[\log b+\log(1-t)\bigr]\notag\\
&=-2b\log b-b\bigl[(1+t)\log(1+t)+(1-t)\log(1-t)\bigr]\notag\\
&=-2b\log b+bg(t).\label{eq:entropy-pair-profile-decomposition}
\end{align}
Moreover, $g'(t)=\log((1-t)/(1+t))\leq0$ for $t\in(0,1)$,
while $g(0)=0$ and $g(1)=-2\log2$. Therefore,
\begin{equation*}
V(\varphi_{H,b})
\ = \
\max_{t\in[0,1]}\varphi_{H,b}(t)-\min_{t\in[0,1]}\varphi_{H,b}(t)
\ = \
b\bigl(g(0)-g(1)\bigr)
\ = \
2b\log2.
\end{equation*}

By~\eqref{eq:entropy-pair-profile-decomposition}, and since adding a constant does not change the best polynomial approximation error, we have
\begin{equation}
\label{eq:entropy-hb-approximation-scaling}
E_L(\varphi_{H,b})=bE_L(g),
\end{equation}
where we used $b>0$.
We next bound $E_L(g)$ from below using the known approximation error of $\varphi_H$~\cite[Appendix~C]{WY} as follows. We have $\varphi_H(t)=g(1-t)+h(t)$, where $h(t):=(2-t)\log(2-t)$ for $t\in[0,1]$ has no singularity on this interval and is therefore easier to approximate.
Since $t\mapsto g(1-t)$ and $g$ have the same best polynomial approximation error on $[0,1]$, the triangle inequality gives
\begin{equation}
\label{eq:entropy-approximation-decomposition}
E_L(g)
\ \geq \
E_L(\varphi_H)-E_L(h).
\end{equation}
For the first term, as in~\cite[Appendix~C]{WY}, we use the classical approximation result in~\cite[Section~7.5.4]{Timan1963}, which states that $L^2E_L(\varphi_H)$ converges to a positive constant as $L\to\infty$. Consequently, there exist a constant $c_0>0$ and an integer $L_0\geq1$ such that
\begin{equation}
\label{eq:entropy-phi-approximation-lower}
E_L(\varphi_H)
\ \geq \
\frac{2c_0}{L^2},
\qquad \text{for all } L\geq L_0.
\end{equation}

To bound the second term $E_L(h)$, we use $\log(2-t)=\log2+\log(1-t/2)$ and the power series $\log(1-x)=-\sum_{k=1}^{\infty}x^k/k$, valid for $|x|<1$. We obtain
\[
\begin{aligned}
h(t)
=(2-t)\left(\log2-\sum_{k=1}^{\infty}\frac{t^k}{k2^k}\right)
&=2\log2-t\log2
-\sum_{k=1}^{\infty}\frac{t^k}{k2^{k-1}}
+\sum_{k=1}^{\infty}\frac{t^{k+1}}{k2^k}\\
&=2\log2-(1+\log2)t
+\sum_{k=2}^{\infty}
\left(\frac1{k-1}-\frac1k\right)\frac{t^k}{2^{k-1}}\\
&=2\log2-(1+\log2)t
+\sum_{k=2}^{\infty}\frac{t^k}{k(k-1)2^{k-1}},
\end{aligned}
\]
for $t\in[0,1]$. Truncating this series at degree $L\geq1$ yields
\[
E_L(h)
\leq\sum_{k=L+1}^{\infty}\frac{1}{k(k-1)2^{k-1}}
\leq\sum_{k=L+1}^{\infty}2^{-k}
=2^{-L}.
\]

Combining this last bound with \eqref{eq:entropy-approximation-decomposition} and \eqref{eq:entropy-phi-approximation-lower} gives $E_L(g)\geq 2c_0/L^2-2^{-L}$ for all $L\geq L_0$. Choosing $L_0$ sufficiently large, we also have $2^{-L}\leq c_0/L^2$ for all $L\geq L_0$. Going back to \eqref{eq:entropy-hb-approximation-scaling} yields $E_L(\varphi_{H,b})\geq c_0b/L^2$ for all $L\geq L_0$, which completes the proof of Lemma~\ref{lem:entropy-approximation}.
\end{proof}


\subsection{Proof of Corollary~\ref{thm:entropy-lower}}
\label{proof-LB-entropy-coro}

Set $m:=\lfloor(S-1)/2\rfloor$, and choose
\begin{equation}
\label{eq:entropy-lower-parameters}
L=\left\lceil\log_2(2S)\right\rceil,
\qquad
b=\frac1S
\wedge
\sqrt{\frac{L+1}{4e(e-1)^2n\alpha^2}}.
\end{equation}
Since $2m+1\leq S$ and $2mb\leq1$, the probability vector
\[
\boldsymbol p=\bigl(\underbrace{b,\ldots,b}_{2m\text{ entries}},1-2mb,\underbrace{0,\ldots,0}_{S-2m-1\text{ entries}}\bigr)
\]
belongs to $\mathcal P_{m,b}^{\mathrm{pair}}$, so this class is nonempty. These choices also satisfy the conditions in \eqref{eq:general-lb-parameters} of Theorem~\ref{thm:general-additive-lower-bound}.
Using $R_{n,\alpha}(F_{\varphi_H})\geq R_{n,\alpha}(F_{\varphi_H};\mathcal P_{m,b}^{\mathrm{pair}})$, and then applying Theorem~\ref{thm:general-additive-lower-bound} with $\varphi=\varphi_H$, we obtain
\[
R_{n,\alpha}(F_{\varphi_H})
\geq
\left[
\frac{3m^2}{32}E_L(\varphi_{H,b})^2
-\frac m4V(\varphi_{H,b})^2
\right]_+
\geq
\left[
\frac{3S^2}{512}E_L(\varphi_{H,b})^2
-\frac S8V(\varphi_{H,b})^2
\right]_+,
\]
where the second inequality follows from $S/4\leq m\leq S/2$ for $S\geq4$.
By Lemma~\ref{lem:entropy-approximation}, there exist universal constants $c_0>0$ and $L_0\geq1$ such that $V(\varphi_{H,b})=2b\log2$ and $E_L(\varphi_{H,b})\geq c_0b/L^2$ for every integer $L\geq L_0$; hence, the above risk bound yields
\[
R_{n,\alpha}(F_{\varphi_H})
\geq
Sb^2\left[
\frac{3c_0^2S}{512L^4}
-\frac12(\log2)^2
\right]_+.
\]
The choice of $L$ in \eqref{eq:entropy-lower-parameters} ensures that $L\geq L_0$ for all $S\geq S_0$, provided that the universal constant $S_0$ is chosen sufficiently large.

We now simplify the risk lower bound obtained above. By \eqref{eq:entropy-lower-parameters}, for $S\geq4$, we have
\begin{equation}
\label{eq:entropy-parameter-bounds}
L\leq\log_2(2S)+1=\log_2 S+2\leq2\log_2 S.
\end{equation}
Consequently,
\[
\frac12\left(\frac{3c_0^2S}{512L^4}\right)
\geq\frac{3c_0^2S}{1024(2\log_2 S)^4}
\geq\frac12(\log2)^2,
\]
where the last inequality holds for every $S\geq S_0$, provided that $S_0$ is chosen sufficiently large. Thus, for every $S\geq S_0$,
\[
R_{n,\alpha}(F_{\varphi_H})
\geq\frac{3c_0^2S^2b^2}{1024L^4}
=\frac{3c_0^2}{1024}
\left(
\frac1{L^4}
\wedge
\frac{S^2(L+1)}{4e(e-1)^2n\alpha^2L^4}
\right),
\]
where the equality follows from the value of $b$ in \eqref{eq:entropy-lower-parameters}. Using \eqref{eq:entropy-parameter-bounds} yields
\[
R_{n,\alpha}(F_{\varphi_H})
\geq c\left(
\frac1{(\log S)^4}
\wedge
\frac{S^2}{n\alpha^2(\log S)^3}
\right),
\]
for a universal constant $c>0$. This proves Corollary~\ref{thm:entropy-lower}.


\section{Proofs of Corollaries~\ref{cor:entropy-complementary-lower} and~\ref{cor:entropy-combined-lower}}
\label{appendix:cmplement_lower_bound_entropy}

\subsection{Proof of Corollary~\ref{cor:entropy-complementary-lower}}
\label{subsec:proof-entropy-complementary-lower}

We establish the following non-asymptotic bound. Corollary~\ref{cor:entropy-complementary-lower} then directly follows by considering the asymptotic regime $n\alpha^2\to\infty$ with $S$ fixed.

\begin{corollary}
\label{cor:entropy-complementary-nonasymptotic}
There exists a universal constant $c>0$ such that, for all $n\geq1$, $S\geq3$, and $\alpha\in(0,1)$ satisfying $n\alpha^2\geq S$, we have
\[
R_{n,\alpha}(F_{\varphi_H})
\ \geq \
c\,\frac{S}{n\alpha^2}
\log^2\!\left(\frac{e\,n\alpha^2}{S}\right).
\]
\end{corollary}

\begin{proof}[Proof of Corollary~\ref{cor:entropy-complementary-nonasymptotic}]
Set $m:=\lfloor S/2\rfloor$ and $b:=1/(2m)$. If $S$ is even, let $\boldsymbol p$ be the uniform probability vector on $[S]$. If $S$ is odd, let $\boldsymbol p$ be uniform on $[S-1]$, with $p_S=0$. In both cases, $\boldsymbol p\in\mathcal P_{m,b}^{\mathrm{pair}}$, so this class is nonempty. Since $n\alpha^2\geq S$, the parameter in Theorem~\ref{thm:orientation-lower-bound} satisfies
\begin{equation}
\label{eq:entropy-complementary-delta}
\delta_{m,b}
=1\wedge\sqrt{\frac{m}{4e(e-1)^2n\alpha^2}}
=\sqrt{\frac{m}{4e(e-1)^2n\alpha^2}}
\leq1.
\end{equation}

Since $R_{n,\alpha}(F_{\varphi_H})\geq R_{n,\alpha}(F_{\varphi_H};\mathcal P_{m,b}^{\mathrm{pair}})$, Theorem~\ref{thm:orientation-lower-bound} applied with $\varphi=\varphi_H$ gives
\begin{equation}
\label{eq:entropy-complementary-risk-from-theorem}
R_{n,\alpha}(F_{\varphi_H})
\geq \frac{m^2}{16}
\sup_{\substack{t,t'\in[0,1]\\ |t-t'|\leq\delta_{m,b}}}
\left|\varphi_{H,b}(t)-\varphi_{H,b}(t')\right|^2
\geq \frac{m^2}{16}
\left|\varphi_{H,b}(1)-\varphi_{H,b}(1-\delta_{m,b})\right|^2,
\end{equation}
where the second inequality follows by taking $t=1$ and $t'=1-\delta_{m,b}$, which belong to $[0,1]$ and satisfy $|t-t'|=\delta_{m,b}$.

We now bound \eqref{eq:entropy-complementary-risk-from-theorem} from below. Write $\varphi_{H,b}$ for the pair profile associated with the entropy function $\varphi_H(x)=-x\log x$, namely
\[
\varphi_{H,b}(t)
:=\varphi_H(b(1+t))+\varphi_H(b(1-t))
=-b(1+t)\log\!\bigl(b(1+t)\bigr)-b(1-t)\log\!\bigl(b(1-t)\bigr).
\]
Since $b=1/(2m)$, we have $m\varphi_{H,b}(1)=\log m$, and
\[
\begin{aligned}
m\varphi_{H,b}(1-\delta_{m,b})
&=-\left(1-\frac{\delta_{m,b}}2\right)
\log\!\left(\frac{1-\delta_{m,b}/2}{m}\right)
-\frac{\delta_{m,b}}2\log\!\left(\frac{\delta_{m,b}}{2m}\right)\\
&=\log m
-\left(1-\frac{\delta_{m,b}}2\right)\log\!\left(1-\frac{\delta_{m,b}}2\right)
+\frac{\delta_{m,b}}2\log\!\left(\frac{2}{\delta_{m,b}}\right).
\end{aligned}
\]
Consequently,
\[
m\bigl(\varphi_{H,b}(1-\delta_{m,b})-\varphi_{H,b}(1)\bigr)
\geq\frac{\delta_{m,b}}2\log\!\left(\frac{2}{\delta_{m,b}}\right),
\]
since $\left(1-\delta_{m,b}/2\right)\log\left(1-\delta_{m,b}/2\right)\leq0$.
Substituting this bound into \eqref{eq:entropy-complementary-risk-from-theorem} gives
\begin{equation}
\label{eq:entropy-complementary-separation-bound}
R_{n,\alpha}(F_{\varphi_H})
\geq\frac{\delta_{m,b}^2}{64}
\log^2\!\left(\frac{2}{\delta_{m,b}}\right).
\end{equation}

By the definition of $\delta_{m,b}$ in \eqref{eq:entropy-complementary-delta}, we have
\[
\log\!\left(\frac{2}{\delta_{m,b}}\right)
=\frac12\log\!\left(\frac{16e(e-1)^2n\alpha^2}{m}\right)
\geq\frac12\log\!\left(\frac{e\,n\alpha^2}{m}\right).
\]
Substituting this inequality and the value of $\delta_{m,b}^2$ from \eqref{eq:entropy-complementary-delta} into \eqref{eq:entropy-complementary-separation-bound} gives
\[
R_{n,\alpha}(F_{\varphi_H})
\geq c\,\frac{m}{n\alpha^2}
\log^2\!\left(\frac{e\,n\alpha^2}{m}\right),
\]
for a universal constant $c>0$. Using $S/3\leq m\leq S$, we obtain the lower bound in Corollary~\ref{cor:entropy-complementary-nonasymptotic}.
\end{proof}

\subsection{Proof of Corollary~\ref{cor:entropy-combined-lower}}
\label{subsec:proof-entropy-combined-lower}

Since $n\alpha^2\geq S^2\log S$, Corollary~\ref{thm:entropy-lower} gives
\[
R_{n,\alpha}(F_{\varphi_H})
\geq c_1\frac{S^2}{n\alpha^2(\log S)^3}.
\]
Moreover, Corollary~\ref{cor:entropy-complementary-nonasymptotic} gives
\[
R_{n,\alpha}(F_{\varphi_H})
\geq c_2\frac{S}{n\alpha^2}
\log^2\!\left(\frac{e\,n\alpha^2}{S}\right),
\]
where $c_1,c_2>0$ are universal constants. Combining these two bounds using $\max\{a,b\}\geq(a+b)/2$ proves Corollary~\ref{cor:entropy-combined-lower} with $c=\min\{c_1,c_2\}/2$.

\section{Proof of Corollary~\ref{cor:power-bias-lower-combined}}
\label{section:proof_power_sum_general_lb}
We prove Corollary~\ref{cor:power-bias-lower-combined} by establishing the following two lower bounds, with all log factors written explicitly. The first covers $\gamma\in(0,1)$, and the second covers $\gamma\in(1,2)$.

\begin{corollary}
\label{cor:power-sum-lower-bound}
Let $\gamma\in(0,1)$.
There exist constants $c_\gamma>0$ and $S_\gamma\geq4$,
depending only on $\gamma$, such that, for all
$n\geq1$, $S\geq S_\gamma$, and $\alpha\in(0,1)$, we have
\begin{equation*}
R_{n,\alpha}(F_{\varphi_\gamma})
\geq c_\gamma
\left(
\frac{S^2}{(\alpha^2 n)^\gamma(\log S)^{3\gamma}}
\wedge
\frac{S^{2(1-\gamma)}}{(\log S)^{4\gamma}}
\right).
\end{equation*}
\end{corollary}

\begin{corollary}
\label{cor:power-sum-lower-superlinear}
Let $\gamma\in(1,2)$.
There exist constants $c_\gamma>0$ and $S_\gamma\geq4$,
depending only on $\gamma$, such that, for all
$n\geq1$, $S\geq S_\gamma$, and $\alpha\in(0,1)$ satisfying
$\alpha^2 n\geq S_\gamma^2$, we have
\begin{equation*}
R_{n,\alpha}(F_{\varphi_\gamma})
\geq
\frac{c_\gamma}
{\log^{4\gamma}\!\bigl(S\wedge\sqrt{\alpha^2 n}\bigr)}
\left(
\frac{S^2}{(\alpha^2 n)^\gamma}
\wedge
\frac{1}{(\alpha^2 n)^{\gamma-1}}
\right).
\end{equation*}
\end{corollary}

The proofs of Corollaries~\ref{cor:power-sum-lower-bound} and~\ref{cor:power-sum-lower-superlinear} follow the same structure as that of Corollary~\ref{thm:entropy-lower}. We apply Theorem~\ref{thm:general-additive-lower-bound} with $\varphi=\varphi_\gamma$, where $\varphi_\gamma(t)=t^\gamma$ and $\gamma\in(0,1)\cup(1,2)$. For $b\in(0,1/2]$, we denote by $\varphi_{\gamma,b}$ the pair profile associated with $\varphi_\gamma$, namely
\[
\varphi_{\gamma,b}(t)
:=\bigl(b(1+t)\bigr)^\gamma+\bigl(b(1-t)\bigr)^\gamma,
\qquad t\in[0,1].
\]
We first bound its best polynomial approximation error $E_L(\varphi_{\gamma,b})$ from below and compute its range width $V(\varphi_{\gamma,b})$ (Section~\ref{subsection-approx-error-power-sum}). We then apply Theorem~\ref{thm:general-additive-lower-bound} to obtain a lower bound on $R_{n,\alpha}(F_{\varphi_\gamma})$ in terms of the parameters $(L,b,m)$ (Section~\ref{appli-thm-lb-powersum}). We conclude by choosing the values of these parameters to prove Corollaries~\ref{cor:power-sum-lower-bound} and~\ref{cor:power-sum-lower-superlinear} (Section~\ref{proof-LB-power-sum-coro}).


\subsection{Approximation error and range width for power sums}
\label{subsection-approx-error-power-sum}

The following lemma gives an exact expression for $V(\varphi_{\gamma,b})$ and a lower bound on $E_L(\varphi_{\gamma,b})$.

\begin{lemma}
\label{lem:power-sum-approximation}
There exist constants $c_0>0$ and $L_\gamma\geq1$, depending only on $\gamma$, such that $V(\varphi_{\gamma,b})=b^\gamma|2^\gamma-2|$ and $E_L(\varphi_{\gamma,b})\geq c_0b^\gamma/L^{2\gamma}$, for all $b\in(0,1/2]$ and all integers $L\geq L_\gamma$.
\end{lemma}

\begin{proof}[Proof of Lemma~\ref{lem:power-sum-approximation}]
Let $L\geq1$ be an integer and $b\in(0,1/2]$. Define $g_\gamma(t):=(1+t)^\gamma+(1-t)^\gamma$ for $t\in[0,1]$. For $t\in(0,1)$, we have
$g_\gamma'(t)=\gamma\bigl((1+t)^{\gamma-1}-(1-t)^{\gamma-1}\bigr)$.
Since $1-t<1+t$, the derivative is negative for $\gamma\in(0,1)$ and positive for $\gamma\in(1,2)$. Thus, $g_\gamma$ is monotone on $[0,1]$, with $g_\gamma(0)=2$ and $g_\gamma(1)=2^\gamma$. Using $\varphi_{\gamma,b}=b^\gamma g_\gamma$, we obtain
\[
V(\varphi_{\gamma,b})
=b^\gamma\lvert g_\gamma(1)-g_\gamma(0)\rvert
=b^\gamma\lvert 2^\gamma-2\rvert.
\]

Since $\varphi_{\gamma,b}=b^\gamma g_\gamma$ and $b>0$, we have
\begin{equation}
\label{eq:power-sum-hb-approximation-scaling}
E_L(\varphi_{\gamma,b})
\ = \
b^\gamma E_L(g_\gamma).
\end{equation}
We next bound $E_L(g_\gamma)$ from below using the known approximation error of $\varphi_\gamma$~\cite{JiaoVenkatHanWeissman2015} as follows. We have $\varphi_\gamma(t)=g_\gamma(1-t)-h(t)$, where $h(t):=(2-t)^\gamma$ for $t\in[0,1]$ has no singularity on this interval and is therefore easier to approximate.
Since $t\mapsto g_\gamma(1-t)$ and $g_\gamma$ have the same best polynomial approximation error on $[0,1]$, the triangle inequality gives
\begin{equation}
\label{eq:power-sum-approximation-decomposition}
E_L(g_\gamma)
\ \geq \
E_L(\varphi_\gamma)-E_L(h).
\end{equation}
For the first term, we use \cite[Lemma~17]{JiaoVenkatHanWeissman2015},
which states that $L^{2\gamma}E_L(\varphi_\gamma)$ converges to a positive
constant as \(L\to\infty\).
Consequently, there exist a constant \(c_0>0\) and
an integer \(L_\gamma\geq1\), depending only on \(\gamma\),
such that
\begin{equation}
\label{eq:power-sum-phi-approximation-lower}
E_L(\varphi_\gamma)
\ \geq \
\frac{2c_0}{L^{2\gamma}},
\qquad \text{for all } L\geq L_\gamma.
\end{equation}

To bound the second term $E_L(h)$, we use the binomial series $(1-x)^\gamma=\sum_{j=0}^{\infty}(-1)^j\binom{\gamma}{j}x^j$, valid for $|x|<1$. This gives
\begin{equation*}
h(t)=(2-t)^\gamma
\ = \
2^\gamma(1-t/2)^\gamma
\ = \
2^\gamma\sum_{j=0}^{\infty}
(-1)^j\binom{\gamma}{j}\left(\frac t2\right)^j,
\qquad t\in[0,1],
\end{equation*}
where $\binom{\gamma}{0}=1$, and $\binom{\gamma}{j}=\gamma(\gamma-1)\cdots(\gamma-j+1)/j!$ for $j\geq1$.
Observe that \(\left|\binom{\gamma}{j}\right|\leq2\)
for all \(j\geq0\); a proof of this elementary inequality is given at the end of this subsection.
Truncating the above series expansion
at degree \(L\geq1\), and then using
\(\left|\binom{\gamma}{j}\right|\leq2\), yields
\[
E_L(h)
\ \leq \
2^\gamma\sum_{j=L+1}^{\infty}
\left|\binom{\gamma}{j}\right|2^{-j}
\ \leq \
2^{\gamma+1}\sum_{j=L+1}^{\infty}2^{-j}
\ = \
2^{\gamma+1}2^{-L}.
\]
Combining this bound with
\eqref{eq:power-sum-approximation-decomposition} and
\eqref{eq:power-sum-phi-approximation-lower} gives
$E_L(g_\gamma)
 \geq 
(2c_0/L^{2\gamma})-2^{\gamma+1-L}$.
Taking \(L_\gamma\) sufficiently large,
we get
\(2^{\gamma+1-L}\leq c_0/L^{2\gamma}\)
for every \(L\geq L_\gamma\).
Going back to
\eqref{eq:power-sum-hb-approximation-scaling} yields
$E_L(\varphi_{\gamma,b})\geq c_0b^\gamma/L^{2\gamma}$
for \(L\geq L_\gamma\),
which completes the proof of
Lemma~\ref{lem:power-sum-approximation}.
\end{proof}

\begin{proof}[Proof of \(\left|\binom{\gamma}{j}\right|\leq2\)
for all \(j\geq0\)]
The first three coefficients satisfy
\[
\left|\binom{\gamma}{0}\right|=1,
\qquad
\left|\binom{\gamma}{1}\right|=\gamma<2,
\qquad
\left|\binom{\gamma}{2}\right|
=\frac{\gamma|\gamma-1|}{2}<1,
\]
while for every \(j\geq3\),
\[
\left|\binom{\gamma}{j}\right|
=
\frac{j-1-\gamma}{j}
\left|\binom{\gamma}{j-1}\right|
\leq
\left|\binom{\gamma}{j-1}\right|.
\]
So we conclude by induction that $\left|\binom{\gamma}{j}\right|\leq2$
for all \(j\geq0\).
\end{proof}

\subsection{Application of Theorem~\ref{thm:general-additive-lower-bound}}
\label{appli-thm-lb-powersum}

Let $m$ be an integer such that $1\leq m\leq\lfloor(S-1)/2\rfloor$, to be chosen in Section~\ref{proof-LB-power-sum-coro}. Set
\begin{equation}
\label{eq:power-sum-lower-parameters}
L=\left\lceil\log_2(2m)\right\rceil,
\qquad \qquad
b=\frac{1}{2m+1}
\wedge
\sqrt{\frac{L+1}{4e(e-1)^2\alpha^2 n}}.
\end{equation}
Since $2m+1\leq S$ and $2mb\leq1$, the probability vector
\[
\boldsymbol p=\bigl(\underbrace{b,\ldots,b}_{2m\text{ entries}},1-2mb,\underbrace{0,\ldots,0}_{S-2m-1\text{ entries}}\bigr)
\]
belongs to $\mathcal P_{m,b}^{\mathrm{pair}}$, so this class is nonempty. The choices of $L$ and $b$ also satisfy the parameter conditions \eqref{eq:general-lb-parameters} of Theorem~\ref{thm:general-additive-lower-bound}.
Using $R_{n,\alpha}(F_{\varphi_\gamma})\geq R_{n,\alpha}(F_{\varphi_\gamma};\mathcal P_{m,b}^{\mathrm{pair}})$, and then applying Theorem~\ref{thm:general-additive-lower-bound} with $\varphi=\varphi_\gamma$, we obtain
\[
R_{n,\alpha}(F_{\varphi_\gamma})
\ \geq \
\left[
\frac{3m^2}{32}E_L(\varphi_{\gamma,b})^2
-\frac m4 V(\varphi_{\gamma,b})^2
\right]_+.
\]
By Lemma~\ref{lem:power-sum-approximation}, we have $V(\varphi_{\gamma,b})=b^\gamma|2^\gamma-2|$ and $E_L(\varphi_{\gamma,b})\geq c_0b^\gamma/L^{2\gamma}$ for $L\geq L_\gamma$. Hence
\[
R_{n,\alpha}(F_{\varphi_\gamma})
\geq
\frac{mb^{2\gamma}}4
\left[
\frac{3c_0^2m}{8L^{4\gamma}}
-(2^\gamma-2)^2
\right]_+.
\]
The choice of $L$ in \eqref{eq:power-sum-lower-parameters} ensures that $L\geq L_\gamma$ for all $m\geq m_\gamma$, provided that the  constant $m_\gamma$ (which depends only on $\gamma$) is chosen sufficiently large.

We now simplify the risk lower bound obtained above. By \eqref{eq:power-sum-lower-parameters}, for $m\geq4$, we have
\begin{equation}
\label{eq:power-sum-parameter-bounds}
L\leq\log_2(2m)+1=\log_2 m+2\leq2\log_2 m.
\end{equation}
Consequently,
\[
\frac{3c_0^2m}{16L^{4\gamma}}
\geq
\frac{3c_0^2m}{16(2\log_2 m)^{4\gamma}}
\geq
(2^\gamma-2)^2,
\]
where the last inequality holds for all $m\geq m_\gamma$ after increasing $m_\gamma$ if necessary. Thus
\[
R_{n,\alpha}(F_{\varphi_\gamma})
\geq \frac{3c_0^2m^2b^{2\gamma}}{64L^{4\gamma}}
=
\frac{3c_0^2}{64}
\left(
\frac{m^2}{(2m+1)^{2\gamma}L^{4\gamma}}
\wedge
\frac{m^2(L+1)^\gamma}
{[4e(e-1)^2]^\gamma(\alpha^2 n)^\gamma L^{4\gamma}}
\right),
\]
where the equality follows from the choice of $b$ in \eqref{eq:power-sum-lower-parameters}. Since $2m+1\leq3m$ and $L+1\geq L$, this gives
\begin{equation}
\label{eq:power-sum-lower-k}
R_{n,\alpha}(F_{\varphi_\gamma})
\geq
c_\gamma
\left(
\frac{m^{2-2\gamma}}{L^{4\gamma}}
\wedge
\frac{m^2}{(\alpha^2 n)^\gamma L^{3\gamma}}
\right),
\end{equation}
where the constant $c_\gamma>0$ depends only on $\gamma$.


\subsection{Proofs of Corollaries~\ref{cor:power-sum-lower-bound} and~\ref{cor:power-sum-lower-superlinear}}
\label{proof-LB-power-sum-coro}

\textit{Case $\gamma\in(0,1)$.}
Take $m=\lfloor(S-1)/2\rfloor$. For $S\geq4$, this gives $m\geq(S-2)/2\geq S/4$. Using this inequality and $L\leq2\log_2 S$ from \eqref{eq:power-sum-parameter-bounds} in the risk lower bound \eqref{eq:power-sum-lower-k}, we obtain
\[
R_{n,\alpha}(F_{\varphi_\gamma})
\geq c_\gamma'\left(
\frac{S^{2(1-\gamma)}}{(\log S)^{4\gamma}}
\wedge
\frac{S^2}{(\alpha^2 n)^\gamma(\log S)^{3\gamma}}
\right)
\]
for $S\geq S_\gamma$, with the constant $S_\gamma$ chosen sufficiently large, and for some constant $c_\gamma'>0$ depending only on $\gamma$. This proves Corollary~\ref{cor:power-sum-lower-bound}.

\smallskip

\noindent \textit{Case $\gamma\in(1,2)$.}
Take
\begin{equation}
\label{eq:power-sum-superlinear-support}
m=\left\lfloor\frac{(S\wedge\sqrt{\alpha^2 n})-1}{2}\right\rfloor.
\end{equation}
For $S\geq S_\gamma$ and $\alpha^2 n\geq S_\gamma^2$, this choice gives $m\geq m_\gamma$ once $S_\gamma$ is chosen sufficiently large. Applying \eqref{eq:power-sum-lower-k} and using $L\geq1$, we obtain
\[
R_{n,\alpha}(F_{\varphi_\gamma})
\geq
\frac{c_\gamma}{L^{4\gamma}}
\left(
m^{2-2\gamma}\wedge\frac{m^2}{(\alpha^2 n)^\gamma}
\right) \ = \ \frac{c_\gamma}{L^{4\gamma}}
\frac{m^2}{(\alpha^2 n)^\gamma},
\]
where the minimum equals its second term since $m^2\leq\alpha^2 n$.
By the choice of $m$ in \eqref{eq:power-sum-superlinear-support}, we have $m\geq\frac14(S\wedge\sqrt{\alpha^2 n})$, so
\[
\frac{m^2}{(\alpha^2 n)^\gamma}
\geq
\frac1{16}
\left(
\frac{S^2}{(\alpha^2 n)^\gamma}
\wedge
\frac1{(\alpha^2 n)^{\gamma-1}}
\right).
\]
Also, \eqref{eq:power-sum-parameter-bounds} gives
$L\leq2\log_2 m
\leq\frac2{\log2}\log\bigl(S\wedge\sqrt{\alpha^2 n}\bigr)$.
Combining the last three bounds yields
\[
R_{n,\alpha}(F_{\varphi_\gamma})
\geq
\frac{c_\gamma'}
{\log^{4\gamma}\!\bigl(S\wedge\sqrt{\alpha^2 n}\bigr)}
\left(
\frac{S^2}{(\alpha^2 n)^\gamma}
\wedge
\frac1{(\alpha^2 n)^{\gamma-1}}
\right),
\]
for some constant $c_\gamma'>0$ depending only on $\gamma$.
This proves Corollary~\ref{cor:power-sum-lower-superlinear} and completes the proofs of both corollaries.


\section{Proof of Corollary~\ref{cor:power-variance-lower}}
\label{appendix:power-variance-lower}

We first prove the following lower bound for $\gamma\in(1,3/2)$. Together with the lower bound of order $1/(\alpha^2 n)$ from~\cite{ButuceaIssartel2021}, it yields Corollary~\ref{cor:power-variance-lower}.

\begin{corollary}
\label{lem:power-sum-complementary-lower}
Let $\gamma\in(1,3/2)$.
There exists a constant $c_\gamma>0$,
depending only on $\gamma$, such that, for all
$S\geq3$, $n\geq1$, and $\alpha\in(0,1)$
satisfying $\alpha^2 n>T_\gamma(S)$, we have $R_{n,\alpha}(F_{\varphi_\gamma})\geq c_\gamma S^{3-2\gamma}/(\alpha^2 n)$.
\end{corollary}

\begin{proof}[Proof of Corollary~\ref{lem:power-sum-complementary-lower}]
Set $m:=\lfloor S/2\rfloor$ and $b:=1/(2m)$. If $S$ is even, let $\boldsymbol p$ be the uniform probability vector on $[S]$. If $S$ is odd, let $\boldsymbol p$ be uniform on $[S-1]$, with $p_S=0$. In both cases, $\boldsymbol p\in\mathcal P_{m,b}^{\mathrm{pair}}$, so this class is nonempty. Since $\alpha^2 n>T_\gamma(S)\geq S$, the parameter $\delta_{m,b}$ in Theorem~\ref{thm:orientation-lower-bound} satisfies
\begin{equation}
\label{eq:power-sum-complementary-delta}
\delta_{m,b}
=1\wedge\sqrt{\frac{m}{4e(e-1)^2\alpha^2 n}}
=\sqrt{\frac{m}{4e(e-1)^2\alpha^2 n}}
\leq\frac{1}{\sqrt{8e}(e-1)}
\leq\frac14.
\end{equation}
Since $R_{n,\alpha}(F_{\varphi_\gamma})\geq R_{n,\alpha}(F_{\varphi_\gamma};\mathcal P_{m,b}^{\mathrm{pair}})$, Theorem~\ref{thm:orientation-lower-bound} applied with $\varphi=\varphi_\gamma$ gives
\begin{equation}
\label{eq:power-complementary-risk-from-theorem}
R_{n,\alpha}(F_{\varphi_\gamma})
\geq\frac{m^2}{16}
\sup_{\substack{t,t'\in[0,1]\\|t-t'|\leq\delta_{m,b}}}
\left|\varphi_{\gamma,b}(t)-\varphi_{\gamma,b}(t')\right|^2
\geq\frac{m^2}{16}
\left|\varphi_{\gamma,b}\!\left(\frac12+\delta_{m,b}\right)
-\varphi_{\gamma,b}\!\left(\frac12\right)\right|^2,
\end{equation}
where the second inequality follows by taking  $t=1/2+\delta_{m,b}$ and $t'=1/2$, which belong to $[0,1]$ and satisfy $|t-t'|=\delta_{m,b}$.

We now bound the functional separation. Write $\varphi_{\gamma,b}$ for the pair profile associated with $\varphi_\gamma(x)=x^\gamma$, namely
\[
\varphi_{\gamma,b}(t)
:=\varphi_\gamma(b(1+t))+\varphi_\gamma(b(1-t))
=b^\gamma\bigl((1+t)^\gamma+(1-t)^\gamma\bigr).
\]
Since $b=1/(2m)$, we have
\[
m\varphi_{\gamma,b}'(t)
=
\frac{\gamma m^{1-\gamma}}{2^\gamma}
\bigl((1+t)^{\gamma-1}-(1-t)^{\gamma-1}\bigr).
\]
For $t\in[1/2,3/4]$, this gives
\[
m\varphi_{\gamma,b}'(t)
\geq d_\gamma m^{1-\gamma},
\qquad \quad 
d_\gamma
:=
\frac{\gamma}{2^\gamma}
\left[
\left(\frac32\right)^{\gamma-1}
-
\left(\frac12\right)^{\gamma-1}
\right]
>0.
\]
Consequently, using $\delta_{m,b}\leq1/4$ from \eqref{eq:power-sum-complementary-delta}, we obtain
\[
m\left(\varphi_{\gamma,b}\!\left(\frac12+\delta_{m,b}\right)-\varphi_{\gamma,b}\!\left(\frac12\right)\right)
=\int_{1/2}^{1/2+\delta_{m,b}}m\varphi_{\gamma,b}'(t)\,\mathrm dt
\geq d_\gamma\delta_{m,b}m^{1-\gamma}.
\]
Substituting this bound into \eqref{eq:power-complementary-risk-from-theorem}, and then using the value of $\delta_{m,b}$ in \eqref{eq:power-sum-complementary-delta}, we obtain
\[
R_{n,\alpha}(F_{\varphi_\gamma})
\geq \frac{d_\gamma^2}{16}\delta_{m,b}^2m^{2-2\gamma}
=\frac{d_\gamma^2}{64e(e-1)^2}\frac{m^{3-2\gamma}}{\alpha^2 n}.
\]
Finally, using $m=\lfloor S/2\rfloor\geq S/3$ and $3-2\gamma>0$, we obtain the lower bound in Corollary~\ref{lem:power-sum-complementary-lower}.
\end{proof}

\end{document}